\documentclass[12pt]{amsart}
\usepackage{epsfig}
\newfont{\sheaf}{eusm10 scaled\magstep1}

\usepackage{graphics}
\usepackage{color}

\newcommand{\ra}{\ensuremath{\rightarrow}}

\def\eea{\end{eqnarray*}}
\def\bea{\begin{eqnarray*}}
\def\Bbb{\bf}

\def\CC{{\Bbb C}}

\def\X{{\mathcal{X}}}

\def\K{{\mathcal{K}}}
\def\MM{{\mathcal{M}}}
\def\de{{\delta}}
\def\De{{\Delta}}
\def\ga{{\gamma}}
\def\Ga{{\Gamma}}

\def\NN{{\Bbb N}}

\newcommand{\Proof}{{\it Proof. }}
\newcommand{\QED}{{\hfill $Q.E.D.$}}

\newtheorem{teo}{Theorem}[section]
\newtheorem{df}[teo]{Definition}
\newtheorem{lem}[teo]{Lemma}
\newtheorem{cor}[teo]{Corollary}
\newtheorem{ex}[teo]{Example}
\newtheorem{oss}[teo]{Remark}
\newtheorem{prop}[teo]{Proposition}
\newtheorem{conj}[teo]{Conjecture}
\newtheorem{question}[teo]{Question}

\newcommand{\U}{\ensuremath{\mathbb{U}}}
\newcommand{\C}{\ensuremath{\mathbb{C}}}
\newcommand{\R}{\ensuremath{\mathbb{R}}}
\newcommand{\Z}{\ensuremath{\mathbb{Z}}}
\newcommand{\Q}{\ensuremath{\mathbb{Q}}}
\newcommand{\F}{\ensuremath{\mathbb{F}}}
\newcommand{\M}{\ensuremath{\mathbb{M}}}
\newcommand{\N}{\ensuremath{\mathbb{N}}}
\newcommand{\hol}{\ensuremath{\mathcal{O}}}
\newcommand{\Hol}{ \ensuremath{\mathcal{O}}^h}
\newcommand{\HH}{\ensuremath{\mathbb{H}}}
\newcommand{\LL}{\ensuremath{\mathbb{L}}}
\newcommand{\PP}{\ensuremath{\mathbb{P}}}
\newcommand{\RRR}{\ensuremath{\mathcal{R}}}
\newcommand{\BB}{\ensuremath{\mathcal{B}}}
\newcommand{\FFF}{\ensuremath{\mathcal{F}}}
\newcommand{\A}{\ensuremath{\mathcal{A}}}

\newcommand{\cD}{\ensuremath{\mathcal{D}}}
\newcommand{\CCC}{\ensuremath{\mathcal{C}}}

\newcommand{\HHH}{\ensuremath{\mathcal{H}}}

\newcommand{\I}{\ensuremath{\mathcal{I}}}
\newcommand{\SSS}{\ensuremath{\mathcal{S}}}
\newcommand{\PPP}{\ensuremath{\mathcal{P}}}
\newcommand{\VV}{\ensuremath{\mathbb{V}}}

\newcommand{\sT}{{\mathcal T}}

\newcommand{\sG}{{\mathcal G}}

\newcommand{\sQ}{{\mathcal Q}}
\newcommand{\sM}{{\mathcal M}}
\newcommand{\sK}{{\mathcal K}}
\newcommand{\sV}{{\mathcal V}}
\newcommand{\sB}{{\mathcal B}}

\newcommand{\la}{\lambda}

\newcommand{\Om}{\Omega}
\newcommand{\La}{\Lambda}
\newcommand{\e}{\epsilon}

\begin{document}

\title[Differentiable and deformation type of algebraic surfaces]{  Differentiable and
deformation type of algebraic surfaces, real and symplectic structures.}

%\author{Fabrizio Catanese\\
%  Universit\"at Bayreuth \footnote{
%The research of the  author was supported by the
%  SCHWERPUNKT "Globale Methoden in der komplexen Geometrie",
%a visit to IHES  by contract nr RITA-CT-2004-505493.
%}}

\author{Fabrizio Catanese\\
  Universit\"at Bayreuth}
\thanks{
The research of the  author was supported by the
  SCHWERPUNKT "Globale Methoden in der komplexen Geometrie",
a visit to IHES  by contract nr RITA-CT-2004-505493.
}

\date{February 27, 2007}
\maketitle
\tableofcontents

\bigskip

\

\vfill

\begin{minipage}{0.8\linewidth}
\small
Che differenza c'e' fra il palo, il Paolo e la banana?

\hfill {\it (G. Lastrucci, Firenze, 9/11/62)}
\end{minipage}
\bigskip

\newpage

\section*{INTRODUCTION}

  As already announced in Cetraro at the beginning of the C.I.M.E. course,
we deflected from the   broader target `` Classification
and deformation types
of complex and real manifolds'', planned and announced originally.

  First of all,  the lectures actually delivered focused on the intersection of
the above vast area with the theme of the School,
``Algebraic surfaces and symplectic 4-manifolds''.

Hence the title of the Lecture Notes has been changed accordingly.

Moreover, the Enriques classification  of real algebraic surfaces
is not touched upon here, and complex conjugation and real structures
appear mostly through their relation to
deformation types of complex manifolds, and in particular
through their relation with strong and weak rigidity theorems.

In some sense then this course is a continuation of the C.I.M.E. course
I held some 20 years ago in Montecatini (\cite{cat6}),
about ` Moduli of algebraic surfaces '.

But whereas  those Lecture Notes had an  initial part of considerable 
length which was
meant to be a general introduction to complex deformation theory, here
the main results of deformation theory which we need are only
stated.

Nevertheless, because the topic can be of interest not only to 
algebraic geometers,
but also to people working in differential or symplectic topology,
we decided to start dedicating the first lecture to recalling
basic notions concerning projective and K\"ahler
manifolds. Especially, we recall the main principles of 
classification theory, and state the Enriques
classification of algebraic surfaces of special type.

Since surfaces of general type and their moduli spaces are a major 
theme for us here,
it seemed worthwhile to recall in detail in lecture two the structure 
of their canonical models,
in particular of their singularities, the socalled Rational Double Points, or
Kleinian quotient singularities. The rest of lecture two is devoted to proving
Bombieri's  theorem on pluricanonical embeddings, to the analysis of 
other quotient
singularities, and to the deformation equivalence relation (showing that
two minimal models are deformation equivalent iff the respective 
canonical models are).
Bombieri's theorem is proven in every detail for the case of an ample 
canonical divisor,
with the hope that some similar result may soon be proven  also in 
the symplectic case.

In lecture three we  show first that deformation equivalence
implies diffeomorphism, and then, using a result concerning symplectic
approximations of projective varieties with isolated singularities 
and Moser's theorem,
we show that a surfaces of general type has a 'canonical symplectic structure',
i.e., a symplectic structure whose class is the class of the canonical divisor,
and which is unique up to symplectomorphism.

In lecture three and the following ones we  thus enter ' in medias res',
since one of the main problems that we discuss in these Lecture Notes is
the comparison of
differentiable and deformation type of minimal surfaces of general type,
keeping also
in consideration the  canonical
symplectic structure (unique up to symplectomorphism and
   invariant for smooth deformation) which these surfaces possess.

We present several counterexamples  to the
  DEF = DIFF speculation of Friedman and Morgan (\cite{f-m1}) that deformation
type and  diffeomorphism type  should coincide for complex algebraic 
surfaces. The
first ones were obtained
by Manetti (\cite{man4}), and exhibit non simply connected surfaces which are
pairwise  not deformation equivalent. We were later  able to show
that they are canonically symplectomorphic (see \cite{cat02} and also 
\cite{cat06}).
An account of these
results is to be found in chapter 6, which is an extra chapter
with title `Epilogue' (we hope however that this title
may soon turn out to be inappropriate in view of future further
developments) .

In lecture 4, after discussing  some classical results (like the 
theorem of Castelnuovo
and De Franchis) and some `semi-classical'
results (by the author) concerning the topological characterization
of irrational pencils on K\"ahler manifolds and algebraic surfaces,
we discuss orbifold fundamental groups and triangle covers.

We use the above results to describe varieties isogenous to a product.
These yield
several examples of surfaces not deformation equivalent to
their complex conjugate surface. We  describe in particular the
examples  by the present author (\cite{cat4}), by
Bauer-Catanese-Grunewald (\cite{bcg}), and then the ones by Kharlamov-Kulikov
(\cite{k-k}) which yield ball quotients.
  In this lecture we discuss complex conjugation
and real structures, starting from elementary examples and ending 
with a survey of
recent results and with open problems on  the theory of `Beauville surfaces'.

The beginning of lecture 5 is again rather elementary, it discusses
connected sums and other surgeries, like fibre sums, and recalls basic
definitions and results on braid groups, mapping class groups and
Hurwitz equivalence.

  After recalling the theory of Lefschetz pencils, especially the
differentiable viewpoint introduced by Kas (\cite{kas}), we
recall Freedman's basic results on the topology
of simply connected compact (oriented) fourmanifolds (see \cite{f-q}).

  We finally devote ourselves to  our main objects of investigation, namely,
  the socalled `(abc)-surfaces' (introduced in \cite{cat02}),
which are simply connected.
We finish Lecture 5 explaining our joint work with Wajnryb 
(\cite{c-w}) dedicated to
the proof that  these last surfaces
are diffeomorphic to each other when the two integers $b$ and $ a+c$ are fixed.

In Chapter 6 we sketch
the proof that these, under suitable numerical conditions, are not
deformation equivalent. A  result which is only very slightly weaker
is explained in the Lecture Notes by Manetti, but with  many more details;
needless to say, we hope  that the combined synergy of the two Lecture Notes
may turn out to be very useful for the reader in order to appreciate 
the long chain
of arguments leading to the theorem that the abc-surfaces give us the
simply connected counterexamples to a weaker version of the DEF= DIFF
question raised
by Friedman and Morgan in \cite{f-m1}.

An interesting question left open (in spite of previous optimism)
concerns the canonical symplectomorphism of the (abc)-surfaces.
We discuss this and other problems, related to the connected
components of moduli spaces of surfaces of general type, and to the
corresponding symplectic structures, again in chapter 6.

The present text  not only expands the contents of the five lectures
actually held in Cetraro. Indeed,
since otherwise we would not have reached a satisfactory target,
we added the extra chapter 6.

As we already mentioned,  since the course by Manetti does not explain
the construction of his examples
(which  are here called Manetti surfaces),
we give a very brief overview of the construction, and sketch a proof
of the canonical symplectomorphism of these examples.

\vfill
\pagebreak

\section{Lecture 1: Projective and K\"ahler Manifolds, the Enriques
classification,  construction techniques.}
\label{first}

\subsection{Projective manifolds, K\"ahler and symplectic structures. }

     The basic interplay between complex algebraic geometry,
theory of complex manifolds, and theory of real symplectic
manifolds starts with projective manifolds.

We consider a closed connected
$\C$-submanifold $X^n \subset\PP^N :=  \PP^N_{\C}$.

This means that, around each point $p \in X$, there is 
a neighbourhood $U_p$ of $p$ and a permutation of the
homogeneous coordinates such that, setting
%\begin{multline}
$$ x_0 =1, \  x' := (x_1, \dots x_n) ,
  \  x'' := (x_{n+1}, \dots x_N) ,$$\label{coord}
%\end{multline}
  the intersection  $X \cap U_p$  coincides with
the graph of a holomorphic map $\Psi$:
$$ X \cap U_p = \{ (x', x'') \in U_p | x'' = \Psi (x') \}. $$
We can moreover assume, after a linear change of the homogeneous
coordinates, that the Taylor
development of
$\Psi$ starts with a second order term (i.e., $p$ is the point
$(1,0, \dots 0)$ and the projective tangent space to $X$ at $p$
is the complex subspace $\{ x'' = 0 \}$.

\begin{df}
The Fubini-Study form is the differential 2-form
$$\omega_{FS} := \frac{i}{2 \pi } \partial
\overline{\partial} log |z|^2,$$ where $z$ is the homogeneous coordinate
vector representing a point of $\PP^N $.

  In fact the above 2- form on $\C^{N+1} \setminus \{0\}$ is invariant

1) for the action of $\U(N, \C)$ on homogeneous coordinate vectors,

2) for multiplication of the vector $z$ by a nonzero holomorphic
scalar function $f(z)$ ($z$ and $ f(z) z$ represent the same point in
$\PP^N $), hence

3) $\omega_{FS}$ descends to a differential form on $\PP^N $
(being $\C^*$-invariant).
\end{df}

The restriction $\omega$ of the Fubini-Study form to a submanifold
$X $ of $\PP^n$ makes
the pair $(X,\omega)$
a K\"ahler manifold according to the following

\begin{df}
A pair $(X,\omega)$ of a complex manifold $X$, and a real differential
2- form $\omega$ is called a {\bf K\"ahler pair} if

i) $\omega$ is closed ($ d \omega  = 0$)

ii) $\omega$ is of type (1,1) $\Leftrightarrow$ for each pair of tangent vectors
$v,w$ one has ($J$ being the operator on complex
tangent vectors given by multiplication by $i = \sqrt -1$),
$$ \omega (Jv, Jw) = \omega (v, w)$$

iii) the associated Hermitian form is strictly positive definite
$\Leftrightarrow$ the real symmetric bilinear form
$\omega ( v, J w)$ is positive definite.

\end{df}

The previous definition becomes clearer if one recalls
the following  easy
  bilinear algebra lemma.

\begin{lem}
Let $V$ be a complex vector space, and $H$ a Hermitian form.
Then, decomposing $H$ in real and imaginary part,
  $$ H = S +  \sqrt -1 A ,$$ we have that $S$ is
symmetric,
$A$ is alternating,  $ S ( u,v) = A ( u, Jv)$ and $ A (Ju , J v) = A (u,v)$.

Conversely, given a real bilinear and alternating form $A$ , $A$ is
the imaginary part of a Hermitian form
$ H ( u,v) = A ( u, Jv) +  \sqrt -1 A (u,v)$ if and only if
$A$ satisfies the socalled first {\bf Riemann bilinear relation}:
$$ A (Ju , J v) = A (u,v).$$
\end{lem}

Observe that property iii) implies that $\omega$ is nondegenerate
(if in the previous lemma $S$ is positive definite, then $A$ is
nondegenerate), thus a K\"ahler pair yields a symplectic manifold
according to the standard definition

\begin{df}
A pair $(X,\omega)$ consisting of a real manifold $X$, and a real differential
2- form $\omega$ is called a {\bf symplectic pair} if

i) $\omega$ is a symplectic form, i.e., $\omega$ is closed ($ d \omega  = 0$)
  and  $\omega$ is nondegenerate at each point (thus $X$ has even dimension).

A symplectic pair $(X,\omega)$ is said to be {\bf integral} iff the De Rham
cohomology class of $\omega$ comes from $ H^2 ( X, \Z)$, or, equivalently,
there is a complex line bundle $L$ on $X$ such that $\omega$ is a first Chern
form of $L$.

An {\bf almost complex structure} $J$ on $X$ is a differentiable endomorphism
of the real tangent bundle of $X$ satisfying $ J^2 = - Id$. It is said to be

ii) {\bf compatible with $\omega$ } if
$$ \omega (Jv, Jw) = \omega (v, w),$$

iii) {\bf tame }
if the quadratic form
$\omega ( v, J v)$ is strictly positive definite.

Finally, a symplectic manifold is a manifold admitting a symplectic form
$\omega$.

\end{df}

Observe that compatibility and tameness are the symplectic geometry
translation of the two classical Riemann bilinear relations which ensure
the existence of a hermitian form, respectively the fact that
the latter is positive definite: the point of
view changes mainly in the order of the choice for $J$, resp. $\omega$.

\begin{df}
A submanifold $Y$ of a symplectic pair ($X , \omega$)  is
a {\bf symplectic} submanifold if $\omega |_Y$ is nondegenerate.

Let ($X' , \omega'$) be another symplectic pair. A diffeomorphism
$ f : X \ra X'$ is said to be a {\bf symplectomorphism}
if $ f^* (\omega' ) = \omega$.
\end{df}

Thus, unlike the K\"ahler property for complex submanifolds, the symplectic
property is not automatically inherited by submanifolds of even
real dimension.

A first intuition about symplectic submanifolds
is given by the following result, which holds more generally
on any K\"ahler manifold, and says that a good differentiable
approximation of a complex submanifold is a symplectic
submanifold.

\begin{lem}\label{approx}
Let $W \subset \PP^N$ be a differentiable submanifold of
even dimension ($ dim W = 2n$), and assume that the tangent space of $W$
is `close to be complex' in the  sense that for each  vector
$v$ tangent to $W$ there is another  vector $v'$ tangent to $W$ such that
$$  J v = v' + u, | u | <  | v| .$$

Then the restriction to $W$ of the Fubini Study form $ \omega_{FS}$
makes $W$ a symplectic submanifold of $\PP^N$.
\end{lem}

\Proof
Let $A$ be the symplectic form on projective space, so that for each
vector $v$ tangent to $W$ we have:

$  | v|^2 = A ( v, J v) = A (v, v') + A (v, u) .$

Since $ | A (v, u) | < |v|^2 $, $ A (v, v') \neq 0$ and
$A$ restricts to a nondegenerate form.

\qed

The above intuition does not hold globally, since it
  was observed by Thurston (\cite{th}) that there are symplectic
complex manifolds which are not K\"ahler. The first example of this
situation was indeed given by Kodaira (\cite{kod}) who described
the socalled
Kodaira surfaces $ \C^2 / \Ga$, which are principal
holomorphic bundles with base and fibre an elliptic curve (they are
not K\"ahler since their first Betti number equals 3).
Many more examples have been given later on.

To close the circle  between the several notions, there is
the following characterization of a K\"ahler manifold
(the full statement is very often referred to as `folklore',
but it  follows from the
statements contained in theorem 3.13 , page 74 of \cite{vois},
and proposition 4.A.8, page 210 of \cite{huy}).

{\bf K\"ahler manifolds Theorem} {\em Let ($X, \omega$)be a symplectic pair,
and let $J$ be an almost complex structure which is compatible and tame for
$\omega$. Let $g (u,v) : = \omega (u, Jv)$ be the associated Riemannian metric.
Then $J$ is parallel for the Levi Civita connection of $g$ (i.e., its
covariant derivative is zero in each direction) if and only if $J$ is 
integrable
(i.e., it yields a complex structure) and $\omega$ is a K\"ahler form. }

Returning to the Fubini-Study form, it has an important normalization
property, namely,
if we consider  a  linear subspace $\PP^m \subset \PP^N$
(it does not matter which one, by the unitary invariance mentioned in
1)  above),
then integration in pluripolar coordinates yields
$$ \int_{\PP^m} \frac1{m!} \omega_{FS}^{m}  = 1.$$

The above equation, together with  Stokes' Lemma, and
a multilinear algebra calculation for which we refer
for instance to  Mumford's book \cite{mum} imply

{\bf Wirtinger's Theorem} {\em  Let $X:= X^n$ be a complex submanifold
of $\PP^N$. Then $X$ is a volume minimizing submanifold for the
n-dimensional Riemannian volume function  of submanifolds $M$
of real dimension $2n$, $$ vol (M) := \int  d Vol_{FS},$$

where $d Vol_{FS} = \sqrt {det (g_{ij})(x)} \  |dx|$ is the volume 
measure of the Riemannian metric
$g_{ij}(x)$ associated to the Fubini Study form.
Moreover, the global volume of $X$ equals a positive integer, called the
{\bf degree} of $X$.}

The previous situation is indeed quite more general:

Let ($X , \omega$) be a symplectic manifold, and let $Y$ be an
oriented submanifold of even dimension $= 2m$:
then the global symplectic volume of $Y$
$ vol (Y) := \int_Y  \frac1{n!} \omega^{m}$ depends only on the homology class
of $Y$, and will be an integer if the pair ($X , \omega$) is integral
(i.e., if the De Rham class of $\omega$ comes from $ H^2 (X, \Z)$).

If moreover $X$ is K\"ahler, and $Y$ is a complex submanifold,
then $Y$ has a natural orientation, and one has the

{\bf Basic principle of K\"ahler geometry:}
{\em Let $Y$ be a compact submanifold  of a K\"ahler manifold $X$: then
$ vol (Y) := \int_Y  \omega^{m} > 0$, and in particular the cohomology class
of $Y$ in $ H^{2m} (X, \Z)$ is nontrivial.}

The main point of the basic principle is that the integrand of
$ vol (Y) := \int_Y  \omega^{m}$ is pointwise positive,
because of condition iii). So we see that a similar principle
holds more generally if we have a symplectic manifold $X$
and a compact submanifold $Y$ admitting an
almost complex structure compatible and tame
for the restriction of $\omega$ to $Y$.

Wirtinger's  theorem  and  the following theorem of Chow provide
the link with algebraic geometry mentioned in the beginning.

{\bf Chow's Theorem} {\em  Let $X:= X^n$ be a (connected) complex submanifold
of $\PP^N$. Then $X$ is an algebraic variety, i.e., $X$ is the locus of zeros
of a homogeneous prime ideal $\PPP$ of the polynomial ring
$\C[x_0, \dots x_N]$.}

We would now like  to show how Chow's theorem is a consequence of
another result:

{\bf Siegel's Theorem} {\em  Let $X:= X^n$ be a compact
(connected) complex manifold of (complex) dimension $n$.
  Then the field $\C^{Mer}(X)$ of meromorphic functions on $X$
is finitely generated, and its
transcendence degree over $\C$ is at most $n$.}

The above was proven by Siegel  just using the lemma of Schwarz
and an appropriate choice of a  finite cover of a compact complex manifold
made by polycylinder charts (see  \cite{sieg}, or \cite{corn}).

{\em Idea of proof of Chow's theorem.}

Let $p \in X$ and take coordinates as in \ref{coord}: then
we have an injection $ \C (x_1, \dots x_n ) \hookrightarrow
\C^{Mer}(X)$, thus $\C^{Mer}(X)$ has transcendency degree
$n$ by Siegel's theorem.

Let $ Z$ be the Zariski closure of $X$: this means that $Z$ is the set of
zeros of the homogeneous ideal $\I_X \subset \C [x_0, \dots x_N]$
generated by the homogeneous polynomials vanishing on $X$.

Since $X$ is connected, it follows right away, going to
nonhomogeneous coordinates and using that the ring
of holomorphic functions on a connected open set is
an integral domain, that the ideal $\I_X = \I_Z$ is a prime ideal.

We consider then the homogeneous coordinate ring $\C [Z] : =
  \C [x_0, \dots x_N] / \I_X $ and the field of rational
functions $\C (Z)$, the field of the fractions  of the integral domain $\C [Z] $
which are homogeneous of degree
$0$. We observe that we
have an injection $\C (Z) \hookrightarrow \C^{Mer}(X) $.

Therefore $ \C (x_1, \dots x_n ) \hookrightarrow
\C (Z)   \hookrightarrow \C^{Mer}(X) $. Thus
the field of rational functions $\C (Z)$ has transcendency
degree $n$ and $Z$ is an irreducible algebraic subvariety
of $\PP^N$ of dimension $n$. Since the smooth locus $Z^* : = Z \setminus Sing (Z)$
is dense in $Z$ for the Hausdorff topology,  is connected, and
contains $X$, it follows that $ X = Z$.

\qed

The above theorem extends to the singular case: a closed
complex analytic subspace of $\PP^N$ is also a closed set
in the Zariski topology, i.e., a closed algebraic set.

We have seen in the course of the proof that the dimension of
an irreducible projective variety is given by the
transcendency degree over $\C$ of the field $\C (Z)$ (which, by a further
extension of Chow's theorem, equals $\C^{Mer} (Z)$).

The  degree of $Z$ is then  defined through the

{\bf Emmy Noether Normalization Lemma.} {\em Let $Z$ be an irreducible
subvariety of $\PP^N$ of dimension $n$: then for general choice
of independent linear forms $(x_0, \dots x_n)$ one has that the
homogeneous coordinate ring of $Z$, $\C [Z] : =
  \C [x_0, \dots x_N] / \I_Z $  is an integral extension of
  $ \C [x_0, \dots x_n] $. One can view  $\C [Z]$ as a torsion free
$ \C [x_0, \dots x_n] $-module, and its rank is called the degree $d$ of $Z$.}

The geometrical consequences of Noether's normalization are
(see \cite{shaf}):

\begin{itemize}
\item
The linear projection with centre $ L : = \{ x | x_0= \dots x_n = 0)$,
$ \pi_L : \PP^N \setminus L \ra \PP^n$ is defined on $Z$ since
$ Z \cap L =\emptyset$, and $\pi: = \pi |_L : X \ra \PP^n$
is surjective and finite.
\item
For $ y \in \PP^n$, the finite set $\pi^{-1} (y)$
has cardinality at most $d$, and equality holds for
$ y$ in a Zariski open set $ U \subset \PP^n$.

\end{itemize}

The link between the volume theoretic and the algebraic notion
of degree is easily obtained via the Noether projection $\pi_L$.

In fact, the formula $ (x_0, x', x'') \ra (x_0, x',  (1-t) x'')$
provides a homotopy between the identity map of $Z$
and a covering of $\PP^n$ of degree $d$, by which it follows
that  $ \int_{Z^*} \omega_{FS}^n $ converges and equals precisely $d$.

We end this subsection by fixing  the standard notation:
for $X$ a projective variety, and $x$ a point in $X$ we denote by $\hol_{X,x}$
the local ring of  algebraic functions on $X$ regular in $x$, i.e.,
$$ \hol_{X,x} : = \{ f \in \C(X) |  \exists  a , b \in \C [X],
{\rm  homogeneous},
s.t. f = a/b  \  {\rm and } \  b (x) \neq 0 \}.$$

This local ring is contained in the local ring of restrictions
of local holomorphic functions from $\PP^N$, which we
denote by $ \Hol_{X,x} $.

The pair $ \hol_{X,x} \subset \Hol_{X,x} $ is a faithfully flat ring extension,
according to the standard

\begin{df}\label{flat}
   A ring extension $A \ra B$  is said to be {\bf flat},
respectively {\bf faithfully flat},
  if the following property holds: a complex of $A$-modules
$ ( M_i, d_i) $ is exact only if (respectively, if and only if)
$ ( M_i \otimes_A B, d_i \otimes_A B) $ is exact.
\end{df}

This basic algebraic property underlies the so called (see \cite{gaga})

{\bf G.A.G.A. Principle. } {\em  Given a projective variety, and
a coherent (algebraic) $\hol_X$-sheaf $\FFF$,
let $\FFF^h : = \FFF \otimes_{\hol_X}\Hol_X$  be the corresponding
holomorphic coherent sheaf: then one has a natural isomorphism
of cohomology groups $$ H^i (X_{Zar}, \FFF) \cong
H^i (X_{Haus}, \FFF^h) ,$$
where the left hand side stands for $\check{C}ech$ cohomology taken
in the Zariski topology, the right hand side stands for
$\check{C}ech$ cohomology taken
in the Hausdorff topology. The same holds replacing $\FFF$ by
$\hol_X^*$.}

Due to the GAGA principle, we shall sometimes make some abuse
of notation, and simply write, given a divisor $D$ on $X$,
$ H^i (X, D)$ instead of $ H^i (X, \hol_X (D))$.

\subsection{The birational equivalence of algebraic varieties}

A rational map of a (projective) variety $\phi : X \dasharrow \PP^N$ is given
through $N$ rational functions $ \phi_1, \dots \phi_N$.

Taking a common multiple $s_0$ of the denominators $b_j$ of $ \phi_j = a_j/ b_j$,
we can write $\phi_j = s_j  / s_0$, and write $ \phi = (s_0, \dots s_N)$,
where the $s_j $'s are all homogeneous of the same degree,
whence they define a graded homomorphism $\phi^* :  \C[\PP^N] \ra  \C[X]$.

The kernel of $\phi^*$ is a prime ideal, and its zero locus, denote it by
$Y$, is called the image of $\phi$, and we say that $X$ dominates $Y$.

One says that $\phi$ is a morphism in $p$ if there is such a representation
$ \phi = (s_0, \dots s_N)$ such that some $s_j(p) \neq 0$. One can see
that there is a maximal open set $ U \subset X$ such that $\phi$ is  a morphism
on  $U$, and that $ Y = \overline {\phi (U)  }$.

If  the local rings $\hol_{X,x}$ are factorial, in particular if $X$ is smooth,
then one can take at each point $x$ relatively prime elements $a_j, b_j$,
let $s_0$ be the least common multiple of the denominators, and it follows
then that the {\bf Indeterminacy Locus $ X \setminus U$} is a closed set
of codimension at least $2$. In particular, every rational map
of a smooth curve is a morphism.

\begin{df}
Two algebraic varieties $X,Y$ are said to be {\bf birational }
iff their function fields
$\C(X), \C(Y)$ are isomorphic, equivalently if there are two
dominant rational maps $ \phi : X \dasharrow Y, \psi : Y \dasharrow X$,
which are inverse to each other. If $\phi , \psi = \phi^{-1}$ are morphisms,
then $X$ and $Y$ are said to be {\bf isomorphic}.
\end{df}

By Chow's theorem, biholomorphism and isomorphism is the same notion for
projective varieties (this ceases to be true in the non compact case,
cf.\cite{ser2}).

Over the complex numbers, we have (\cite{hiro})

{\bf Hironaka's theorem on resolution of singularities. } {\em  Every
projective variety is birational to a smooth projective variety.}

As we already remarked, two birationally equivalent curves are isomorphic,
whereas for a smooth surface $S$, and a point $p \in S$, one may consider
the blow -up of the point $p$, $ \pi : \hat{S} \ra S$. $\hat{S}$ is obtained
glueing together $S \setminus \{p\}$ with the closure
of the projection with centre $p$,  $ \pi_{p} : S\setminus \{p\}
\ra \PP^{N-1}$.
One can moreover show that $\hat{S}$ is projective.
The result of blow up is that the point $p$ is replaced by
the projectivization of the tangent plane to $S$ at $p$, which
is a curve $E \cong \PP^1$, with normal sheaf $ \hol_E (E) \cong \hol_{\PP^1} (-1)$.
In other words, the selfintersection of $E$, i.e., the
degree of the normal bundle of $E$, is $ -1$, and we simply say that $E$
is an {\bf Exceptional curve of the I Kind}.

{\bf Theorem of Castelnuovo and Enriques.} {\em Assume that a smooth projective
surface $Y$ contains an irreducible curve $E \cong \PP^1$ with
selfintersection $E^2 = -1$: then there is a birational morphism
$ f : Y \ra S$ which is  isomorphic to the blow up $ \pi : \hat{S} \ra S$
of a point $p$ (in particular $E$ is the only curve contracted to a point by $f$).}

The previous theorem justifies the following

\begin{df}
A smooth projective surface is said to be  {\bf minimal } if it does not
contain any exceptional curve of the I kind.
\end{df}

One shows then that every birational transformation is a composition of
blow ups and of inverses of blow ups, and each surface $X$ is birational to
a smooth minimal surface $S$. This surface $S$ is unique, up to isomorphism,
if $X$ is not ruled (i.e., not birational to a product $C \times \PP^1$),
by the classical

{\bf Theorem of Castelnuovo.} {\em Two birational minimal models $S$, $S'$
are isomorphic unless they are birationally {\bf ruled}, i.e., birational to a product $C
\times
\PP^1$, where $C$ is a smooth projective curve.
In the ruled case, either $S$  $\cong \PP^2$, or
$S$ is isomorphic to the projectivization $\PP(V)$
of a rank 2 vector bundle $V$
on $C$.}

Recall now that  a variety $X$ is smooth if and only if the
sheaf of differential forms $\Omega^1_X$ is locally free, and
locally generated by $ dx_1, \dots d x_n$, if $x_1, \dots x_n$ yield
local holomorphic coordinates.

The vector bundle (locally free sheaf) $\Omega^1_X$ 
and its associated bundles provide
birational invariants in view of the classical (\cite{arcata})

{\bf K\"ahler's lemma.} {\em Let $ f : X^n \dasharrow Y^m$ be a
dominant rational map between smooth projective varieties of
respective dimensions $n, m$. Then one has injective pull back
linear maps $H^0 ( Y, {\Om^1_Y }^{\otimes r} ) \ra
H^0 ( X, {\Om^1_X }^{\otimes r} ) $. Hence the vector spaces
$ H^0 ( X, {\Om^1_X }^{\otimes r_1} \otimes \dots \otimes
{\Om^n_X }^{\otimes r_n}) $ are birational invariants.}

Of particular importance is the top
exterior power $\Omega^n_X= \La^n (\Omega^1_X)$, which is locally free
of rank 1, thus can be written as $ \hol_X (K_X)$ for a suitable
Cartier divisor $K_X$, called the {\bf canonical divisor}, and well defined only
up to linear equivalence.

\begin{df}
The {\bf i-th pluriirregularity} of a smooth projective variety $X$ is
the dimension $ h^{0,i} : = dim  ( H^i (X, \hol_X)),$ which by Hodge Theory
equals $ dim (H^0 ( X, \Om^i_X ))$. The {\bf m-th plurigenus $P_m$ }
is instead the dimension $ P_m (X) : = dim (H^0 ( X, {\Om^n_X }^{\otimes m})) =
h^0 ( X , m K_X).$
\end{df}

A finer birational invariant is the canonical ring of $X$.
\begin{df}
The {\bf canonical ring} of a smooth projective variety $X$ is
the graded ring $$ \RRR (X) : = \bigoplus_{m=0}^{\infty} H^0 ( X , m K_X) .$$

If  $ \RRR (X) = \C$  one defines $ Kod (X) = - \infty$, otherwise
the {\bf Kodaira dimension} of $X$ is defined as
the transcendence degree over $\C$ of the {\bf canonical subfield of $\C(X)$},
  given by the field $\sQ (X)$ of homogeneous fractions of
degree zero of $ \RRR (X)$.

$X$ is said to be of {\bf general type} if its Kodaira dimension is
maximal (i.e., equal to the dimension $n$ of $X$).
\end{df}

As observed in \cite{andrcime} $\sQ (X)$ is algebraically closed
inside $\C(X)$, thus one obtains that $X$ is of general type if and only if
there is a positive integer $m$ such that $H^0 ( X , m K_X) $
yields a birational map onto its image $\Sigma_m$.

  One of the more crucial questions in classification theory is
whether the canonical ring of a variety of general type is
finitely generated, the answer being affirmative (\cite{mum1},
\cite{mori}) for dimension $ n \leq 3$ \footnote{The question seems 
to have been settled for varieties of general type, and with a positive answer.}.

\subsection{ The Enriques classification: an outline. }

The main discrete invariant of smooth projective curves $C$ is the
{\bf genus } $ g (C) : = h^0 (K_C) = h^1 (\hol_C)$.

It determines easily the Kodaira dimension, and the Enriques classification
of curves is the subdivision

\begin{itemize}
\item
$ Kod (C) = - \infty \Leftrightarrow g (C) = 0 \Leftrightarrow  C \cong \PP^1.$
\item
$ Kod (C) = 0 \Leftrightarrow g (C) = 1 \Leftrightarrow  C \cong \C /
(\Z + \tau \Z)$, with $ \tau \in \C , Im (\tau )> 0 $
$ \Leftrightarrow C$ is an elliptic curve.
\item
$ Kod (C) = 1 \Leftrightarrow g (C) \geq 2
  \Leftrightarrow  C$ is of general type.

\end{itemize}

Before giving  the Enriques classification of projective surfaces
over the complex numbers, it is convenient to discuss further the birational
invariants of surfaces.

\begin{oss}
An important birational invariant of smooth varieties $X$ is the fundamental
group $\pi_1 (X)$.

For surfaces, the most important invariants are :

\begin{itemize}
\item
the {\bf irregularity} $ q : = h^1 (\hol_X)$
\item
the {\bf geometric genus} $ p_g  : = P_1 : = h^0 (X, K_X)$, which for surfaces
  combines with the
irregularity to give the {\bf holomorphic Euler-Poincar\'e characteristic}
$ \chi(S) : =  \chi(\hol_S) : =1 -q + p_g$
\item
the {\bf bigenus} $ P_2 : = h^0 (X, 2 K_X)$ and
especially the {\bf twelfth plurigenus} $ P_{12} : = h^0 (X, 12 K_X)$.
\end{itemize}
If $S$ is a non ruled minimal surface, then also the following
are birational invariants:
\begin{itemize}
\item
the selfintersection of a canonical divisor $K^2_S$, equal to $c_1(S)^2$,
\item
the {\bf topological Euler number} $ e(S)$, equal to $c_2(S)$ by the Poincar\'e
Hopf theorem, and which by
{\bf Noether's theorem} can also be expressed as
  $$ e(S) =  12 \chi(S) - K^2_S = 12 (1 - q + p_g) -
K^2_S,$$
\item
the {\bf topological index} $\sigma (S)$ (the index of the quadratic form

$q_S : H^2 (S, \Z)  \times H^2 (S, \Z) \ra \Z$), which, by the
{\bf Hodge index theorem}, satisfies the equality
  $$ \sigma (S) = \frac13  (K^2_S - 2 e(S)),$$
\item
  in particular, all the Betti numbers $b_i(S)$ and
\item
  the positivity $b^+ (S)$ and the negativity  $b^-(S)$ of $q_S$
(recall that $b^+(S)  + b^-(S) = b_2(S)$).
\end{itemize}

\end{oss}

The Enriques classification of complex algebraic surfaces gives a  very simple
description of the surfaces with nonpositive Kodaira dimension:

\begin{itemize}
\item
$S$ is a ruled surface of irregularity $g$ $\iff : $

$\iff : $ $S$ is birational to a product $C_g
\times \PP^1$, where $C_g$ has genus $g$ $\iff$

$\iff$ $ P_{12} (S) = 0, q (S)  = g$$\iff$

$\iff$ $ Kod (S) = - \infty$, $ q(S) = g$.

\item
$S$ has $ Kod (S) = 0$ $\iff$ $ P_{12} (S) = 1$.
\end{itemize}

\medskip

\noindent
There are four classes of such surfaces with $ Kod (S) = 0$:
\begin{itemize}
\item
Tori $\iff$ $ P_{1} (S) = 1 , q(S) = 2$,
\item
K3 surfaces $\iff$ $ P_{1} (S) = 1 , q(S) = 0$,
\item
Enriques surfaces $\iff$ $ P_{1} (S) = 0 , q(S) = 0$, $ P_{2} (S) = 1$,
\item
Hyperelliptic surfaces $\iff$ $ P_{12} (S) = 1 , q(S) = 1.$
\end{itemize}

Next come the surfaces with strictly positive Kodaira dimension:

\begin{itemize}
\item
$S$ is a properly elliptic surface  $\iff : $

$\iff : $ $ P_{12} (S) > 1$, and $
H^0(12 K_S)$ yields a map to a curve with fibres elliptic 
curves $\iff$

$\iff$ $S$ has $ Kod (S)
= 1$ $\iff$

$\iff$ assuming that $S$  is minimal: $ P_{12} (S) > 1$ and $ K^2_S = 0$.
\item
$S$ is a  surface  of general type $\iff : $

$\iff :$ $S$ has $ Kod (S)
= 2$ $\iff$

$\iff  $ $ P_{12} (S) > 1$, and $
H^0(12 K_S)$ yields a  birational map onto its image 
$\Sigma_{12}$ $\iff$

$\iff$ assuming that $S$  is minimal: $ P_{12} (S) > 1$ and $ K^2_S \geq 1$.

\end{itemize}

\subsection{ Some constructions of projective varieties. }

Goal of this subsection is first of all to illustrate concretely the meaning of
the concept 'varieties of general type'. This means, roughly speaking, that if
we have a construction of varieties of a fixed dimension involving some integer 
parameters, most of the time we get varieties of  general type 
when these parameters are
all sufficiently large.

[1] {\bf Products. }

Given projective varieties $X \subset \PP^n$ and  $Y\subset \PP^m$,
their product $ X \times Y$ is also projective.
This is an easy consequence of the fact that  the product
$ \PP^n \times  \PP^m$ admits the Segre embedding in $\PP^{mn + n + m} \cong
\PP  (Mat (n,m ))$ onto the subspace of rank one matrices,
given by the morphism $ (x,y) \ra   x  \cdot  \ ^t y$.

[2] {\bf Complete intersections. }

Given a smooth variety $X$, and divisors $ D_1 = \{ f_1 = 0\}, \dots , D_r
= \{ f_r = 0\}$ on $X$,
their intersection $Y =  D_1 \cap \dots  \cap D_r$ is said to be a complete intersection
if $Y$ has codimension $r$ in $X$. If $Y$ is smooth, or, more generally, reduced, 
locally its ideal is generated by the local equations of the $D_i$'s ( $ \I_Y = ( f_1, \dots
f_r)$).

$Y$ tends to inherit much from the geometry of $X$, for instance,
if $ X = \PP^N$ and $Y$ is smooth of dimension $ N- r \geq 2$,
then $Y$ is simply connected by the theorem of Lefschetz.

[3] {\bf  Finite coverings according to Riemann, Grauert and Remmert. }

Assume that $Y$ is a normal variety (this means that each local ring
$ \hol_{X,x} $ is integrally closed in the function field $ \C (X)$), and
that $ B$  is a closed subvariety of $Y$ (the letter $B$ stands for 
'branch locus).

Then there is (cf. \cite{g-r}) a correspondence between

[3a] subgroups $ \Ga \subset \pi_1 ( Y \setminus B)$ of finite index, and

[3b] pairs $ (X, f) $ of a normal variety $X$ and a finite map $ f : X \ra Y$
which,  when restricted to $ X \setminus f^{-1} (B)$,
is a local biholomorphism and a topological covering space of $   Y \setminus B$.

The datum of the covering is equivalent to the datum of the sheaf of
$\hol_Y$-algebras $ f_* \hol_X$.  As an $\hol_Y$-module $ f_* \hol_X$
is locally free if and only if $f$ is {\bf flat} (this means that, $\forall x \in X$, 
$\hol_{Y, f(x)} \ra \hol_{X,x} $ is flat), and this
is indeed the case when $f$ is finite and $Y$ is smooth.

[4] {\bf  Finite Galois coverings. }

Although this is just a special case of the previous one, namely when
$\Ga$ is a normal subgroup with factor group $ G : = \pi_1 ( Y 
\setminus B) / \Ga$,
in the more special case (cf. \cite{Pardini}) where $G$ is Abelian and $Y$ 
is smooth, one can
give explicit equations for the covering.  This is due to the fact that all
irreducible representations of an abelian group are 1-dimensional,
so we are in the {\bf split case} where $ f_* \hol_X$ is a direct sum of
invertible sheaves.

The easiest  example is the one of

[4a] Simple cyclic coverings of degree $n$.

In this case there is

i) an invertible sheaf $\hol_Y (L)$ such that
$$ f_* \hol_X = \hol_Y \oplus \hol_Y (- L) \oplus \dots \oplus \hol_Y 
(-(n-1)L).$$
ii) A section $ 0 \neq \sigma \in H^0 (\hol_Y (n L))  $ such that
$X$ is the divisor, in the geometric line bundle $\LL$ whose sheaf of 
regular sections is
$\hol_Y (L)$, given by the equation $  z^n = \sigma (y)$.

Here, $z$ is the never vanishing
section of $ p^* (\hol_Y (L))$ giving a
  tautological linear form on the fibres of $\LL$: in other words, one has
an  open cover $ U_{\alpha}$ of $Y$ which is trivializing for $\hol_Y (L)$,
and $X$ is obtained by glueing together the local equations
$  z_{\alpha}^n = \sigma_{\alpha} (y)$, since
$  z_{\alpha} = g _{\alpha, \beta}(y) z_{\beta},$ $  \sigma_{\alpha}(y) =
g _{\alpha,
\beta}(y)^n \sigma_{\beta}(y)$.

One has as branch locus $B = \De :=  \{ \sigma = 0\}$, at least if 
one disregards the
multiplicity (indeed $ B = (n-1) \Delta$). Assume $Y$ is smooth: then $X$ is
smooth iff $\De$ is smooth, and, via direct image,
all the calculations of cohomology groups of basic sheaves on  $X$ 
are reduced to
calculations for corresponding sheaves on $Y$. For instance, since $ 
K_X = f^*( K_Y +
(n-1) L)$, one has:
$$ f_* (\hol_X (K_X)) = \hol_Y (K_Y) \oplus \hol_Y (K_Y + L) \oplus \dots \oplus
\hol_Y ( K_Y + (n-1) L)$$
(the order is exactly as above according to the characters   of the 
cyclic group).

We see in particular that $X$ is of general type if $L$ is sufficiently positive.

[4b] Simple iterated cyclic coverings.

Suppose that we take a simple cyclic covering $ f : Y_1 \ra Y$ as 
above, corresponding
to the pair $ (L, \sigma )$, and we want to consider again a simple 
cyclic covering
of $Y_1$. A small calculation shows that it is not so easy to describe
$H^1(\hol^*_{Y_1})$ in terms of the triple $ (Y, L, \sigma )$; but in any case
$H^1(\hol^*_{Y_1}) \supset H^1(\hol^*_{Y})$.
Thus one defines an {\bf iterated simple cyclic covering }
as the composition of a chain of simple cyclic coverings
$ f_i : Y_{i+1} \ra Y_i$, $i=0, \dots k-1$
(thus $ X : = Y_k$, $ Y : = Y_0$) such that at each step
 the divisor $ L_i$ is the pull back
of a divisor on $ Y = Y_0$.

In the case of iterated double coverings, considered in  \cite{man3},
we have at each step $ ( z_i)^2  = \sigma_i$ and each $ \sigma_i$ is
written as $ \sigma_i =  b_{i,0} + b_{i,1} z_1 + b_{i,2} z_2 +  \dots + b_{i, 1, \dots i-1} 
z_1 \dots z_{i-1}$,
where, for $ j_1< j_2 \dots <  j_h $, we are given a section $ b_{i,j_1, \dots j_h} \in
H^0 ( Y,
\hol_Y (2  L_i - L_{j_1} -
\dots - L_{j_h} ))$.

In principle, it looks like one could describe the Galois covers
with solvable Galois group
$G$ by considering iterated cyclic coverings, and
then imposing the Galois condition.
But this does not work without resorting
to more complicated cyclic covers and to special geometry.

[4c] Bidouble covers (Galois with group $ ( \Z/2)^2$).

The {\bf  simple bidouble covers} are  simply the fibre product of 
two double covers,
thus here $X$ is the complete intersection of the following two divisors
$$ z^2 = \sigma_0, \  w^2 = s_0 $$
in the vector bundle $ \LL \oplus \M$.

These are the examples we shall  mostly  consider.

More generally, a  bidouble cover of a smooth variety $Y$ occurs (\cite{cat1})
as the subvariety $X$ of the direct sum of 3 line bundles $ \LL_1 \oplus  \LL_2
\oplus  \LL_3$, given by equations

$$
Rank  \left( \begin{array}{ccc}
        x_1 & w_3 & w_2 \\
        w_3 & x_2 & w_1 \\
        w_2 & w_1 & x_3
      \end{array} \right) = 1 \qquad \qquad (*)
$$

Here, we have $3$ Cartier divisors  $D_j = div (x_j)$ on $Y$ and
  $3$ line bundles $\LL_i$, with fibre coordinate $w_i$,
such that  the following linear
equivalences hold on $Y$,

$$ L_i + D_i \equiv L_j + L_k , $$

  {\rm for each permutation } $(i,j,k)$  of $(1,2,3)$.

One has  :
$f_* {\hol}_X = {\hol}_Y \bigoplus  (\oplus_i {\hol}_Y (- L_i)).$

Assume in addition that $Y$ is a smooth variety, then :
\begin{itemize}
\item $X$ is normal if and only if the divisors $D_j$ are reduced
and have no common components .
\item $X$ is smooth if and only if the divisors $D_j$ are smooth, they
do not have a common intersection and have pairwise transversal
intersections.
\item $X$ is Cohen - Macaulay and for its dualizing sheaf $\omega_X$
(which, if $Y$ is normal, equals the sheaf of Zariski differentials that we shall discuss
later) we have

  $f_* \omega_X = {\mathcal H}om_{\hol_Y}(f_* {\hol}_X , \omega_Y )
  = \omega_Y \bigoplus (\oplus_i \omega_Y (L_i))$.
\end{itemize}

[5] Natural deformations.

One should in general consider Galois covers as `special varieties'.

For instance, if we have a line bundle $\LL$ on $Y$, we consider
in it the divisor $X$ described by an equation
$$ z^n + a_2 z^{n-2} + \dots a_{n-1} z + a_n = 0 , \ for \  a_i
\in H^0 ( Y, \hol_Y (i
L)).$$

It is clear that we obtain a simple cyclic cover if we set $ a_n = - \sigma_0$,
and, for $ j \neq n$, we set $a_j =  0 $.

The family of above divisors (note that we may  assume $ a_1 = 0$ 
after performing a Tschirnhausen transformation) is called the  family of
{\bf natural deformations} of a simple cyclic cover.

One can define more generally a similar concept for any Abelian covering.
In particular, for simple bidouble covers, we have the following family
of natural deformations
  $$ z^2 = \sigma_0 (y) + w \sigma_1 (y) , \  w^2 = s_0 (y)  + z s_1 (y),$$
where $ \sigma_0  \in H^0 ( Y, \hol_Y (2
L)), \sigma_1  \in H^0 ( Y, \hol_Y (2
L - M)), s_0  \in H^0 ( Y, \hol_Y (2
M)), s_1  \in H^0 ( Y, \hol_Y (2
M - L)) $.

[6] Quotients.

In general, given an action of a finite group $G$ on the function
field $\C (X)$ of a variety $X$, one can always take the birational
quotient, corresponding to the invariant subfield   $\C (X)^G$.

Assume that $X \subset \PP^N$ is a projective variety and
that we have a finite group $ G \subset \PP GL (N+1, \C),$
such that $ g (X) = X , \ \forall g \in G.$

We want then to construct a biregular quotient $ X / G$
with a projection
  morphism $\pi : X \ra X / G$.

For each point $x \in X$ consider a hyperplane $ H $ such that
$ H \cap Gx = \emptyset$, and let $ U : = X \setminus (\cup_{g \in G} \ g(H))$.

$U$ is an invariant affine subset, and we consider on the quotient
set $U/G$ the ring of invariant polynomials $ \C [U]^G$,
which is finitely generated since we are in characteristic zero and
we have a projector onto the subspace of invariants.

It follows that if $X$ is normal, then also $X / G$ is normal, \and 
moreover projective
since there are very ample $g$-invariant divisors on $X$.

If $X$ is smooth, one has that $X / G$ is smooth if

1) $G$ acts freely  or,  

more generally, if and only if

2) for each point $ p \in X$, the stabilizer subgroup $G_p : = \{ g | 
g(p) = p \}$
is generated by pseudoreflections (theorem of Chevalley, cf. 
for instance \cite{dolg}).

To explain the meaning of a {\bf pseudoreflection}, observe that, if 
$ p \in X$
is a smooth point, by a theorem of Cartan (\cite{cartan}), one can 
{\em linearize} the
action of $G_p$, i.e., there exist local holomorphic coordinates 
$z_1, \dots z_n$ such that
the action in these coordinates is linear. Thus,  $g \in G_p$ acts by
$ z \ra A(g) z$, and one says that $g$ is a pseudoreflection if $ 
A(g)$ (which is
diagonalizable, having finite order) has $(n-1)$ eigenvalues equal to $1$.

[7] {\bf Rational Double Points = Kleinian singularities}.

These are exactly the quotients  $ Y = \C^2 / G$ by the action of a finite
group $ G \subset SL (2, \C)$.  Since $ A(g) \in SL (2, \C)$ it
follows that $G$ contains no pseudoreflection, thus $Y$ contains exactly one
singular point $p$, image of the unique point with
a nontrivial stabilizer, $ 0 \in \C^2$.

These singularities  $ (Y,p)$ will play a prominent  role
in the next section.

In fact, one of their properties is due to the fact  that the differential
form $ d z_1 \wedge d z_2$ is
$G$-invariant (because
$ det (A(g)) = 1$), thus the sheaf $ \Omega^2_Y$ is trivial on $ Y 
\setminus \{ p \}$.

Then the dualizing sheaf  $\omega_Y = i_* (\Omega^2_{Y \setminus \{ p \}})$
is also trivial.

\newpage

\section{Lecture 2: Surfaces of general type and their Canonical models.
  Deformation equivalence
and singularities.}
\label{second}

\subsection{Rational double points}

Let us take up again the Kleinian singularities introduced in the 
previous section

\begin{df}
A Kleinian singularity is a singularity $ (Y,p)$ analytically isomorphic to
a quotient singularity $  \C^n / G$  where $G$ is a finite
subgroup $ G \subset SL (n, \C)$.
\end{df}

\begin{ex}. The surface singularity $A_n$ corresponds to the cyclic
group $ \mu_n \cong \Z / n$ of  n-th roots of unity acting with
characters $1$ and $(n-1)$.

I.e., $\zeta \in  \mu_n$ acts by $ \zeta (u,v) :=  ( \zeta u,  \zeta ^{n-1}v)$, and the ring
of  invariants is
the ring $ \C [x,y,z] / (xy - z^n)$, where
$$ x : = u^n , \ y : = v^n , z : = u v.   $$
\end{ex}

\begin{ex}\label{Riem} One has more generally the cyclic quotient surface
singularities corresponds to the cyclic group $ \mu_n \cong \Z / n$ 
of  n-th roots
  of unity acting with
characters $a$ and $b$, which are denoted by $ \frac1n (a,b)$.

Here, $ \zeta (u,v) :=  ( \zeta^a  u,  \zeta ^{b}v)$.

We compute the ring of invariants in the case $n=4, a=b=1$:
   the ring of invariants
is   generated by
$$ y_0 : = u^4 ,  y_1 : = u^3 v, \  y_2 : = u^2 v^2, \  y_3 : = u v^3 
, y_4 : = v^4 ,   $$
and the ring is $ \C [y_0, \dots , y_4] / J$, where $J$ is the ideal 
of  $2 \times 2$
minors of the matrix
$\begin{pmatrix}y_0 & y_1 & y_2& y_3\\  y_1 & y_2& y_3 & y_4 \end{pmatrix}$,
or equivalently of the matrix
$\begin{pmatrix}y_0 & y_1 & y_2\\  y_1 & y_2& y_3 \\ y_2& y_3 & y_4
\end{pmatrix}$.
The first realization of the ideal $J$ corresponds to the identification
of the singularity $Y$ as the cone over a rational normal curve of degree $4$
(in $\PP^4$), while in the second $Y$ is viewed as a linear section of
the cone over the Veronese surface.

We observe that $ 2 y_2$ and $y_0 + y_4$ give a map to $\C^2$ which is
finite of degree 4. They are invariant for the group of order 16
generated by 
$$ (u,v) \mapsto (iu, iv), \ (u,v) \mapsto (iu, - iv), \ (u,v) \mapsto (v, u), $$
hence $Y$ is a bidouble cover of  \ $\C^2$ branched on three lines passing
through the origin (cf. (*), we  set $ x_3 := x_1 - x_2$ and we 
choose as branch
divisors
$x_1, x_2, x_3 := x_1 - x_2$).
\end{ex}

In dimension two, the classification of Kleinian singularities is a nice
chapter of geometry  ultimately going back to Thaetethus' Platonic solids.
Let us briefly recall it.

First of all, by averaging the positive definite Hermitian product in $\C^n$,
one finds that a finite subgroup $ G \subset SL (n, \C)$ is conjugate
to a finite subgroup $ G \subset SU (n, \C)$. Composing the inclusion
$ G \subset SU (n, \C)$  with the surjection
$SU (n, \C) \ra \PP SU (n, \C) \cong SU (n, \C) / \mu_n$
  yields a finite group $G'$ acting
on $\PP^{n-1}$.

Thus, for $n=2$, we get $G' \subset \PP SU (2, \C) \cong SO(3)$ acting on
the Riemann sphere $\PP^1 \cong S^2$.

The consideration of the Hurwitz formula for
the quotient morphism $\pi : \PP^1 \ra \PP^1 /G'$,
and the fact that $\PP^1 /G'$ is a smooth curve of genus $0$, (hence
$\PP^1 /G' \cong \PP^1$) allows the classification of such groups $G'$.

Letting in fact $p_1, \dots p_k$ be the branch points of $\pi$,
and $m_1, \dots m_k$ the respective multiplicities (equal to the order in $G'$
of the element corresponding to the local monodromy), we have
 {\bf  Hurwitz's formula} (expressing the degree of the canonical divisor
$K_{\PP^1}$ as the sum of the degree of the pull back of $K_{\PP^1}$
with the degree of the ramification divisor)
$$ - 2 = | G'| ( -2 + \sum_{i=1}^k  [ 1 - \frac1{m_i} ] ).$$
Each term in the square bracket is $ \geq \frac12$, and the left hand side
is negative: hence $ k  \leq 3.$

The situation to classify is the datum of a ramified covering of
$\PP^1 \setminus  \{ p_1, \dots p_k\}$, Galois with group $G'$.

By the Riemann existence theorem, and since $\pi_1 (\PP^1 \setminus 
\{ p_1, \dots
p_k\}) $ is the socalled infinite polygonal group $ T (\infty^k) =  T 
(\infty, \dots, \infty)
$  generated by
  simple geometric loops
$\alpha_1,\dots , \alpha_k$, satisfying the relation
$\alpha_1 \cdot  \dots  \cdot \alpha_k = 1$,
the datum of such a covering amounts to the datum of an epimorphism
$\phi : T (\infty, \dots, \infty)  \ra G' $ such that, for each $ 
i=1, \dots ,k$,
$  a_i : = \phi (\alpha_i) $ is an element of order $m_i$.

The group  $ T (\infty^k)$ is trivial for $k=1$,  infinite cyclic
for $ k=2$, in general a free group of rank $k-1$.

Since $  a_i : = \phi (\alpha_i) $ is an element of order $m_i$,
the epimorphism factors through the polygonal group
$$ T (m_1, \dots, m_k) : = \langle \alpha_1,   \dots  , \alpha_k| 
\alpha_1 \cdot  \dots \cdot
\alpha_k =  \alpha_1^{m_1}  =  \dots = \alpha_k^{m_k}  =1\rangle.$$

If $k=2$, then  we may assume  $ m_1 = m_2 = m$ and we have a cyclic
subgroup $G'$ of order $m$ of $\PP SU (2, \C)$, which, up to 
conjugation, is generated
by a transformation $ \zeta (u,v) :=  ( \zeta u,  \zeta ^{n-1}v)$,
with $\zeta$ a primitive m-th  root of $1$ for $m$ odd, and
a primitive 2m-th  root of $1$ for $m$ even. Thus, our group $G$ is
  a cyclic group of order $n$, with $ n=2 m$ for $m$ even,
and with $ n=2 m$ or  $ n=m$
for $m$  odd. $G$ is generated
by a transformation $ \zeta (u,v) :=  ( \zeta u,  \zeta ^{n-1}v)$
(with $\zeta$ a primitive n-th  root of $1$), and we have the singularity $A_n$
previously considered.

If $k=3$, the only numerical solutions for the Hurwitz' formula
are
$$ m_1 = 2, m_2 = 2, m_3 = m \geq 2,$$
$$m_1 = 2, m_2 = 3, m_3 = 3,4,5.$$

Accordingly the order of the group $G'$ equals $2m, 12, 24, 60$.
Since $m_3$, for $m_3 \geq 3$, is not the least common multiple of $m_1, m_2$, the
group $G'$ is not abelian, and it follows (compare \cite{klein})
that $ G'$ is respectively isomorphic to $  D_m, \A_4, \SSS_4 , \A_5$.

Accordingly, since as above the lift of an element in $G'$
of even order $k$ has necessarily order $2k$, it follows
that $G$ is the full inverse image of  $G'$, and $G$ is respectively
called the binary dihedral group, the binary tetrahedral group, the binary
octahedral group, the binary icosahedral group.

Felix Klein computed explicitly the ring of polynomial invariants
for the action of $G$, showing that $ \C [u,v]^G$ is
a quotient ring $ \C [x,y,z] /( z^2 - f(x,y))$ , where
\begin{itemize}
\item
$  f(x,y)= x^2 + y ^{n+1}  $ for the $A_n$ case
\item
$  f(x,y)= y (x^2 + y ^{n-2} ) $ for the $D_n$ case ($n \geq 4$)
\item
$  f(x,y)= x^3 + y ^{4}  $ for the $E_6$ case, when $G' \cong \A_4$
\item
$  f(x,y)= y (x^2 + y ^{3} ) $ for the $E_7$ case, when $G' \cong \SSS_4$
\item
$  f(x,y)=x^3 + y ^{5} $ for the $E_8$ case, when $G' \cong \A_5$.
\end{itemize}

We refer to \cite{durfee} for several equivalent characterizations
of Rational Double points, another name for the Kleinian singularities.
An important property ( cf. \cite{reid1} and \cite{reid2}) is that these singularities 
may be resolved
just by a sequence of point blow ups: in this procedure no points
of higher multiplicity than $2$ appear, whence it follows
once more that the canonical divisor of the minimal resolution
is the pull back of the canonical divisor of the singularity.

A simpler way to resolve these singularities (compare \cite{bpv}, pages
86 and following)
is to observe that they are expressed as double covers branched over
the curve $ f(x,y) = 0$. Then the standard method, explained
in full generality by Horikawa in \cite{quintics} is to resolve the branch curve
by point blow ups, and keeping as new branch curve at each step $B'' - 2D''$,
where $B''$ is the total transform of the previous branch curve $B$,
and $ D''$ is the maximal effective divisor such that $B'' - 2D''$
is also effective. One obtains the following

\begin{teo}
The minimal resolution of a Rational Double Point has as exceptional
divisor a finite union of curves $E_i \cong \PP^1$, with selfintersection $-2$,
intersecting pairwise transversally in at most one point, and moreover such
that no three curves pass through one point. The dual graph of the
singularity, whose vertices correspond to the components $E_i$,
and whose edges connect $E_i$ and $E_j$ exactly when
$ E_i \cdot E_j = 1$, is a tree, which is a linear tree with $n-1$ vertices
exactly in the $A_n$ case. In this way one obtains exactly all the
  Dynkin diagrams corresponding to the simple Lie algebras.
\end{teo}
\begin{oss}
i) See the forthcoming theorem \ref{dynkin} for a list of these Dynkin diagrams.

ii) The relation to simple Lie algebras was clarified by Brieskorn in \cite{nice}:
these singularities are obtained by intersecting the orbits of the coadjoint
action with a  three dimensional submanifold in general position.
\end{oss}

We end this subsection with an important observation concerning 
the automorphisms of
a Rational Double Point $ (X, x_0)$.

Let $H$ be a finite group of automorphisms of the germ $ (X, x_0) = 
(\C^2,0) / G$.

Then the quotient $ (X, x_0) /H$ is a quotient of $ (\C^2,0)$ by a group
$H'$ such that $ H' / G \cong H$. Moreover, by the usual averaging trick
(Cartan's lemma, see \cite{cartan}) we may assume that $ H ' \subset
GL(2, \C)$. Therefore $H'$ is contained in the normalizer $N_G$ of $G$
inside $GL(2, \C)$. Obviously, $N_G$ contains the centre $\C^*$ of
$GL(2, \C)$, and $\C^*$ acts on the graded ring $\C[x,y,z]/ (z^2 - f(x,y))$
by multiplying homogeneous elements of degree $d$ by $ t^d$.
Therefore $H$ is a finite subgroup of the group $H^*$ of graded automorphisms
of the ring $\C[x,y,z]/ (z^2 - f(x,y))$, which is determined as follows
  (compare \cite{cat2})
\begin{teo}\label{autrdp}
The group $H^*$ of graded automorphisms of a RDP is:

1) $\C^*$ for $E_8, E_7$

2) $\C^* \times \Z/2$ for $E_6, D_n  (n \geq 5 )$

3)  $\C^* \times \SSS_3$ for $ D_4$

4) $ (\C^* )^2 \times \Z/2$ for $A_n  (n \geq 2 )$

5) $GL(2, \C) / \{\pm 1\}$ for $A_1$.
\end{teo}
{\em Idea of proof .}
The case of $A_1$ is clear because $G = \{ \pm 1\}$ is contained in the centre.
In all the other cases, except $D_4$, $y$ is the generator of smallest
degree, therefore it is an eigenvector, and, up to using $\C^*$, we 
may assume that
$y$ is left invariant by an automorphism $h$.
Some calculations allow to conclude that $h$ is the identity in case 1), or
the {\bf trivial involution} $ z \mapsto -z$ in case of $E_6$ and of $D_n $
for $n$ odd; while for $D_n $
with $n$ even the extra involution is $ y \mapsto -y$.

Finally, for $D_4$, write the equation as $ z^2 = y (x+iy)(x - i y)$ 
and permute
the three lines which are the  components of the branch locus.
For $A_n$, one finds that the normalizer is the semidirect product of the
diagonal torus with the involution given
  by  $ (u,v) \mapsto (v,u)$.

One may also derive the result from  the symmetries of the Dynkin diagram.

\qed

\subsection{Canonical models of surfaces of general type.}

Assume now that $S$ is a smooth minimal (projective) surface of general
type.

We have (as an easy consequence of the Riemann Roch theorem)
 that $S$ is minimal of
general type if $K_S^2 > 0$ and
$K_S$
  is nef (we recall that a divisor $D$ is said to be 
{\bf nef} if, for each irreducible curve $C$, we have $ D \cdot C \geq 0$).

In fact, $S$ is minimal of
general type iff $K_S^2 > 0$ and
$K_S$
  is nef. Since, if $D$ is nef and, for $ m >0$, we write $ |m D| = |M| + \Phi$ as
the sum of its movable part and its fixed part, then $ M^2  = m^2 D^2  - m 
D \cdot \Phi  - M \cdot \Phi \leq m^2 D^2 $. Hence, if $D^2 \leq 0$,
the linear system  $ |m D|$ yields a rational map whose image has
dimension at most 1.

Recall further  that the Neron-Severi group
$NS(S) = Div(S) / \sim$ is the group of divisors
modulo numerical equivalence ($D$ is numerically equivalent to $0$,
and we write $D \sim 0$, $\Leftrightarrow$ $ D \cdot C = 0$
for every irreducible curve $C$ on $S$).

The Neron Severi group is a discrete subgroup of the  vector space
$ H^1 ( \Omega^1_S)$, and indeed on a projective manifold $Y$
it equals the intersection $(H^2(Y, \Z) / {\it Torsion} )\cap H^{1,1} (Y).$

By definition, the intersection form is non degenerate on the
Neron Severi group, whose rank $\rho$ is called the {\bf Picard number}.
But the Hodge index theorem implies the

{\bf Algebraic index theorem} {\em The intersection form on
$NS(S)$ has positivity index precisely $1$ if $S$ is an algebraic surface.}

The criterion of Nakai-Moishezon  says that a
divisor $L$ on a surface $S$ is ample if and only if $L^2 > 0$ and
$L\cdot C > 0$ for each irreducible curve $C$ on $S$. Hence:

{\em  The canonical divisor $K_S$ of a minimal surface of general type $S$ 
 is ample iff
there does not exist an irreducible curve $C$ ($\neq 0$) on $S$ with $K\cdot
C=0$.}

\begin{oss} Let $S$ be a minimal surface of general type and
$C$ an irreducible curve on $S$ with $K\cdot C=0$. Then, by the
  index theorem,  $C^2 < 0$ and by the adjunction
formula we see that $2p(C)-2=K\cdot C+C^2=C^2 < 0$.

In general $p(C): = 1 - \chi (\hol_C)$ is the {\bf arithmetic
genus} of $C$, which is equal to the sum $ p(C) = g (\tilde{C}) + \delta$ of the {\bf
geometric genus} of $C$, i.e., the genus of the normalization $p : \tilde{C} \ra C$ of $C$,
with the number $\delta$ of double points of $C$, defined as
$\delta : =  h^0 ( p_* \hol_{\tilde{C}} / \hol_C)$.

Therefore here $p(C)=0$,
so that $C\cong \PP^1$, and $C^2=-2$.

These curves are called $(-2)$-{\bf
curves}.
\end{oss}

Thus $K_S$ is not ample  if and only if there
exists a $(-2)$-curve on $S$. There is an upper bound for the number of these
$(-2)$-curves.

\begin{lem}  Let $C_1,\ldots,C_k$ be
irreducible $(-2)$-curves on a minimal surface $S$ of general type.
We have:
$$
(\Sigma n_iC_i)^2\le 0,
$$
and
$$
(\Sigma
n_iC_i)^2=0\,  \mbox{\it if and only if $n_i=0$ for all
$i$}.
$$
Thus their images in the Neron-Severi group $NS(S)$ are independent
and in particular $k\le \rho-1$ ( $\rho$ is the rank of
$NS(S)$), and $k\le h^1(\Omega^1_S)-1.$
\end{lem}

\proof Let $\Sigma n_iC_i=C^+-C^-$, ($C^+$ and
$C^-$ being effective divisors  without common components) be the
(unique) decomposition of $\Sigma n_iC_i$ in its positive and its
negative part. Then  $K\cdot C^+= K\cdot C^-=0$ and $C^+\cdot C^-\ge
0$, whence $(C^+-C^-)^2=(C^+)^2+(C^-)^2-2(C^+\cdot C^-)\le (C^+)^2 +
(C^-)^2$. By the index theorem  $ (C^+)^2 +
(C^-)^2 $ is $ \le 0$ and $=0$ iff $C^+=C^-=0$.
\qed

We can classify all possible configurations of $(-2)$-curves on
a minimal surface $S$ of general type by the following argument.

If $C_1$ and $C_2$ are two
$(-2)$-curves on $S$, then:
$$
0 > (C_1+C_2)^2=-4+2C_1\cdot
C_2,
$$
hence $C_1.C_2 \leq 1$, i.e., $C_1$ and $C_2$ intersect
transversally in at most one point.

If $C_1,C_2,C_3$ are $(-2)$-curves on $S$, then again
we have
$$
0 > (C_1+C_2+C_3)^2=2(-3+C_1\cdot C_2+C_1\cdot
C_3+C_2\cdot C_3).
$$
Therefore no three curves meet in one point, nor do they form a triangle.

We associate to a  configuration $\cup C_i$ of $(-2)$-curves on $S$ its
{\it Dynkin graph}: the vertices correspond to the $(-2)$-curves
$C_i$, and two vertices (corresponding to $C_i$, $C_j$) are connected
by an edge if and only if $C_i\cdot C_j=1$.

Obviously the Dynkin
graph of a configuration $\cup C_i$ is connected iff $\cup C_i$ is
connected. So, let us assume that
$\cup C_i$ is connected.

\begin{teo}\label{dynkin}  Let $S$ be a minimal
surface of general type and $\cup C_i$ a (connected) configuration of
$(-2)$-curves on $S$. Then the associated (dual) Dynkin graph of $\cup C_i$
is one of those listed in figure $1$.
\end{teo}

\begin{figure}[htbp]
\begin{center}
\input{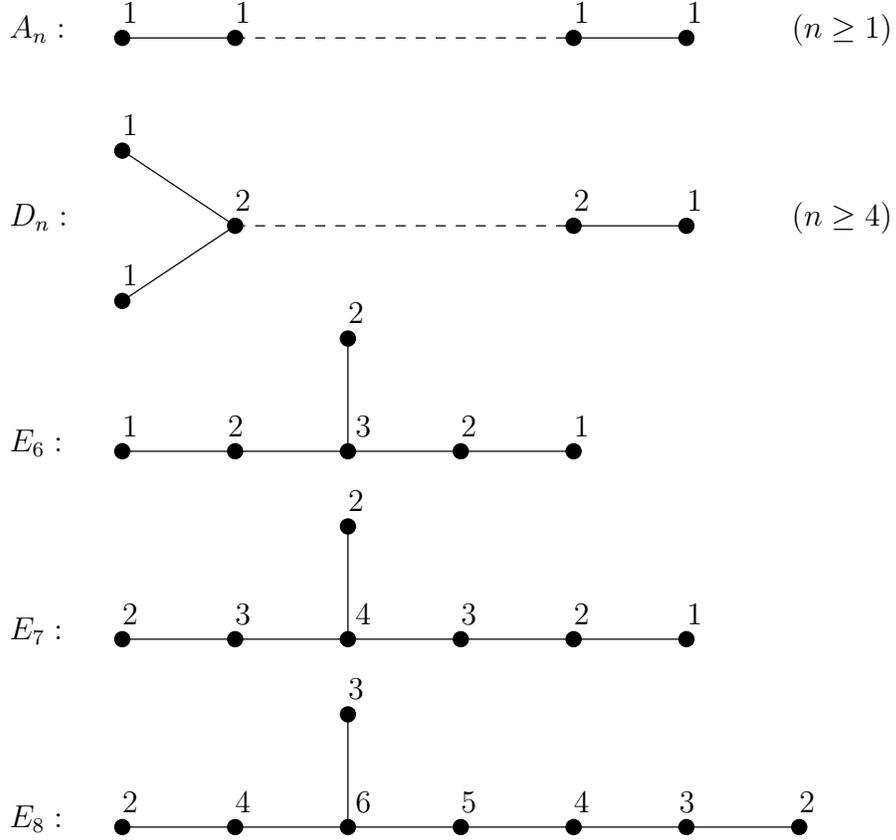}
\end{center}
\caption{The Dynkin-Diagrams of $(-2)$-curves
configurations
(the index $n$ stands for the number of vertices,
i.e. of curves).
The labels for the vertices are the coefficients of the 
fundamental cycle.}
\label{figure:yourreference}
\end{figure}

\begin{oss}
The figure indicates also the weights $n_i$  of the vertices of the 
respective trees.
These weights correspond to a divisor, called
{\it fundamental cycle}
$$
Z:=\Sigma n_iC_i
$$
  defined (cf.\cite{artin}) by the
properties
$$
Z\cdot C_i\le 0 \ \mbox{\rm for all \ $i$,  \ $Z^2=-2$, and
$n_i > 0$.}\leqno(**)
$$
\end{oss}

{\em Idea of proof of \ref{dynkin}}. The simplest proof is obtained
considering the above set of {\em Dynkin-Diagrams}
$\cD:=\{A_n,D_n,E_6,E_7,E_8\}$ and  the corresponding set of {\em
Extended-Dynkin-Diagrams} $\tilde{\cD}:=\{\tilde{A}_n,
\tilde{D}_n,\tilde{E}_6,\tilde{E}_7,\tilde{E}_8\}$ which
classify the {\bf divisors of elliptic type} 
 made of
$(-2)$-curves and are listed in Figure 2 (note that 
the divisors of elliptic type classify all
the possible nonmultiple fibres $F$ of
elliptic fibrations).   Notice that each graph
$\Gamma$ in $\cD$ is a subgraph of a corresponding graph
$\tilde{\Gamma}$ in $\tilde{\cD}$, obtained by adding exactly
a $(-2)$-curve:
$\Gamma=\tilde{\Gamma}-C_{\mbox{\scriptsize end}}$.
In this correspondence the fundamental cycle equals
$Z=F-C_{\mbox{\scriptsize end}}$ thus
(**) is proven since $ F \cdot C_i = 0 $ for each $i$. Moreover,
  by Zariski's Lemma (\cite{bpv}  ) the intersection form on $\Gamma$ is
negative definite. If moreover $\Gamma$ is a graph with a negative
definite intersection form, then $\Gamma$ does not contain
as a subgraph a graph in
$\tilde{\cD}$,
since  $F^2=0$.
The proof can now be  easily concluded.
\qed

\begin{figure}[htbp]
\begin{center}
\input{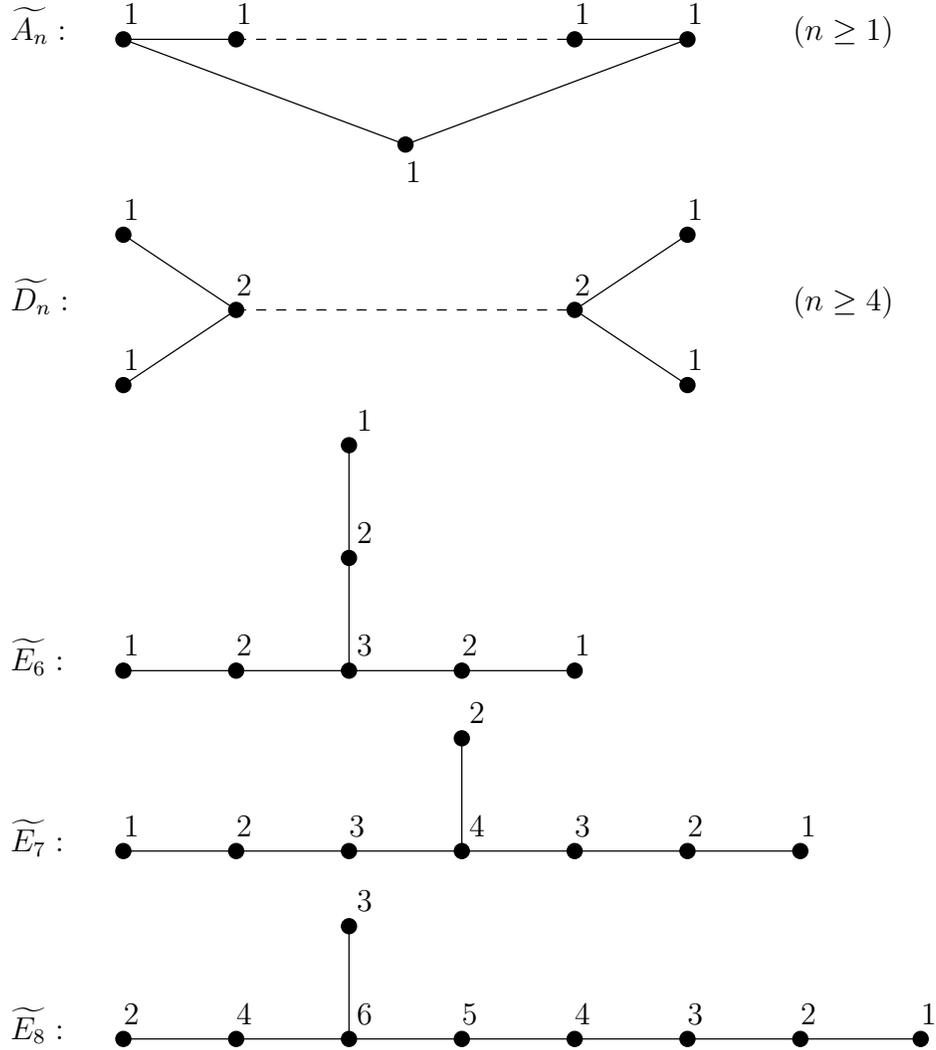}
\end{center}
\caption{The extended Dynkin-Diagrams of $(-2)$-curves
configurations.
The labels for the vertices are the coefficients of the 
divisor $F$ of elliptic type.}
\label{figure:yourreference}
\end{figure}

Artin (\cite{artin}) showed indeed that the above configurations
can be holomorphically contracted to Rational Double Points, and that
the fundamental cycle is indeed the inverse image of the maximal ideal
in the local ring of the singularity.
By applying these contractions to the minimal model $S$ of a surface of general type 
one obtains in this way a normal surface
$X$ with Rational Double Points as singularities, 
called the {\bf canonical model of $S$.}

We prefer however to sketch briefly how the canonical model
is more directly obtained from the  pluricanonical maps of $S$,
and ultimately it can be defined as the Projective Spectrum (set of homogeneous
prime ideals) of the canonical ring $\RRR (S)$.
We need first of all Franchetta's theory of numerical connectedness.

\begin{df}
An effective divisor $D$ is said to be $m$-connected if, each time
we write $D=A+B$, with
$A,B > 0$, then
$$
A\cdot B\ge m.\leqno(*)
$$
\end{df}

\begin{lem}\label{1-conn}
Let $D$ be a nef divisor on a smooth surface $S$,
with $D^2 > 0$. Then, if $D$ is effective, then $D$ is $1$-connected.
\end{lem}

\proof
Since $D$ is nef,
$$
A^2 + A \cdot B = D\cdot A \ge 0,  B^2 + A \cdot B =D\cdot B\ge 0.$$
Assume $A\cdot B \le 0$: then  $ A^2,B^2\ge -(AB) \ge 0 $
$\Longrightarrow\,
A^2\cdot B^2\ge (AB)^2$.

But, by the Index Theorem, $A^2B^2\le
(AB)^2$. Thus  equality holds in the Index theorem
$\Longleftrightarrow\exists L$ such that
$A\sim aL$, $B\sim bL$,
$D\sim (a+b)L$. Moreover, since $D^2 >0$ we have $L^2\ge 1$,
and we may
assume
$a,b > 0$ since $A,B$ are effective. Thus
$A\cdot B=a\cdot b\, \, L^2\ge 1$, equality
holding
$$
\Longleftrightarrow a=b=1 (\Longrightarrow
\, D\sim 2L),\\
    L^2=1.
$$
\qed

\begin{oss}
Let $A\cdot B=1$ and assume $A^2B^2 < (AB)^2\, \, \Longrightarrow\,
\, A^2\cdot B^2\le 0$, but $A^2,B^2\ge -1$. Thus, up to
exchanging $A$ and $B$, either $A^2=0$, and
then $D\cdot A=1, A^2=0$; or $A^2 > 0, \, B^2=-1$, and  then $D\cdot B=0$,
$B^2=-1$.
\end{oss}
Hence the following
\begin{cor}\label{2-conn} Let $S$ be minimal of general type, $D\sim mK$, $m\ge 1$:
then
$D$ is $2$-connected except possibly if $K^2=1$, and  $m=2$, or $m=1$, and $K\sim
2L$, $L^2=1$.
\end{cor}
Working a little more one finds

\begin{prop}  Let $K$ be
nef and big as before, $D\sim mK$ with $m\ge 2$. Then $D$ is
$3$-connected except possibly if
\begin{itemize}
\item
$ D=A+B, A^2=-2,\,
A\cdot K=0\, \, (\Longrightarrow A\cdot
B=2)$
\item
$m=2, K^2=1,2$
\item
$ m=3,
  K^2=1.$
\end{itemize}
\end{prop}

We use now the {\bf Curve embedding
Lemma} of \cite{cf}, improved in
\cite{4auth} to the more general case of any curve $ C$ (i.e., a pure
1-dimensional scheme) .

\begin{lem} (Curve-embedding lemma)  Let $C$ be a curve contained in a smooth
algebraic surface $S$, and let $H$ be a divisor on $C$. 
Then $H$ is very ample if, for each
length
$2$
$0$-dimensional subscheme $\zeta$ of $C$ and  for each
effective divisor $\, B\le C$,
we have
$$
{\rm Hom} (\I_\zeta,\omega_B(-H))=0.
$$

In particular $H$ is very ample on $C$ if $\forall \ B \leq C$, $ H\cdot B > 2 
p(B)-2+$ length $ \zeta \cap B $, where length $ \zeta \cap B : = $colength
$(\I_\zeta\hol_B)$. A fortiori, $H$ is very ample on $C$ if, $\forall \ B \leq C$,
$$
H\cdot B\ge
2 p(B)+1.\eqno(*)
$$
\end{lem}

\proof It suffices to show the surjectivity $H^0(\hol_C(H))\,
-\!\!\!>\!> H^0(\hol_\zeta(H))$ $\forall$ such $\zeta$.
In fact, we
can take either $\zeta=\{x,y\}$ [2 diff.\ points], or $\zeta=(x,\xi)$,
$\xi$ a tangent \ vector at $x$. The surjectivity is implied by
$H^1(\I_\zeta\hol_C(H))=0$.

By Serre-Grothendieck duality, and
since $\omega_C=\hol_C(K_S+C)$, we have, in case of nonvanishing,
$0 \neq H^1(\I_\zeta\hol_C(H))^{\vee}\cong
{\rm Hom} (\I_\zeta\hol_C(H),\hol_C(K_S+C))\ni \sigma\neq 0$.

Let $Z $ be the maximal
subdivisor of $C$ such that $\sigma$ vanishes on $Z$ (i. e., $Z=\mbox{\rm
div}(z)$, with $z|\sigma$) and  let $B=V(\mbox{\rm
Ann}(\sigma))$. Then $B+Z=C$ since, if
$C=\{(\beta\cdot z)=0\}$, Ann$(\sigma)= (\beta )$.

\medskip
Indeed, let $f\in
\I_\zeta$ be a non zero divisor: then $\sigma$ is identified with
the rational function $\sigma=\frac{\sigma(f)}{f}$; we can lift
everything to the local ring $\hol_S$, then $f$ is coprime with the
equation $\gamma : =(\beta z)$ of $C$, and $ z =\mbox{\rm
G.C.D.}(\sigma(f),\gamma)$. Clearly now Ann$(\sigma)=\{u\mid
u\sigma(f)\in(\beta z)\}=(\beta)$.

\medskip
Hence $\sigma$ induces

$$
{\hat{\sigma}} := \frac{\sigma}{z}:\I_\zeta\hol_B(H) \to \hol_B(K_S+C-Z)\\
$$
which is ``good'' (i.e., it is injective and with finite cokernel),
thus we get
$$
0\to
\I_\zeta\hol_B\stackrel{\hat{\sigma}}{\longrightarrow}\hol_B(K_B-H)\to
\Delta \ra 0
$$
where
supp$(\Delta)$ has $\dim=0$.\\
Then, taking the Euler Poincar\'e characteristics $\chi$
of the sheaves in question, we obtain

$$
0 \leq  \chi (\Delta) =\chi(\hol_B(K_B-H))- \chi(\I_\zeta\hol_B)
= - H \cdot B + 2 p(B) - 2 + length (\zeta \cap B) < 0,
$$
a contradiction.
\qed

The basic-strategy for the study of pluricanonical
maps is then to find, for every length 2 subscheme of $S$,
a divisor $C\in |(m-2)K| $ such that $\zeta\subset C$.

Since then,  in
characteristic $=0$ we have the vanishing theorem (\cite{ram})

\begin{teo} (Kodaira, Mumford, Ramanujam).
Let $L^2> 0$ on $S$, $L$ nef $\Longrightarrow H^i(-L)=0, i\ge
1$. In particular, if $S$ is minimal of general type, then 
$H^1( - K_S)= H^1(2K_S)=0$.
\end{teo}

  As shown by Ekedahl (\cite{ekedahl}) this vanishing theorem
is false in positive characteristic, but only if $\mbox{\rm char}=2$,
and for 2 very special cases of surfaces!

\begin{cor} If $C\in |(m-2)K| $, then $H^0(\hol_S(mK))\,
-\!\!\!>\!> H^0(\hol_C(mK))$. Therefore, $|mK_S|$ is very ample on $S$ if
$h^0((m-2)K_S)\ge 3$ and if the hypothesis on $H = mK_S$ in the curve embedding
Lemma is verified for any $C\equiv (m-2)K$.
\end{cor}

We shall limit ourselves here to give the proof of a weaker
version of Bombieri's theorem (\cite{bom})

{\bf Theorem on Pluricanonical-Embeddings.
(Bombieri)}. $(mK)$ is almost very ample (it embeds $\zeta$ except if
$\exists\, B$ with $\zeta\subset B$, and $B\cdot K=0$) if $m\ge 5$, $m=4$
and $K^2\ge 2$, $m=3$, $p_g\ge3$, $K^2\ge 3$.

One first sees when $h^0((m-2)K)\ge 2$.

\begin{lem} For $ m \geq 3$ we have $h^0((m-2)K)\ge 3$ except if
$m=3$ $p_g\le 2$, $m=4$, $\chi=K^2=1$ (then $q=p_g=0$) and $\ge 2$
except if $m=3$, $p_g\le 1$.
\end{lem}

\proof $p_g=H^0(K)$, so let
us assume $m\ge 4$.
$$
h^0((m-2)K)\ge\chi((m-2)K)\ge
\chi+\frac{(m-2)(m-3)}{2}K^2.
$$
Now, $\chi\ge 1$ and $K^2\ge 1$, so
we are done unless $m=4$, $\chi=K^2=1$.

\qed

The possibility that
$K_S$ may not be ample is contemplated in the following

\begin{lem} Let $H=mK,\, B\le C\equiv
(m-2)K$ and assume $ K \cdot B > 0$.
Then
$$
H\cdot B\ge 2 p(B)+1\, \, \mbox{\rm  except possibly if
}
$$

(A)\, \, $m=4$ and $K^2=1$, or $m=3$ and $K^2\le
2$.
\end{lem}

\proof Let $ C = B + Z$ as above. Then we want
$$
mK\cdot B \ 
{\ge} 2 p(B) -2 + 3= ( K + B) \cdot B + 3 =
(K+C-Z)\cdot B + 3 =
    [(m-1)K-Z] \cdot B + 3,
$$
i.e., $$
K\cdot
B+B\cdot Z \geq 3.$$
Since we assumed $K\cdot B \geq 1$, if $Z=0$ we use $ K^2 \geq 2 $ 
if $ m \geq 4$, and $ K^2 \geq 3 $ if $m=3$, else
it suffices to have
$B\cdot Z \geq 2,$ which is implied by the previous corollary \ref{2-conn}
( if $m=3$, $ B \sim Z \sim L,  \ L^2 = 1$, then   $ K \cdot B = 2$).

\qed

{\it Remark:} Note that then $\zeta$ is
contracted iff
$\exists\, B$ with $\zeta \subset B$, $K\cdot B=0$!
Thus, if there are no $(-2)$ curves, the theorem says that we have
an embedding of $S$. Else, we have a birational morphism which exactly
contracts the fundamental cycles $Z$  of $S$.
To obtain the best  technical result one has to replace
the subscheme $\zeta$ by  the subscheme $2Z$,
and use that a fundamental cycle $Z$ is
$1$-connected. We will not do it here,
we simply refer to \cite{4auth}.

The following is the more precise  theorem of Bombieri (\cite{bom})

\begin{teo}
  Let $S$ be a minimal surface of general type,
and consider the linear system $|mK|$ for $m\ge 5$, for $m=4$
when $K^2\ge 2$, for $m=3$ when $p_g\ge3$, $K^2\ge 3$.

Then $|mK|$ yields a birational morphism onto its image,
which is a normal surface $X$ with at most Rational Double Points
as singularities.  For each singular point $ p \in X$
  the inverse image of the maximal ideal $\frak M_p \subset \hol_{X,p}$
is a fundamental cycle.
\end{teo}

Here we sketch another way to look at the above surface $X$ 
(called canonical model of
$S$).

\begin{prop}   If $S$ is a surface of general type the canonical ring
$\RRR(S)$ is a graded $\C$-algebra of finite
type.
\end{prop}

\proof We choose a natural number such that $|mK|$
is without base points, and consider a pluricanonical
morphism which is birational onto its image
$$
\phi_m :S \ra
\Sigma_m=\Sigma\subset \PP^N.
$$
For $r=0,\ldots,m-1$, we set
$\FFF_r:=\phi_*(\hol_S(rK))$. \\
The Serre correspondence (cf. \cite{fac})
associates to $\FFF_r$ the
module
$$
M_r :=\bigoplus\limits^\infty_{i=1}H^0(\FFF_r(i))=
\bigoplus\limits^\infty_{i=1}H^0(\phi_*(\hol_S(rK))(i)) =$$
$$= \bigoplus\limits^\infty_{i=1}H^0(\phi_*(\hol_S((r+im)K)))=
\bigoplus\limits^\infty_{i=1}H^0(\hol_S((r+im)K)) =
\bigoplus\limits^\infty_{i=1}\RRR_{r+\mbox{\scriptsize
im}}.
$$

$M_r$ is finitely generated over the ring
$\A=\C [y_0,\ldots,y_N]$, hence $\RRR=\bigoplus\limits^{m-1}_{r=0}M_r$ is
a finitely generated $\A$-module.

We consider the natural
morphism $\alpha:\A\to \RRR, \, y_i\mapsto s_i\in \RRR_m$, (then the
$s_i$ generate a subring $B$ of $\RRR$ which is a quotient of $A$). If
$v_1,\ldots,v_k$ generate $\RRR$ as a graded $A$-module, then
$v_1,\ldots,v_k$, $s_0,\ldots,s_N$ generate $\RRR$ as a
$\C$-algebra.

\qed

The relation between the {\it
canonical ring} $\RRR(S,K_S)$ and the image of pluricanonical maps
for $ m \geq 5$ is then that
  $X=\mbox{\rm Proj}(\RRR(S,K_S))$.

In practice, since
$\RRR$ is a finitely generated graded
$\C$-algebra, generated by elements $x_i$ of degree $r_i$,
there is a surjective
morphism
$$
\lambda:\C[z_0,\ldots,z_N]-\!\!\!>\!> \RRR,\,
\lambda(z_i)=x_i.
$$
If we decree that $z_i$ has degree $r_i$, then
$\lambda$ is a graded surjective homomorphism of degree
zero.

With this grading (where
$z_i$ has degree $r_i)$
one defines (see \cite{dolg}) the {\bf weighted
  projective space} $\PP (r_0, \dots r_n)$ as
  $\mbox{\rm Proj}(\C [z_0,\ldots,z_N])$.

  $\PP (r_0, \dots r_n)$ is simply the quotient $:=\C^{N+1}-\{0\} / \C^*$, where
$\C^*$ acts on
$\C^{N+1}$ in the following
way:
$$
t(z)=(z_0 t^{r_0},\ldots, z_N t^{r_N}).
$$

The surjective homomorphism $\lambda$ corresponds to an embedding of $X$ into
$\PP (r_0, \dots r_n)$.

With the above notation, one can
easily explain some classical examples
which show that Bombieri's theorem is the best possible result.

{\bf Ex. 1}\label{ex1} $m\ge 5$ is needed.
Take a hypersurface $X_{10}\subset\PP(1,1,2,5) $
with Rational Double Points
defined by a (weighted) homogeneous polynomial $F_{10}$
of degree $10$. Then $\omega_X=\hol_X(10 -\Sigma e_i)=\hol_X(1)$,
$ K^2_X = 10 / \prod
e_i = 1$, and any $m$-canonical map with $m \leq4$
is not birational.

In fact here the quotient 
ring $
\C [y_0,y_1,x_3,z_5]/ (F_{10})$, where $\deg y_i=1,\deg x_3=2, \deg
z_5=5$ is exactly the canonical ring $\RRR(S) $.

{\bf Ex. 2:} $m=3$,
$K^2=2$ is also an exception.\\
Take $S=X_8\subset\PP(1,1,1,4)$. Here $S$
was classically described as a double cover
$S\to \PP^2$ branched on  a curve $B$ of degree $8$
( since $F_8=z^2-f_8(x_0,x_1,x_2))$.

The canonical ring, since also here $\omega_S \cong \hol_S (1)$, equals
$$
\RRR(S)=\C [x_0,x_1,x_2,z]/(F_8). $$
Thus $p_g=3$,
$K^2= 8 / 4 = 2$ but $|3K|$ factors through the double cover of $\PP^2$.

\subsection{Deformation equivalence of surfaces}

The first important consequence of the theorem on pluricanonical embeddings
is the finiteness, up to deformation, of the minimal surfaces $S$
of general type with fixed invariants $K^2$ and $\chi$.

In fact, their 5-canonical models $\Sigma_5$ are surfaces with
Rational Double Points and of degree $ 25 K^2$ in a fixed
projective space $ \PP^N$, where $ N + 1 = P_5 = h^0 (5 K_S) =
\chi + 10 K^2$.

In fact, the Hilbert polynomial of $\Sigma_5$ equals
$$ P (m) :=  h^0 (5 m K_S) =
\chi + \frac{1}{2}( 5 m -1)  5 m  K^2 .$$

Grothendieck (\cite{groth})  showed that there is

i) an integer $d$ and

ii) a subscheme $\HHH = \HHH_P$ of the
Grassmannian of codimension $P(d)$- subspaces of $ H^0 (\PP^N, \hol (d))$,
called Hilbert scheme, such that

iii) $\HHH $  parametrizes the degree $d$ pieces $ H^0
(\I_{\Sigma}(d))$ of the homogeneous ideals of all the subschemes $\Sigma 
\subset \PP^N$
having the given Hilbert polynomial $P$.

Inside $\HHH$ one has the open set
$$ \HHH ^0 : = \{ \Sigma | \Sigma {\rm \ is\ reduced \ with \ only \ 
R.D.P. 's \ as \
singularities }\}
$$

and one defines

\begin{df}The 5-pseudo moduli space of surfaces of general type with
given invariants $ K^2$, $\chi$ is the closed subset
$ \HHH _0 \subset \HHH^0$,
$$  \HHH _0 (\chi, K^2) : = \{ \Sigma \in \HHH^0| \omega_{\Sigma}^{ 
\otimes 5} \cong
  \hol_{\Sigma}(1)   \} $$

\end{df}

\begin{oss}The group $ \PP GL (N+1 , \C)$ acts on $\HHH_0$ with 
finite stabilizers
(corresponding to the groups of automorphisms of each surface)
and the orbits correspond to the isomorphism classes of minimal 
surfaces of general
type with invariants $ K^2$, $\chi$. A quotient by this action exists 
as a complex
analytic space. Gieseker showed in \cite{gieseker} that if one replaces
the 5-canonical embedding by an m-canonical embedding with 
much higher $m$, then
the corresponding quotient exists as a quasi-projective  scheme.
\end{oss}

Since $\HHH_0$ is a quasi-projective scheme, it has a finite number 
of irreducible
components (to be precise, these are the irreducible components of $( 
\HHH_0)_{red}$).

\begin{df}
The connected components of $  \HHH _0 (\chi, K^2)$ are called the 
{\bf deformation
types } of the surfaces of general type with
given invariants $ K^2$, $\chi$.
\end{df}

The above deformation types  coincide with the equivalence classes for
the relation of deformation equivalence (a more
general definition introduced by Kodaira and Spencer), in view of the following

\begin{df}
1) A {\bf deformation}  of a compact complex space $X$ is a pair
consisting of

1.1) a flat morphism $ \pi : \X \ra T$ between connected complex spaces (i.e.,
$\pi^* : \hol_{T,t} \ra \hol_{\X,x}$ is a flat ring extension for each
$ x$ with $ \pi (x) = t$)

1.2) an isomorphism $ \psi : X \cong \pi^{-1}(t_0) : = X_0$ of $X$ with a fibre
$X_0$ of $\pi$.

2) Two compact complex manifolds $X,Y$ are said to be {\bf direct deformation
equivalent} if there are a deformation $ \pi : \X \ra T$ of $X$ with $T$
 irreducible and where all the fibres are smooth,
and an isomorphism $ \psi' : Y \cong \pi^{-1}(t_1) : = X_1$ of $Y$ with a
fibre $X_1$ of $\pi$.

3) Two canonical models $X,Y$ of surfaces of general type are said to 
be {\bf direct
deformation equivalent} if there are a deformation $ \pi : \X \ra T$ 
of $X$ where $T$
is irreducible and where all the fibres have at most Rational Double Points
as singularities ,
and an isomorphism $ \psi' : Y \cong \pi^{-1}(t_1) : = X_1$ of $Y$ with a
fibre $X_1$ of $\pi$.

4) {\bf  Deformation equivalence} is the equivalence relation generated by
direct deformation equivalence.

5) A {\bf small deformation} is the germ $ \pi : (\X, X_0) \ra (T, t_0)$ of
a deformation

6) Given a deformation  $ \pi : \X \ra T$ and a morphism $ f : T' \ra T$
with $ f (t'_0) = t_0$, the {\bf pull-back} $ f^* (\X)$
is the fibre product $ \X': = \X \times_T T'$ endowed with the projection onto
the second factor $T'$ (then $ X \cong X'_0$).

\end{df}

The two definitions 2) and 3) introduced above do not conflict with each other
in view of the following

\begin{teo}\label{can=min}

Given two minimal surfaces of general type $S, S'$ and their respective
canonical models $X, X'$, then

$S$ and $S'$ are  deformation equivalent (resp.: direct deformation
equivalent) $\Leftrightarrow$
$X$ and
$X'$ are deformation equivalent (resp.: direct deformation
equivalent).
\end{teo}

We shall highlight the idea of proof of the above proposition in
the next subsection: we observe here that the proposition
implies that the deformation equivalence classes of surfaces
of general type correspond to the deformation types introduced above
(the connected components of $\HHH_0$), since over $\HHH$ lies a
natural family $ \X \ra \HHH $, $ \X \subset \PP^N \times \HHH$,
and the fibres over $\HHH^0 \supset \HHH_0$ have at most RDP's  as
singularities.

A simple but powerful observation is that, in order to analyse
deformation equivalence, one may restrict oneself to the case
where $ dim (T) = 1$: since two points in a complex space $ T \subset \C^n$
belong to the same irreducible component of $T$ if and only if they
belong to an irreducible curve $ T ' \subset T$.

One may further reduce to the case where  $T$ is smooth simply by taking the
normalization $ T^0 \ra T_{red} \ra T$ of the reduction $T_{red}$
of $T$, and taking the pull-back of the family  to $T^0$.

This procedure is particularly appropriate in order to study the 
closure of subsets
of the pseudomoduli space. But in order to show openness of certain
subsets, the optimal strategy is to consider the small deformations
of the canonical models (this is like Columbus' egg: the small
deformations of the minimal models are sometimes too complicated
to handle, as shown by Burns and Wahl \cite{b-w} already for surfaces
in $\PP^3$).

The basic tool is the generalization due to Grauert of
Kuranishi's theorem (\cite{grauert})

\begin{teo}
{\bf Grauert's Kuranishi type theorem for complex spaces.} Let $X$ 
be a compact
complex space: then

I) there is a semiuniversal deformation $ \pi : (\X, X_0) \ra (T, t_0)$ of
$X$, i.e., a deformation such that every other small deformation
$ \pi' : (\X', X'_0) \ra (T', t'_0)$ is the pull-back of $\pi$ for
an appropriate morphism $f : (T', t'_0) \ra (T, t_0)$ whose
differential at $t'_0$ is uniquely determined.

II) $(T, t_0)$ is unique up to isomorphism, and is a germ of analytic 
subspace of the
vector space
${\rm Ext }^1 (\Omega^1_X, \hol_X ),$ inverse image of the origin under a local
holomorphic map (called obstruction map and denoted by $ob$) $ ob : 
{\rm Ext }^1
(\Omega^1_X,
\hol_X ) \ra {\rm Ext }^2 (\Omega^1_X, \hol_X) $ whose differential vanishes
at the origin (the point corresponding to the point $t_0$).
\end{teo}

The theorem of Kuranishi (\cite{kur}, \cite{kur2}) dealt with the case of compact
complex manifolds, and in this case  ${\rm Ext }^j (\Omega^1_X,
\hol_X ) \cong  H^j (X, \Theta_X)$, where $\Theta_X : = Hom ( \Omega^1_X,
\hol_X ) $ is the sheaf of holomorphic vector fields.
In this case the quadratic term in the Taylor development of
$ob$, given by  the cup product $ H^1 (X, \Theta_X) \times
  H^1 (X, \Theta_X) \ra  H^2 (X, \Theta_X)$, is easier to calculate.

\subsection{Isolated singularities,
simultaneous resolution}

The main reason in the last subsection to consider deformations of 
compact complex
spaces was the aim to have a finite dimensional base $T$ for the
semiuniversal deformation (this would not have been the case in general).

Things work in a quite parallel way if one considers germs of
isolated singularities of complex spaces $ (X, x_0)$.
The definitions are quite similar, and  there is an embedding
$\X  \ra \C^n \times T$ such that  $\pi$ is induced by the second 
projection. There
is again a completely similar general theorem by Grauert (\cite{grauert1})

\begin{teo}
{\bf Grauert's theorem for deformations of isolated singularities .}
Let  $ (X, x_0)$ be a germ of an isolated singularity of a
  complex space: then

I) there is a semiuniversal deformation $ \pi : (\X, X_0, x_0) \ra 
(\C^n,0) \times
(T, t_0)$ of
$X$, i.e., a deformation such that every other small deformation
$ \pi' : (\X', X'_0, x'_0) \ra  (\C^n,0) \times (T', t'_0)$ is the 
pull-back of $\pi$
for an appropriate morphism $f : (T', t'_0) \ra (T, t_0)$ whose
differential at $t'_0$ is uniquely determined.

II) $(T, t_0)$ is unique up to isomorphism, and is a germ of analytic 
subspace of the
vector space
${\rm Ext }^1 (\Omega^1_X, \hol_X ),$ inverse image of the origin under a local
holomorphic map (called obstruction map and denoted by $ob$) $ ob : 
{\rm Ext }^1
(\Omega^1_X,
\hol_X ) \ra {\rm Ext }^2 (\Omega^1_X, \hol_X) $ whose differential vanishes
at the origin (the point corresponding to the point $t_0$).
\end{teo}

One derives easily from the above a previous result of G. Tjurina
concerning the deformations of isolated hypersurface singularities.

For, assume that $(X,0) \subset (\C^{n+1},0)$ is the zero set of a holomorphic
function $f$, $ X = \{ z | f(z)=0\}$ and therefore, if $ f_j = 
\frac{\partial f}
{\partial z_j}$, the origin is the only point in the locus
$ \SSS = \{ z | f_j(z)=0 \  \forall j \}$.

We have then the exact sequence

$$ 0 \ra \hol_X \overset{(f_j)}{\to} \hol_X^{n+1} \ra  \Omega^1_X \ra 0$$
which yields ${\rm Ext }^j (\Omega^1_X, \hol_X
) = 0 $ for $ j \geq 2$, and
$$ {\rm Ext }^1 (\Omega^1_X, \hol_X) \cong \hol_{\C^{n+1}, 0} / (f, f_1, \dots
f_{n+1}) : = T^1 .
$$

In this case the basis of the semiuniversal deformation is just the 
vector space
$T^1$, called the {\bf Tjurina Algebra}, and one obtains the following

\begin{cor}
{\bf Tjurina's deformation.} Given $(X,0) \subset (\C^{n+1},0)$ an isolated
hypersurface singularity $ X = \{ z | f(z)=0\}$, let $g_1, \dots g_{\tau}$ be a
basis of the Tjurina Algebra $T^1  =  \hol_{\C^{n+1}, 0} / (f, f_1, \dots
f_{n+1}) $ as a complex vector space.

Then $\X \subset \C^{n+1} \times \C^{\tau}$, $ \X : = \{ z |F (z,t):=  f(z)
+ \sum_j t_j g_j (z) =0 \}$ is the semiuniversal deformation of $(X,0)$.

\end{cor}

A similar result holds more generally (with the same proof) when $X$ 
is a complete
intersection of $r$ hypersurfaces
   $ X = \{ z | \phi_1(z)= \dots = \phi_r (z) = 0\}$, and then one has 
a semiuniversal
deformation of the form $\X \subset \C^{n+1} \times \C^{\tau}$, $ \X : = \{ z |
F_i (z,t):=  \phi_i(z)
+ \sum_j t_j G_{i,j}(z) =0 , i = 1, \dots r\}$.

In both cases the singularity admits a so-called {\bf smoothing}, given by
the Milnor fibre (cf. \cite{mil})

\begin{df}
Given a hypersurface singularity $ (X,0)$, $ X = \{ z | f(z)=0\}$, the Milnor
fibre $\frak M_{\delta, \epsilon}$ is the intersection of the hypersurface
$  \{ z | f(z)= \epsilon \}$ with the ball $ \overline{B ( 0, 
\delta)}$ with centre
the origin and radius $\delta << 1$, when $ |\epsilon| << \delta $.

$\frak M : = \frak M_{\delta, \epsilon}$ is a manifold with boundary whose
diffeomorphism type is independent of  $ \epsilon , \delta $ when $ 
|\epsilon| <<
\delta << 1$.

More generally, for a complete intersection, the Milnor fibre is the 
intersection
of the ball $ \overline{B ( 0, \delta)}$ with centre
the origin and radius $\delta << 1$ with a smooth level set
$ X_{\epsilon} : = \{ z | \phi_1(z)= \epsilon_1 , \dots  \phi_r (z) = \epsilon_r \}$.

\end{df}

\begin{oss}
Milnor defined the Milnor fibre $\frak M$  in a different way, as the
  intersection of the sphere $ S ( 0, \delta)$ with centre
the origin and radius $\delta << 1$ with the set $  \{ z | f(z)= 
\eta |f(z)| \}$, for $ | \eta | = 1$.

In this way the complement $ S ( 0, \delta) \setminus X$ is fibred over $S^1$
with fibres diffeomorphic to the interiors of the Milnor fibres; 
using Morse theory
Milnor showed that $\frak M$ has the homotopy type of a bouquet
of $\mu$ spheres of dimension $n$, where $\mu$, called the {\bf Milnor number},
is defined as the dimension of the {\bf Milnor algebra }
  $M^1  =  \hol_{\C^{n+1}, 0} / (f_1, \dots
f_{n+1}) $ as a complex vector space.

\end{oss}

The Milnor algebra and the Tjurina algebra coincide in the case of a
weighted homogeneous singularity (this means that there are 
weights $m_0, \dots m_n$
such that $f$ contains only monomials $ z_0^{i_0} \dots z_n^{i_n}$ of
weighted degree $ \sum_j i_j
m_j = d$), by Euler's rule $\sum_j m_j z_j f_j = d f$.

This is the case, for instance, for the Rational Double Points,
the singularities which occur on the canonical models of surfaces
of general type.
Moreover, for these, the Milnor number $\mu$ is easily
seen to coincide with the index $i$ in the label for the singularity
(i.e., $i=n$ for an $A_n$-singularity), which in turn corresponds
to the number of vertices of the corresponding Dynkin diagram.

Therefore, by the  description we gave of the minimal resolution
of singularities of a RDP, we see that this is also homotopy
equivalent to a bouquet of $\mu$ spheres of dimension 2. This is in fact
no accident, it is just a manifestation of the fact that there
is a so-called simultaneous resolution of singularities (cf. \cite{tju},
\cite{brieskorn2}, \cite{nice})

\begin{teo}\label{simultaneous}
{\bf Simultaneous resolution according to Brieskorn and Tjurina.}
Let $T : = \C^{\mu}$ be the basis of the semiuniversal deformation
of a Rational Double Point $ (X,0)$. Then there exists a finite ramified Galois
cover $ T' \ra T$ such that the pull-back $ \X ' : = \X \times _T T'$
admits a simultaneous resolution of singularities
$ p : \SSS' \ra \X'$ (i.e., $p$ is bimeromorphic and
all the fibres of the composition $  \SSS' \ra \X' \ra T'$
are smooth and equal, for $t'_0$, to the minimal resolution
of singularities of $ (X,0)$.
\end{teo}

We shall give Tjurina' s proof for the case of $A_n$-singularities.

\Proof
Assume that we have the $A_n$-singularity
$$ \{ (x,y,z) \in \C^3 | xy = z^{n+1}  \}.$$
Then the semiuniversal deformation is given by
  $$ \X : = \{ ((x,y,z) ,(a_2, \dots a_{n+1}) ) \in \C^3 \times \C^n| xy = z^{n+1} +
a_2 z^{n-1} + \dots a_{n+1}  \} ,$$
the family corresponding to the natural deformations of the
simple cyclic covering.

We take a ramified Galois covering with group $\SSS_{n+1}$ corresponding
to the splitting polynomial of the deformed degree $n+1$ polynomial
$$\X' : =  \{ ((x,y,z), (\alpha_1, \dots \alpha_{n+1})) \in \C^3 \times \C^{n+1}|
\sum_j \alpha_j = 0, \ xy = \prod_j ( z - \alpha_j)  \} .$$
One resolves the new family $\X'$ by defining
$ \phi_i : \X' \dasharrow \PP^1$ as $$ \phi_i : = (x , \prod_{j=1}^i ( z -
\alpha_j))$$
and then taking the closure of the graph of
$ \Phi : = (\phi_1, \dots \phi_n) : \X' \dasharrow (\PP^1)^n$.

\qed

We shall consider now in detail the case of a node, i.e., an $A_1$ singularity.
This singularity and its simultaneous resolution was considered also in
the course by Seidel, and will occur once more
when dealing with Lefschetz pencils (but then in lower dimension).

\begin{ex}
Consider  a node, i.e., an $A_1$ singularity.

Here, we write $ f = z^2 - x^2 - y^2$, and the total space
of the semiuniversal deformation $\X =  \{ (x,y,z,t) | f- t =0\}
= \{ (x,y,z,t) | z^2 - x^2 - y^2 = t \}$ is smooth.
The base change $ t = w^2$ produces a quadratic nondegenerate
singularity at the origin for $\X' =  \{ (x,y,z,w) |  z^2 - x^2 - y^2 = w^2 \}\
=  \{ (x,y,z,w) |  z^2 - x^2 = y^2 + w^2 \}$.

The closure of the graph of $ \psi : = \frac{z-x}{w + i y} =  \frac{w - i y}{z+x}$
yields a so-called small resolution, replacing the origin by a curve
isomorphic to $\PP^1$.

In the Arbeitstagung of 1958 Michael Atiyah made the observation
that this procedure is nonunique, since one may also use the closure
of the rational map  $ \tilde{\psi} : = \frac{z-x}{w - i y} = 
\frac{w + i y}{z+x}$
to obtain another small resolution. An alternative
way to compare the two resolutions is to blow up the origin, getting
the big resolution (with exceptional set $ \PP^1 \times \PP^1$)
and view each of the two small resolutions as the contraction
of one of the two rulings of $ \PP^1 \times \PP^1$.

Atiyah showed in this way (see also \cite{bpv}) that the moduli space for
K3 surfaces is non Hausdorff.
\end{ex}

\begin{oss}
The first proof of theorem \ref{simultaneous} was given by G. Tjurina.
It had been observed that the Galois group $G$ of the covering $ T' \ra T$
in the above theorem is the Weyl group corresponding to the Dynkin diagram
of the singularity, defined as follows. If $\sG$ is the simple algebraic group
corresponding to the Dynkin diagram (see \cite{hum}), and $H$ is a 
Cartan subgroup,
$N_H$ its normalizer, then the Weyl group is the factor group $W: =  N_H / H$.
For example, $A_n$ corresponds to the group $ SL(n+1, \C)$, its Cartan subgroup
is the subgroup of diagonal matrices, which is normalized by the symmetric
group $\SSS_{n+1}$, and $N_H$ is here a semidirect product of $H$ with
$\SSS_{n+1}$.

As we already mentioned, E. Brieskorn (\cite{nice}) found a direct 
explanation of this
interesting phenomenon, according to a conjecture of Grothendieck. He 
proved that an
element
$ x \in \sG$ is unipotent and subregular iff the morphism
$\Psi : \sG \to H /W$, sending $x$ to the conjugacy class of its 
semisimple part
$x_s$, factors around $x$ as the composition of a submersion with
the semiuniversal deformation of the corresponding RDP singularity.
\end{oss}

With the aid of Theorem \ref{simultaneous} we can now prove that
deformation equivalence for minimal surfaces of general type
is the same as restricted deformation equivalence for their
canonical models (i.e., one allows only deformations whose fibres
have at most canonical singularities).

 {\em Idea of the Proof of Theorem \ref{can=min}. }

It suffices to observe  that

0) if we have a family $ p \colon \SSS \ra \De$
where $ \Delta \subset \C$ is the unit disk, and the fibres
are smooth surfaces,  if the central fibre is minimal of general type,
then so are all the others.

1) If we have a family $ p : \SSS \to \Delta $,
where $ \Delta \subset \C$ is the unit disk, and the fibres
are smooth minimal surfaces of general type, then their
canonical models form a flat family $\pi : \X \to \Delta$.

2) If we have a flat family $\pi : \X \to \Delta$ whose fibres
$X_t$ have at most Rational Double Points and $K_{X_t}$
is ample, then for each $t \in \Delta$ there is a ramified  covering
$ f: (\De, 0) \ra (\De, t)$ such that the pull back $ f^* \X$
admits a simultaneous resolution.

0) is a consequence of Kodaira's theorem on the stability of
$-1$-curves by deformation (see\cite{kod-1}) and of the two following facts:

i) that a minimal surface $S$ with $ K^2_S > 0$ is either of
general type, or isomorphic to $\PP^2$ or to a Segre-Hirzebruch surface
$\F_n$ ($n\neq 1$, $\F_0 \cong \PP^1 \times
\PP^1$)

ii) that $\PP^2$ is rigid (every deformation
of $\PP^2$ is a product), while $\F_n$ deforms only to $\F_m$, with $ n \equiv m
\ (mod \ 2)$. 

2) is essentially the above quoted theorem, 1) is a consequence of Bombieri's theorem,
since $p_* (\hol_{\X} ( 5  K_{\X})$ is generated by global sections and
a trivialization of this sheaf
  provides a morphism  $\phi : \X \to \De \times \PP^N$ which
induces the 5-canonical embedding on each fibre.

\qed

We end this section by describing the results of Riemenschneider (\cite{riem}) on
the semiuniversal
deformation of the quotient singularity $\frac{1}{4}(1,1)  $ 
described in Example
\ref{Riem}, and a generalization thereof.

More generally, Riemenschneider considers the singularity
$Y_{k+1}$, a quotient singularity
of the RDP (Rational Double Point) $A_{2k+1}$ $\{ uv - z^{2k+ 2} = 0 \}$ by the
involution multiplying $(u,v,z)$ by $-1$. Indeed, this
is a quotient singularity of type
$ \frac{1}{4k+4} (1, 2k + 1) $, and  the $A_{2k+1}$ singularity is the
quotient by the subgroup $ 2 \Z / (4k+4) \Z$.

We use here the more general concept of Milnor fibre of a smoothing
which the reader can find in definition \ref{Milnor}.

\begin{teo}
{\bf (Riemenschneider)} The basis of the semiuniversal deformation of
the singularity $Y_{k+1}$, quotient
of the RDP $A_{2k+1}$ by multiplication by $-1$,
consists of two smooth components $T_1$, $T_2$ intersecting
transversally. Both components yield smoothings, but
only the smoothing over $T_1$ admits a simultaneous resolution.
The Milnor fibre over $T_1$ has Milnor number $ \mu = k+1$,
the Milnor fibre over $T_2$ has Milnor number $ \mu = k$.
\end{teo}

For the sake of simplicity, we shall explicitly describe the two families
in the case $k=0$ of the quotient singularity $\frac{1}{4}(1,1)  $
  described in Example
\ref{Riem}. We use for this the two determinantal presentations
of the singularity.

1) View the singularity as  $ \C [y_0, \dots , y_4] / J$, where $J$ 
is the ideal
generated by the $2
\times 2$ minors of the matrix
$\begin{pmatrix}y_0 & y_1 & y_2& y_3\\  y_5 & y_6& y_7 & y_4 \end{pmatrix}$
and by the 3 functions $f_i :=  y_i - y_{4+i}$, for $ i= 1,2,3$
(geometrically, this amounts to viewing the rational normal curve of degree
4 as a linear section of the Segre 4-fold $\PP^1 \times \PP^3$).
We get the family $T_1$, with base $\C^3$, by changing the level sets of the
three functions $f_i$ , $ f_i (y) = t_i$, for $ t = (t_1, t_2, t_3) \in \C^3$.

2) View the singularity as  $ \C [y_0, \dots , y_4] / I$, where $I$ is
the ideal generated by the
$2 \times 2$ minors  of the matrix
$\begin{pmatrix}y_0 & y_1 & y_2\\  y_1 & y_5& y_3 \\ y_2& y_3 & y_4
\end{pmatrix}$ and by the function $f : =  y_5 - y_2$.

In this second realization  the cone over a rational normal curve of degree $4$
(in $\PP^4$)  is viewed as a linear section of
the cone over the Veronese surface.

We get the family $T_2$, with base $\C$, by changing the level set of the
function $f$ , $ y_5 - y_2 = t$, for $ t \in \C$.

We see in the latter case that the Milnor fibre is just the complement to
a smooth conic in the complex projective plane $\PP^2$, therefore
its Milnor number (equal by definition to the second Betti number)
is equal to 0. Indeed the Milnor fibre is homotopically equivalent
to the real projective plane, but this is better seen in another way
which allows a great generalization.

In fact, as we already observed, the singularities $Y_k$ are a special case
($n=2, d = k+1, a= 1$) of the following
    $${\bf Cyclic \ quotient \ singularities} \ \frac{1}{d n^2} (1,dna-1)
   = A_{dn-1} / \mu_n.$$
These are quotients of $\C^2$ by a cyclic
group of order $d n^2$ acting with the indicated characters $(1,dna-1)$,
but can also be viewed as quotients  of the Rational Double Point $ A_{dn-1}  $
of equation $ uv - z^{dn} = 0$ by the action of the group $\mu_n$ of
n-roots of unity acting in the following way:

$$\xi \in \mu_n  \ {\bf acts \  by \ :} \
   (u,v,z) \rightarrow  ( \xi u , \xi^{-1} v,\xi^{a} z).$$

This quotient action gives rise to a quotient family
$ \mathcal{X} \rightarrow \C^d$, where

$ \mathcal{X}=
   \mathcal{Y}/ \mu_n$ , $ \mathcal{Y}$ is the hypersurface in
$\C^3 \times \C^d$ of equation
$$(***) \ \  uv - z^{dn} = \Sigma_{k=0}^{d-1} t_k z ^{kn}$$ and the action of
$\mu_n$ is extended trivially on the factor $\C^d$.

We see in this way that the Milnor fibre is the quotient
of the Milnor fibre of the Rational Double Point $ A_{dn-1}  $
by a cyclic group of order $n$ acting freely.
In particular, in the case $n=2, d=1,a= 1$,
it is homotopically equivalent to the quotient
of $S^2$ by the antipodal map, and we get $\PP^2_{\R}$.

Another important observation is that  $\mathcal{Y}$, being
   a hypersurface,
is Gorenstein (this means that the canonical sheaf
   $\omega_{\mathcal{Y}}$ is
invertible). Hence,  such a quotient $ \mathcal{X}=
   \mathcal{Y}/ \mu_n$ by an action which is unramified in codimension $1$,
is (by definition) $\Q$-Gorenstein.

\begin{oss}
These smoothings were considered by Koll\'ar and Shepherd Barron
   (\cite{k-sb}, 3.7-3.8-3.9,  cf. also   \cite{man0}),
who pointed out their relevance in the theory
of compactifications of moduli spaces of surfaces,
and showed that, conversely, any $\Q$-Gorenstein smoothing
of a quotient singularity is induced by the above family (which has a
smooth base, $\C^d$).
\end{oss}

Returning to the cyclic quotient
singularity $\frac14 (1,1)$, the first description
that we gave of the $\Q$-Gorenstein smoothing (which does
obviously not admit a simultaneous resolution since its Milnor number is 0)
makes clear that an alternative way is to view
the singularity (cf. example \ref{Riem}) as a bidouble cover
of the plane branched on three lines
passing through the origin, and then this smoothing ($T_2$) is simply 
obtained by
deforming these three lines till they meet in 3 distinct points.

\newpage

\section{Lecture 3: Deformation and diffeomorphism, canonical 
symplectic structure for
surfaces of general type}
\label{third}

Summarizing some of the facts we saw up to now, given a birational equivalence
class of surfaces of general type, this class contains a unique
(complete) smooth minimal surface $S$, called the {\bf minimal model},
such that $K_S$ is nef ($ K_S \cdot C \geq 0$ for every effective curve $C$);
and a unique surface $X$ with at most Rational Double Points as singularities,
  and such that the invertible sheaf $\omega_X$
is ample, called the {\bf canonical model}.

$S$ is the minimal resolution of the singularities of $X$,
and every pluri-canonical map of $S$ factors through the projection
$ \pi : S \ra X$.

The basic numerical invariants of the birational class are
$\chi : = \chi (\hol_S) = \chi (\hol_X) = 1 - q + p_g$
($p_g = h^0(\hol_S (K_S)) = h^0 (\omega_X)$)
and $K^2_S = K^2_X$ (here $K_X$ is a Cartier divisor
such that $\omega_X \cong \hol_X (K_X)$).

The totality of the canonical models of surfaces with fixed numerical 
invariants
$ \chi = x , K^2 = y$ are parametrized (not uniquely, because of the action of the
projective group) by a quasi  projective scheme
$  \HHH _0 (x, y)$, which we called the {\bf pseudo moduli space}.

The connected components of the pseudo moduli spaces $  \HHH _0 (x, y)$
are the deformation types of the surfaces of general type, and
a basic question is whether one can find some invariant to
distinguish these. While it is quite easy to find
invariants for the irreducible components of the pseudo moduli space,
just by using the geometry of the fibre surface over the generic point,
it is less easy to produce effective invariants for the connected components.
Up to now the most effective invariant to distinguish
connected components has been the {\bf divisibility index r} of the 
canonical class
($r$ is the divisibility of $c_1 (K_S)$ in $H^2(S, \Z)$) (cf. \cite{cat3})

Moreover, as we shall try to illustrate more amply in the next 
lecture, there is
another fundamental difference between the curve and the surface 
case. Given a curve,
the genus $g$ determines the topological type, the differentiable type,
and the deformation type, and the moduli space $\frak M _g$ is irreducible.

In the case of surfaces, the pseudo moduli space $  \HHH _0 (x, y)$ is defined
over $\Z$, whence the absolute Galois group  $ Gal (\bar {\Q}, \Q )$
operates on it. In fact, it operates by possibly changing the topology
of the surfaces considered, in particular the fundamental group may change !

Therefore the algebro-geometric study of moduli spaces cannot be reduced
only to the study of isomorphism classes of complex structures
on a fixed differentiable manifold.

We shall now recall how the deformation type determines the differentiable
type, and later we shall show that each surface of general type $S$ has a
symplectic structure $(S, \omega)$, unique up to symplectomorphism,
such that the cohomology class of $ \omega $ is the canonical class
$c_1 (K_S)$.

\subsection{Deformation implies diffeomorphism.}
Even if  well known,
let us recall the  theorem of Ehresmann (\cite{ehre})

\begin{teo}{\bf (Ehresmann)}
Let $ \pi :  \mathcal {X}\rightarrow T  $ be a proper submersion of
differentiable manifolds with $T$ connected: then $\pi$ is a 
differentiable fibre
bundle, in particular  all the fibre manifolds $X_t$ are
diffeomorphic to each other.
\end{teo}

The idea of the proof is to endow $\X$ with a Riemannian metric,
so that a local vector field $\xi$ on the base $T$ has a unique
differentiable lifting which is orthogonal to the fibres.
Then, in the case where $T$ has dimension 1, one integrates the lifted
vector field. The general case is proven by induction on $ dim_{\R}T$.

The same argument allows  a variant with boundary of
Ehresmann's theorem

\begin{lem}\label{variant}
Let $ \pi :  \sM \rightarrow T  $ be a proper submersion
where $\sM$ is a differentiable manifold with boundary, such that
also the restriction of $\pi$ to $\partial \sM$
is a submersion. Assume that $T$ is
a ball in $\R^n$, and
assume that we are given a fixed trivialization  $\psi$ of
   a  closed family $  \mathcal {N} \rightarrow T  $ of submanifolds with
boundary. Then we can find a
trivialization of $ \pi :  \mathcal {M} \rightarrow T  $ which
induces the given trivialization $\psi$.
\end{lem}

\begin{proof}
It suffices to take on $\mathcal {M}$ a Riemannian metric
where the sections $ \psi (p,T) $, for $ p \in  \mathcal {N}$,
are orthogonal to the fibres of $\pi$. Then we use the customary proof
of Ehresmann's theorem, integrating  liftings orthogonal to the fibres of
   standard vector fields on $T$.
\end{proof}

Ehresmann's theorem implies then the following

\begin{prop}
Let $X, X'$ be two compact complex manifolds which are deformation equivalent.
Then they are diffeomorphic by a diffeomorphism $\phi : X' \ra X$
preserving the canonical class (i.e., such that   $\phi^* c_1( K_X) = c_1(
K_{X'})$).
\end{prop}
\Proof
The result follows by induction once it is established for $X,X'$ fibres
of a family $ \pi: \X \ra  \De$ over a 1-dimensional disk.
Ehresmann's theorem provides a differentiable trivialization
$ \X \cong X \times \De$. Notice that, since the normal bundle to a 
fibre is trivial,
the canonical divisor of a fibre $K_{X_t}$ is the restriction
of the canonical divisor $K_{\X}$ to $X_t$. It follows that the
trivialization provides a diffeomorphism $\phi$ which
preserves the canonical class.
\qed

\begin{oss}Indeed,
by the results of Seiberg Witten theory, an arbitrary diffeomorphism
between differentiable 4-manifolds carries
   $c_1( K_X)$ either to  $c_1( K_{X'})$ or to  $- c_1( K_{X'})$
(cf. \cite{witten} or \cite{mor}). Thus deformation equivalence
imposes only $\e$ more than diffeomorphism only.
\end{oss}

\subsection{Symplectic approximations of projective varieties with
isolated singularities.}

The variant \ref{variant} of Ehresmann's theorem will now be first 
applied to the
Milnor fibres of smoothings of isolated singularities.

Let $ (X,x_0)$ be the germ of an isolated singularity of a complex space,
which is pure dimensional 
of dimension $ n = dim_{\C}X$,
assume $x_0 = 0 \in X \subset \C^{n+m}$, and consider as above
the ball $ \overline{B ( x_0, \delta)}$ with centre
the origin and radius $\delta $.
Then, for all $ 0 < \delta << 1$, the intersection $\sK_0 : =  X \cap S ( x_0,
\delta)$, called the {\bf link} of the singularity, is a smooth manifold of
real dimension $ 2n-1$. 

Consider the semiuniversal deformation $ \pi : (\X, X_0, x_0) \ra 
(\C^{n+m},0) \times
(T, t_0)$ of $X$ and
the family of singularity links
$\sK : =  \X \cap  (S ( x_0, \delta) \times
(T, t_0))$. By a uniform continuity argument it follows that
$\sK \to T$ is a trivial bundle if we restrict $T$ suitably around 
the origin $ t_0$
(it is a differentiably trivial fibre bundle in the sense of
stratified spaces, cf. \cite{math}).

We can now introduce the concept of Milnor fibres of $ (X,x_0)$.

\begin{df}\label{Milnor}
Let $(T, t_0)$ be the basis of the semiuniversal deformation of a germ
of isolated singularity  $ (X,x_0)$, and let $T = T_1 \cup \dots \cup T_r$
be the decomposition of $T$ into irreducible components.
$T_j$ is said to be a smoothing component if there is a $ t \in T_j$ 
such that the
corresponding fibre $X_t$ is smooth.  If $T_j$ is a smoothing component,
then the corresponding Milnor fibre is the intersection
of the ball $ \overline{B ( x_0, \delta)}$ with the fibre $X_t$, for
$ t \in T_j$, $ |t| < \eta << \de << 1$.
\end{df}

Whereas the singularity links form a trivial bundle, the Milnor fibres
form only a differentiable bundle of manifolds with boundary
over the open set
$ T^0_j : = \{ t \in T_j ,  |t - t_0| < \eta |  \ X_t {\rm \ is \ smooth} \}.$

Since however $T_j$ is irreducible, $ T^0_j$ is connected, and the Milnor fibre
is unique up to smooth isotopy, in particular up to diffeomorphism.

We shall now apply again lemma \ref{variant} in order to perform some
surgeries to projective varieties with isolated singularities.

\begin{teo}\label{glueing}
Let $X_0 \subset \PP^N$ be a projective variety with isolated
  singularities admitting a smoothing component.

Assume that for each singular point $x_h \in X$,
we choose a smoothing component $T_{j(h)}$ in the basis of the
semiuniversal deformation of the germ $(X, x_h)$.
Then (obtaining different results for each such choice) $X$ can be
approximated by symplectic submanifolds
$W_t$ of $ \PP^N$, which are diffeomorphic to the glueing
of the 'exterior' of $X_0$ (the complement to the union $B = \cup_h B_{h}$ of
   suitable (Milnor) balls around the singular points) with
the  Milnor fibres $\sM_h$ , glued along the singularity links $\sK_{h,0}$.
\end{teo}

\begin{figure}[htbp]
\begin{center}
\input{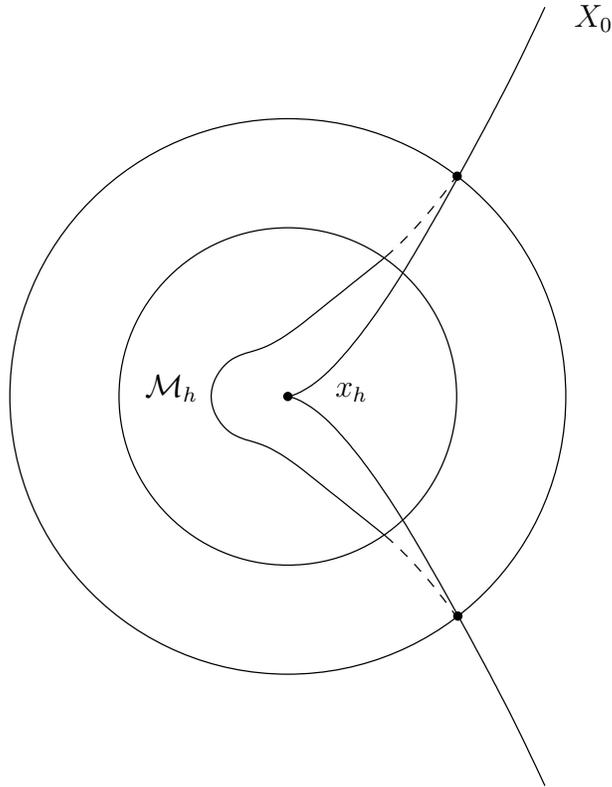}
\end{center}
\caption{Glueing the ``exterior'' of $X_0$ 
(to the Milnor Ball around $x_h$)
with a smaller Milnor fibre $\mathcal{M}_h$}
\label{figure:yourreference}
\end{figure}

A pictorial view of the proof is contained in Figure 3.

\Proof

First of all, for each singular point $x_h \in X$, we choose
a holomorphic path $ \De \to T_{j(h)}$ mapping $0$ to the
distinguished point corresponding to the germ $ (X,x_h)$,
and with image of $ \De \setminus 0$ inside  the smoothing
locus $ T^0_{j(h)} \cap \{ t | | t| < \eta \}.$

We apply then lemma \ref{variant} once more
  in order
to thicken the trivialization of the singularity links to a closed tubular
neighbourhood
in the family $\X$.

Now, in order to simplify our notation, and without loss of generality,
assume that
$X_0$ has only one singular point $x_0$, and let $B: = B(x_0, \de)$
be a Milnor ball around the singularity. Moreover, for $t \neq 0, t 
\in \De \cap B(0,
\eta)$ we consider the Milnor fibre $\MM_{\de,\eta}(t)$, whereas we have the
two Milnor links
   $$\K_0 : = X_0 \cap S(x_0, \de)  \ {\rm and } \ \K_t : =
   \X_t  \cap S(x_0, \de - \e) $$ .

We can consider the Milnor collars $\CCC_0(\e) := X_0 \cap
(\overline{B(x_0, \de)} \setminus B(x_0, \de - \e))$,
and $\CCC_t(\e) := \X_t \cap
(\overline{B(x_0, \de)} \setminus B(x_0, \de - \e))$.

The Milnor collars fill up a complex
submanifold of dimension $ dim X_0 +1 := n + 1$
of $ \C^{n+m} \times \De$.

We glue now $ X \setminus  B(x_0, \de - \e))$  and
the Milnor fibre $\MM_{\de,\eta}(t)$
by identifying the Milnor collars $\CCC_0(\e)$ and $\CCC_t(\e) $.

We obtain in this way an abstract differentiable manifold
$W$ which is independent of $t$, but we want now to give an embedding
$ W \to W_t \subset \C^{n+m}$ such that  $ X \setminus  B(x_0, \de))$
maps through the identity, and the complement of
the collar inside  the Milnor fibre
maps to  $\MM_{\de,\eta}(t)$  via the restriction of the identity.

As for the collar  $\CCC_0(\de)$, its outer boundary will be mapped to
  $\K_0 $, while its inner boundary will be mapped to $ \K_t$
(i.e., we join the two different singularity links  by a differentiable
embedding  of the  abstract Milnor collar).

For $\eta < < \de $ the tangent spaces to the image of the
abstract Milnor collar can be made very close to the tangent spaces
of the Milnor collars $\MM_{\de,\e}(t)$, and we can conclude
the proof via lemma \ref{approx}.

\qed

The following well known theorem of Moser guarantees that,
once the choice of a smoothing component is made for each $x_h \in Sing(X)$,
then the approximating symplectic submanifold $ W_t$ is unique up to
symplectomorphism.

\begin{teo}{\bf ( Moser)}
Let $ \pi :  \mathcal {X}\rightarrow T  $ be a proper submersion of
differentiable manifolds with $T$ connected, and assume that we have
a differentiable $2$-form $\omega$ on $\mathcal {X}$ with the property
that

(*) $\forall t \in T$ $ \omega_t : = \omega |_{X_t}$ yields
a symplectic structure on $X_t$ whose class in $H^2 (X_t, \R)$
   is locally constant on $T$ (e.g., if it lies on $H^2 (X_t, \Z)$).

Then the symplectic manifolds $(X_t, \omega_t)$ are all symplectomorphic.
\end{teo}

The unicity of the symplectic manifold $W_t$ will play a  crucial
role in the next subsection.

\subsection{Canonical symplectic structure for varieties with ample canonical
class and canonical symplectic structure for  surfaces of general type.}

\begin{teo}\label{canonical}A minimal surface of general type
$S$ has a canonical symplectic
   structure, unique up to symplectomorphism, and stable by deformation,
such that the class of
the symplectic form is the class of the canonical sheaf
$   \Omega^2_S = \hol_S (K_S).$
The same result holds for any projective smooth variety with
ample canonical bundle.

\end{teo}

\Proof

Let $V$ be
a smooth projective variety  of dimension $n$ whose  canonical divisor
$ K_V$ is ample.

Then there is a positive integer $m$ (depending only on $n$)  such 
that $m K_V$ is
very ample (any $m\geq 5$ does by Bombieri's theorem in the case
of surfaces, for higher dimension we can use Matsusaka's big theorem,
   cf. \cite{siu} for an effective version).

  Therefore
the $m$-th pluricanonical map $ \phi_m := \phi_{|mK_V|}$
is an embedding of $V$ in a projective space $\PP^{P_m-1}$,
where $P_m := dim H^0(\hol_V (m K_V) )$.

We define then  $\omega_m$ as follows:  $\omega_m := \frac{1}{m} \phi_m
   ^* (FS)$ (where $FS$
is the Fubini-Study form $\frac{i}{2 \pi } \partial
\overline{\partial} log |z|^2$), hence $\omega_m$ yields
a symplectic form as desired.

One needs to show that the symplectomorphism class of $(V, \omega_m)$
is independent of $m$. To this purpose, suppose that the integer $r$ has also
the property that $ \phi_r$  yields an embedding of $V$: the same 
holds also for
$rm$, hence it suffices to show that $(V, \omega_m)$ and
$(V, \omega_{mr})$ are symplectomorphic.

To this purpose we use first the well known and easy fact
that the pull back of the Fubini-Study form under the $r$-th
Veronese embedding $v_r$ equals the $r$-th multiple of the Fubini-Study
form. Second,  since $v_r \circ \phi_m$ is a linear projection
of $\phi_{rm}$, by Moser's Theorem follows the desired symplectomorphism.
Moser's theorem implies also that if we have a deformation
$ \pi : \sV \to T$ where $T$ is connected and all the fibres have 
ample canonical
divisor, then all the manifolds $V_t$, endowed with their
canonical symplectic structure, are symplectomorphic.

Assume now that $S$ is a minimal surface of general type and
that $ K_S$ is not ample: then for any $m\geq 5$ (by Bombieri's
cited theorem) $ \phi_m$  yields an embedding of the canonical model
$X$ of $S$, which is obtained by contracting the finite number of
smooth rational curves with selfintersection number $= -2$ to a finite
number of Rational Double Point singularities. For these, the
base of the semiuniversal deformation is smooth and yields a
   smoothing of the singularity.

By the quoted theorem \ref{simultaneous} on simultaneous resolution,
it follows that

1) $S$ is diffeomorphic to any smoothing $S'$ of $X$
(but it can happen that $X$ does not admit any global smoothing, as shown by
 many  examples which one can find for instance in \cite{cat5}).

2) $S$ is diffeomorphic to the manifold obtained glueing the exterior
$X \setminus B$ ($B$ being the union of Milnor balls of radius $\delta$ around
the singular points of $X$) together  with the respective Milnor fibres,
i.e., $S$ is diffeomorphic to each of the symplectic submanifolds $W$
of projective space which approximate the embedded canonical model 
$X$ according
to theorem \ref{glueing}.

We already remarked that
$W$ is unique up to symplectomorphism, and this fact ensures that we have a
unique canonical symplectic structure on $S$ (up to symplectomorphism).

  Clearly moreover, if
$X$ admits a global smoothing, we can then take $S'$ sufficiently close to $X$
   as our approximation $W$.
Then $S'$ is a surface with ample canonical bundle, and, as we have
seen, the symplectic structure induced by (a submultiple of) the Fubini Study
form is the canonical symplectic structure.

The stability by deformation is again a consequence of Moser's theorem.

\qed

\subsection{Degenerations preserving the
canonical symplectic structure.}

Assume once more that we consider the minimal surfaces $S$
of general type with fixed invariants  $\chi = x$ and $K^2 = y$,
and their 5-canonical models $\Sigma_5$, which  are surfaces with
Rational Double Points and of degree $ 25 K^2$ in a fixed
projective space $ \PP^N$, where $ N  =
\chi + 10 K^2 - 1$.

The choice of $S$ and of a projective basis for
$   \PP H^0 (5  K_S)$ yields, as we saw,
a point in the
5-pseudo moduli space of surfaces of general type with
given invariants $\chi = x$ and $K^2 = y$,
  i. e.,  the locally closed set $  \HHH _0 (x , y) $
of the corresponding Hilbert scheme $ \HHH$,
which is the closed subset
$$  \HHH _0 (x , y) : = \{ \Sigma \in \HHH^0| \omega_{\Sigma}^{ 
\otimes 5} \cong
  \hol_{\Sigma}(1)   \} $$
of the open set
$$ \HHH ^0 (x , y): = \{ \Sigma | \Sigma {\rm \ is\ reduced \ with \ 
only \ R.D.P. 's
\ as \ singularities }\}.
$$

In fact, even if this pseudo moduli space is conceptually clear, it is
computationally more complex than just an appropriate
 open subset of $ \HHH ^0 (x , y)$,
which we denote by  $ \HHH ^{00} (x , y)$ and  parametrizes triples
$$(S, L, \sB)$$
where

i) $S$ is a minimal surface
of general type with fixed invariants  $\chi = x$ and $K^2 = y$

ii) $ L \in Pic^0  (S)$ is a topologically trivial holomorphic line bundle

iii) $\sB$ is a a projective basis for
$   \PP H^0 (5  K_S + L)$.

To explain how to define $ \HHH ^{00} (x , y)$, let $ \HHH ^{n} (x , y)
\subset  \HHH ^{0} (x , y)$ be the open set of surfaces $\Sigma$
with $K^2_{\Sigma} = y$. Let $H$ be the hyperplane divisor, and observe
that by the Riemann Roch theorem $ P_{\Sigma} (m) = \chi (\hol_{\Sigma})
+ 1/2 \ m H \cdot (m H - K_{\Sigma})$, while by definition 
$ P_{\Sigma} (m) = x + 
 1/2 (5m-1) 5 m y$. Hence, $H^2 = 25 y$, $ H \cdot  K_{\Sigma} = 5 y$, 
$\chi (\hol_{\Sigma}) = x$, and by the Index theorem $K^2_{\Sigma} \leq y$,
equality holding if and only if $ H \sim 5 K_{\Sigma}$.

Since the group of linear equivalence classes of divisors which are numerically
equivalent to zero is parametrized by  $Pic^0  (\Sigma) \times Tors
(H^2 (\Sigma, \Z))$, we get that the union of the connected components of  
$ \HHH ^{n} (x , y)$ containing $ \HHH _0 (x , y)$ yields an open set
$ \HHH ^{00} (x , y)$ as described above.

Since $Pic^0  (S)$ is a complex torus of dimension $q = h^1 (\hol_S)$,
it follows that indeed there is a natural bijection, induced by inclusion,
between irreducible (resp. connected) components of $  \HHH _0 (x , y)$
and of $  \HHH ^{00} (x , y)$. Moreover, $  \HHH _0 (x , y) $ and   
$\HHH ^{00} (x , y)$
coincide when $  q= 0$.

As we shall see, there are surfaces of general type which
are diffeomorphic, or even canonically symplectomorphic,
but which are not deformation equivalent.

Even if $ \HHH ^{00} (x , y)$ is highly disconnected, and not pure dimensional,
one knows by a general result by Hartshorne (\cite{hartshilb}), that the Hilbert scheme
$\HHH$ is connected, and one may therefore
  ask

A) is $\overline{ \HHH ^{00} (x , y)}$ connected ?

B) which kind of singular surfaces does one have to consider in order to
connect different components of $ \HHH ^{00} (x , y)$?

The latter question is particular significant, since first of
all any projective variety admits a flat deformation to a scheme supported
on the projective cone over its hyperplane section (iterating this procedure,
one reduces to the socalled {\bf stick figures}, which in this case would
be supported on a finite union of planes. Second, because when going across
badly singular surfaces, then the topology can change drastically
(compare example 5.12, page 329 of \cite{k-sb}).

We refer to \cite{k-sb} and to \cite{vieh} for  a theory
of compactified moduli spaces of surfaces of general type. We would 
only like to
mention that the theory describes  certain  classes of singular surfaces
which are allowed, hence a certain open set in the Hilbert scheme $\HHH$.

One important question is, however,  which degenerations of smooth
surfaces do not change the canonical symplectomorphism class. In other
words, which surgeries do not affect the canonical symplectic structure.

A  positive result is  the following theorem, which is used
in order to show that the Manetti surfaces are canonically symplectomorphic
(cf. \cite{cat06})

\begin{teo}\label{families}
Let $  \mathcal {X}\subset \PP ^N \times \Delta$ and
$  \mathcal {X}' \subset \PP ^N \times \Delta'$ be two flat families
of normal surfaces over the disc of radius 2 in $\C$.

Denote by $ \pi :  \mathcal {X} \to  \Delta$ and
by  $ \pi ' :  \mathcal {X}' \to  \Delta$ the respective projections
and make the following assumptions on the respective fibres of $\pi , \pi'$:

1) the central fibres
$X_0 $ and $ X_0'$ are surfaces with
cyclic quotient singularities and
   the two flat families yield
    $\Q$-Gorenstein smoothings of them.

2) the other fibres $X_t$, $X_t'$, for $ t, t' \neq 0$ are smooth.

Assume moreover that

3) the central fibres
$X_0 $ and $ X_0'$ are projectively equivalent to
respective fibres ($X_0 \cong Y_0$ and $ X_0' \cong Y_1$)
   of an equisingular
projective family $  \mathcal {Y}\subset \PP^N \times \Delta$
of surfaces.

Set $X := X_1$, $X' := X_1'$: then

a) $X$ and $X'$ are diffeomorphic

b) if  $FS $ denotes the symplectic form inherited from the Fubini-Study
K\"ahler metric on $\PP^N$, then the symplectic manifolds
$(X,FS)$ and $(X',FS)$ are symplectomorphic.

\end{teo}

The proof of the above is based on quite similar ideas to those
of the proof of theorem \ref{glueing}.

\begin{oss}

Theorem \ref{families} holds
   more generally for varieties of higher dimension with isolated
singularities under the assumption
that, for each singular point $ x_0$ of $X_0$, letting
$y_0(t)$ be the corresponding singularity of  $Y_t$

i)  $ (X_0, x_0) \cong ( Y_t, y_0(t))$

ii) the two smoothings $  \mathcal {X},   \mathcal {X}'$,
correspond to paths in the same irreducible
component of $Def (X_0, x_0) $.

\end{oss}

\newpage

\section{Lecture 4: Irrational pencils, orbifold fundamental groups, and
surfaces isogenous to a product.}
\label{fourth}

In the previous lecture we considered the possible deformations
and mild degenerations of surfaces of general type. In this lecture we want to
consider a very explicit class of surfaces (and higher dimensional varieties),
those which admit an unramified covering which is a product of curves
(and are said to be {\bf isogenous to a product}).
For these one can reduce the description of the moduli space to the
description of certain moduli spaces of curves with automorphisms.

Some of these varieties are rigid, i.e., they admit no nontrivial deformations;
in any case these surfaces $S$ have the weak rigidity property  that any
surface homeomorphic to them is deformation equivalent either
to $S$ or to the conjugate surface $\bar{S}$.

Moreover, it is quite interesting to see which is the action of complex
conjugation on the  moduli space: it turns out that it interchanges
often two distinct connected components.
In other words,there are surfaces such that the complex conjugate
surface is not deformation equivalent to the surface itself
(this phenomenon has been observed by several authors independently,
cf.  \cite{f-m} Theorem 7.16 and  Corollary 7.17 on p.
208, completed in \cite{friedman} for elliptic surfaces, cf.  \cite{k-k},\cite{cat4},
\cite{bcg} for the case of surfaces of general type). However,
in this case we obtained surfaces which are diffeomorphic to each other, but only
through a diffeomorphism not preserving the canonical class.

Other reasons to include these examples are not only their simplicity
and beauty, but also the fact that these surfaces lend themself quite
naturally to reveal the action of the Galois group $ Gal (\overline{\Q}, \Q)$
on moduli spaces.

In the next section we shall recall some basic results on fibred surfaces
which are used to treat the class of surfaces isogenous to a product.

\subsection{Theorem of Castelnuovo-De Franchis, irrational pencils and
the orbifold fundamental group}

We recall some classical and some new results (see \cite{albanese} and
\cite{fibred}
for more references)

\begin{teo} {\bf Castelnuovo-de Franchis}. Let $X$ be a compact
    K\"ahler manifold and $U \subset H^0(X, \Omega^1_X)$ be an isotropic
    subspace (for the wedge product) of dimension $\geq 2$. Then there
    exists a fibration $f : X \rightarrow B$, where $B$ is  a curve,
    such that $U \subset f^* (H^0(B,\Omega^1_B))$ (in particular,
the genus  $ g (B)$ of
$B$ is at least 2).
\end{teo}

{\em Idea of proof}

Let $\omega_1, \omega_2$ be two $\C$-linearly independent 1-forms
$\in H^0(X, \Omega^1_X)$ such that
$\omega_1 \wedge \omega_2 \equiv 0$. Then their ratio defines a nonconstant
meromorphic function $F$ with $\omega_1 = F \omega_2$.

After resolving the indeterminacy of the meromorphic map $ F : X 
\dasharrow \PP^1$
we get a morphism $ \tilde{F} : \tilde{X}  \to \PP^1$ which does not 
need to have
connected fibres, so we let  $ f : \tilde{X}  \to B$ be its Stein 
factorization.

Since holomorphic forms are closed, $0 = d \omega_1 = d F \wedge  \omega_2$
and the forms $\omega_j$ restrict to zero on the fibres of $f$.
A small ramification calculation shows then that the two forms $\omega_j$
are pull back of holomorphic one forms on $B$, whence $B$ has genus
at least two. Since every map of $\PP^1 \to B$ is constant, we
see that $f$ is indeed holomorphic on $X$ itself.

\qed

\begin{df}
Such a fibration $f$ as above is called an {\em irrational pencil}.
\end{df}

Using Hodge theory and the K\"unneth formula, the Castelnuovo-de
Franchis theorem implies (see \cite{albanese}) the following

\begin{teo}(Isotropic subspace theorem).
1) Let $X$ be a compact K\"ahler manifold and
$U \subset H^1(X, \mathbb{C})$ be an isotropic subspace of
    dimension $\geq 2$. Then there exists an irrational pencil
$f : X \rightarrow
    B$, such that $U \subset f^* (H^1(B,\mathbb{C}))$.

2)There is a 1-1
    correspondence between irrational pencils $f:X \rightarrow B$, $g(B)=b
    \geq 2$, and subspaces $V = U \oplus \bar{U}$, where $U$ is maximal
    isotropic of dimension $b$.
\end{teo}

\Proof

1) Using the fact that $ H^1 (X, \C) =  H^0(X, \Omega^1_X) \oplus
\overline { H^0(X, \Omega^1_X) }$ we write a basis of $U$ as
$ ( \phi_1 = \omega_1 + \overline { \eta_1}, \dots , \phi_b = 
\omega_b + \overline {
\eta_b})$.

  Since again Hodge theory gives us the direct sum
$$ H^2 (X, \C) =  H^0(X, \Omega^2_X) \oplus H^1 (X, \Omega^1_X) \oplus
\overline { H^0(X, \Omega^2_X) }$$
the isotropicity condition $ \phi_i \wedge \phi_j = 0 \in H^2 (X, \C)$
reads:
$$ \omega_i \wedge \omega_j \equiv 0, \eta_i \wedge \eta_j \equiv 0,
  \omega_i \wedge \overline {\eta_j} + \overline{ \eta_i }\wedge 
\omega_j \equiv 0,
  \ \forall i,j.$$
The first two identities show that we are done if
one can apply the theorem of
Castelnuovo de Franchis to the $\omega_j$'s, respectively to the $\eta_j$'s ,
obtaining two irrational pencils $f : X \ra B,f' : X \ra B'$.
In fact, if the image of $ f \times f' : X \to B \times B'$ is a curve, then
the main assertion is proven. Else, $f \times f'$ is surjective and
the pull back $ f^*$ is injective. But then
$ \omega_i \wedge \overline {\eta_j} + \overline{ \eta_i }\wedge 
\omega_j \equiv 0$
contradicts the K\"unneth formula.

Hence, there is only one case left  to consider, namely that, say, all the
$\omega_j$'s are $\C$-linearly dependent. Then we may assume
$\omega_j \equiv 0 ,  \ \forall \ j \geq 2$  and the above equation
yields $ \omega_1 \wedge \overline {\eta_j} = 0, \ \forall j \geq 2$.
But then $ \omega_1 \wedge \eta_j \equiv 0$, since if $\xi$ is the K\"ahler form,
$ |\omega_1 \wedge \eta_j|^2 = \int_X \omega_1 \wedge \overline {\eta_j}
\wedge \overline {\omega_1 \wedge \overline {\eta_j} } \wedge \xi^{n-2} = 0.$

2) follows easily from 1) as follows.

The correspondence is
given by $f \mapsto  V : =f^*(H^1(B, \mathbb{C}))$.

In fact,  since $f : X \to B $ is a
continuous map which induces a surjection of fundamental groups, then the
algebra homomorphism $f^*$ is injective when restricted to
$H^1(B,\C)$ (this statement follows
also  without the K\"ahler
hypothesis) and $ f^* (H^1(B,\mathbb{C})) \subset
H^1(X,\C)$   contains many isotropic subspaces $U$ of dimension $b$
with $ U \oplus \bar{U} = f^*(H^1(B, \mathbb{C}))$.
If such subspace $U$ is not maximal isotropic, then
it is contained in $U'$, which determines an irrational
pencil $f'$ to a curve $B'$ of genus $ > b$, and $f$ factors
through $f'$ in view of the fact that every curve
of positive genus is embedded in
its Jacobian. But this contradicts the fact that $f$ has connected fibres.

\qed

   To give an idea of the power of the above result, let us show how
the following result
due to Gromov ( \cite{gromov}, see also \cite{catgr} for details)
follows  as  a simple
consequence

\begin{cor}
Let $X$ be a compact K\"ahler manifold and assume we have a surjective
morphism $\pi_1 (X) \rightarrow \Gamma$, where $\Gamma$ has a
presentation with $n$ generators, $m$ relations, and with $n-m \geq 2$. Then
there is an irrational pencil $f : X \rightarrow
    B$, such that $2 g(B) \geq n-m$ and $H^1(\Gamma, \mathbb{C})
\subset f^*(H^1(B, \mathbb{C})$.
\end{cor}
\Proof
By the argument we gave in 2) above,
$H^1(\Gamma, \mathbb{C})$ injects into $ H^1(X, \mathbb{C})$
and we claim that each vector $v$ in  $H^1(\Gamma, \mathbb{C})$ is contained
in a nontrivial isotropic
subspace. This follows because  the classifying
space $Y : =  K (\Ga, 1)$ is obtained by attaching
$n$ 1-cells, $m$ 2-cells, and then only
cells of higher dimension. Hence  $ h^2 (\Ga, \Q ) = h^2 (Y, \Q) \leq m$,
and $ w \to w \wedge v$ has a kernel of dimension $ \geq 2$ on
$H^1(\Gamma, \mathbb{C})$. The surjection $\pi_1 (X) \rightarrow \Gamma$
induces a continuous map $ F : X \ra Y$, and each vector in the pull back
of $H^1(\Gamma,
\mathbb{C})$ is contained in a nontrivial maximal isotropic subspace,
thus, by 2) above , in a subspace $  V : =f^*(H^1(B, \mathbb{C}))$
for a suitable irrational pencil $f$.
   Now, the
corresponding subspaces
$V$ are defined over $\mathbb{Q}$ and  $H^1(\Gamma,
\mathbb{C})$ is contained in their union. Hence, by Baire's theorem,
  $H^1(\Gamma,
\mathbb{C})$ is contained in one of them.
\qed

In particular, Gromov's theorem applies to a surjection
$ \pi_1 (X) \to \Pi_g$, where $ g \geq 2 $, and
$ \Pi_g$ is the fundamental group of a compact complex curve of genus $g$.
But in general the genus $b$ of the target curve $B$  will not be equal
to $g$, and we would like to detect $b$ directly
  from the fundamental group $ \pi_1 (X) $.
For this reason (and for others) we need to recall a concept 
introduced by Deligne
and Mostow (\cite{D-M}, see also \cite{isogenous}) in order to 
extend to higher
dimensions some standard arguments about Fuchsian groups.

\begin{df}\label{orbi}
Let $Y$ be a normal complex space and let $D$ be a closed analytic set.
  Let $D_1, \dots , D_r$ be the divisorial
(codimension 1) irreducible components
of $D$, and  attach to each $D_j$ a positive integer $m_j > 1$.

  Then the {\bf orbifold
fundamental group}
$\pi_1^{orb} (Y \setminus D,( m_1, \dots m_r)) $ is defined as the quotient
of  $\pi_1
(Y \setminus ( D_1 \cup \dots \cup D_r)$
by the subgroup normally generated by the $\{ \gamma_1^{m_1}, \dots,
\gamma_r^{m_r}\}$,
where $ \gamma_i$ is a simple geometric loop around the divisor $D_i$
(this means, $ \gamma_i$ is the conjugate via a simple path $\de$  of 
a local loop
$\ga$ which,  in a local coordinate
chart where $Y$ is smooth and $ D_i = \{ (z ) | z_1 = 0 \}$,
is given by  $\ga (\theta) : = (exp (2 \pi i \theta), 0, \dots 0), \ \forall
\theta \in [0,1].$

We observe in fact that another choice for $\ga_i$ gives a conjugate element,
so the group is well defined.

\end{df}

\begin{ex}
Let $Y = \C$, $ D = \{ 0\}$: then $\pi_1^{orb} (\C \setminus \{ 0 \},m)
\cong \Z / m$ and its subgroups correspond to the subgroups
$ H \subset \Z $ such that $ H \supset m \Z$, i.e., $ H = d \Z$, where
$d$ divides $m$.
\end{ex}

The above example fully illustrates the meaning of the orbifold
fundamental group, once we use once more the well known theorem of
Grauert and Remmert ( \cite{g-r})

\begin{oss}
There is a bijection between
$$ {\rm Monodromies}\  \mu : \pi_1^{orb} (Y \setminus D,( m_1, \dots m_r))
\ra \SSS (M)$$
and normal locally finite coverings $ f : X \to Y$, with general 
fibre $\cong M$,
and such
that for each component
$R_i$ of $ f^{-1} ( D_i)$ the ramification index divides $m_i$.
\end{oss}

We have moreover (see \cite{isogenous}) the following

\begin{prop}
Let $X$ be a complex manifold, and $G$ a group of holomorphic automorphisms
of $X$, acting properly discontinuously. Let $D$ be the branch locus of
$ \pi : X \ra Y : = X / G$, and for each divisorial component $D_i$ of $D$
let $m_i$ be the branching index. Then we have an exact sequence
$$ 1 \to \pi_1 (X) \to \pi_1^{orb} (Y \setminus D,( m_1, \dots m_r)) 
\to G \to 1. $$
\end{prop}

\begin{oss}
I) In order to extend the above result to the case where $X$ is only
normal (then $ Y : = X / G$ is again normal), it suffices
to define the orbifold fundamental group of a normal variety $X$ as
$$ \pi_1^{orb} (X) : = \pi_1 (X \setminus Sing (X)).$$

II) Taking the monodromy action of $ \pi_1^{orb} (Y \setminus D,( 
m_1, \dots m_r))$
acting on itself by translations, we see that there exists
a universal orbifold covering space $ \overline {(Y  \setminus D,( m_1, \dots
m_r))}$ with a properly discontinuous action of $ \pi_1^{orb} (Y 
\setminus D,( m_1,
\dots m_r))$ having $Y$ as quotient, and the prescribed ramification.

III) Obviously  the universal orbifold covering space
$ \overline {(Y  \setminus D,( m_1, \dots
m_r))}$ is (connected and) simply connected.
\end{oss}

\begin{ex}
a) Let $Y $ be a compact complex curve of genus $g$, $ D = \{ p_1, 
\dots p_r \}$: then
$ \Ga : = \pi_1^{orb} (Y \setminus \{ p_1, \dots p_r \},(m_1, \dots 
m_r )) $ has
a presentation
$$\Ga : = < \ga_1, \dots, \ga_r, \alpha_1, \beta_1, \dots \alpha_g, \beta_g |
\ga_1 \cdot \dots \cdot \ga_r \cdot \prod_{i=1}^g [\alpha_i, \beta_i 
] =1, \ga_j^{m_j}
= 1>$$

b) $\Ga$ acts on a simply connected complex curve $\Sigma$, with $ 
\Sigma / \Ga \cong
Y$. By the uniformization theorem $\Sigma \cong  \PP^1$ iff  $\Sigma$ 
is compact,
i.e., iff $\Ga$ is finite  (then $ Y\cong  \PP^1$). If instead $\Ga$ 
is infinite,
then there is a finite index subgroup $\Ga'$ acting freely on $\Sigma$.
Then correspondingly we obtain $ C' : = \Sigma / \Ga' \to Y$
a finite covering with prescribed ramification $m_i$ at each point $p_i$.
\end{ex}

\begin{ex}
{\bf Triangle groups}
We let $Y = \PP^1$, $r=3$, without loss of generality $ D = \{ \infty 
, 0 , 1 \}$.
Then the orbifold fundamental group in this case reduces to the previously
defined triangle group
$  T (m_1, m_2,  m_3 ) $ which  has
a presentation
$$T (m_1, m_2,  m_3 ) : = < \ga_1, \ga_2, \ga_3|
\ga_1 \cdot \ga_2  \cdot \ga_3  =1,
\ga_1^{m_1} = 1, \ga_2^{m_2} = 1, \ga_3^{m_3} = 1>.$$

The triangle group is said to be of {\bf elliptic type}
iff $ \Sigma \cong \PP^1$,  of {\bf parabolic type}
iff $ \Sigma \cong \C$,  of {\bf hyperbolic type}
iff $ \Sigma \cong \HH : = \{ \tau | Im (\tau ) > 0\}.$

It is classical (and we have already seen the first alternative as a 
consequence of
Hurwitz' formula in lecture 2) that the three alternatives occur
\begin{itemize}
\item
elliptic $\Leftrightarrow$ $ \sum_i \frac{1}{m_i} > 1$ 
$\Leftrightarrow$ (2,2,n) or
(2,3,n) ($n= 3,4,5$)
\item
parabolic $\Leftrightarrow$ $ \sum_i \frac{1}{m_i} = 1$ 
$\Leftrightarrow$ (3,3,3) or
(2,3,6) or (2,4,4)
\item
hyperbolic $\Leftrightarrow$ $ \sum_i \frac{1}{m_i} < 1$

\end{itemize}
We restrict here to the condition $ 1 < m_i < \infty$, else for instance there
is also the parabolic case ($2,2, \infty$), where the uniformizing function
is $ cos \colon \C \ra \PP^1_{\C}$.

The group $  T (m_1, m_2,  m_3 ) $, which was described for the elliptic
case in lecture 2, is in the parabolic case a semidirect product
of the period lattice $\Lambda$ of an elliptic curve by its group
$\mu_n$ of linear automorphisms
\begin{itemize}
\item
(3,3,3) : $\La = \Z \oplus \zeta_3 \Z $, $\zeta_3$ a generator of $\mu_3$
\item
(2,3,6)  : $\La = \Z \oplus \zeta_3 \Z $, $- \zeta_3$ a generator of $\mu_6$
\item
(2,4,4) : $\La = \Z \oplus i \Z $, $i$ a generator of $\mu_4$.
\end{itemize}

There is a good reason to call the above 'triangle groups'. Look in fact at
the ramified covering $f : \Sigma \to \PP^1$, branched in $\{ \infty 
, 0 , 1 \}$.
Complex conjugation on $\PP^1$ lifts to the covering, as we shall see 
later in more
detail. Consider then a connected component $\De$ of $ f^{-1} (\HH)$. 
We claim that
it is a triangle (in the corresponding geometry: elliptic, resp. Euclidean,
respective hyperbolic) with angles $\pi/ m_1 , \pi/ m_2, \pi/ m_3$.

In fact, take a lift of complex conjugation which is the identity on one
of the three sides of $\De$: then it follows that this side is contained in the
fixed locus of an antiholomorphic automorphism of $\Sigma$, and the 
assertion follows
then easily.

In terms of this triangle (which is unique up to automorphisms of $\Sigma$
in the elliptic and hyperbolic case) it turns out that the three generators
of  $  T (m_1, m_2,  m_3 ) $
are just rotations around the vertices of the triangle, while the
triangle group $  T (m_1, m_2,  m_3 ) $ 
sits as a subgroup of index 2 inside the group
generated by the reflections on the sides of the triangle.

\end{ex}

Let's leave for the moment aside the above concepts, which will be of
the utmost importance
in the forthcoming sections, and let us return to the irrational pencils.

\begin{df}
Let $X$ be a compact K\"ahler manifold and assume we have
a  pencil $f : X \rightarrow B$.  Assume that $t_1, \dots t_r$
are the points of $B$ whose fibres $ F_i : = f^{-1} (t_i) $
are the multiple fibres of $f$. Denote by $m_i$ the multiplicity
of $ F_i$, i.e., the G.C.D. of the multiplicities of the
irreducible components of $F_i$.  Then the {\bf orbifold
fundamental group of the fibration}
$\pi_1 (f) : = \pi_1 (b, m_1, \dots m_r) $ is defined as the quotient
of  $\pi_1
(B \setminus \{ t_1, \dots t_r \})$
by the subgroup normally generated by the $ \gamma_i^{m_i}$'s,
where $ \gamma_i$ is a geometric loop around $t_i$.

The orbifold fundamental group is said to be of {\bf hyperbolic type}
if the corresponding universal orbifold (ramified) covering of $B$ is
the upper half plane.
\end{df}

The orbifold fundamental group of a fibration is a natural object in 
view of the
following result (see \cite{cko}, \cite{fibred})

\begin{prop}
\label{orbex}
Given a fibration   $f : X \rightarrow B$ of
a compact K\"ahler manifold onto a compact complex curve $B$, we have 
the orbifold
fundamental group exact sequence
$ \pi_1 (F) \rightarrow \pi_1 (X) \rightarrow  \pi_1 (b, m_1, \dots m_r)
\rightarrow 0,$  where $F$ is a smooth fibre of $f$.
\end{prop}

The previous exact sequence leads to following result, which is a small
generalization of Theorem 4.3. of \cite{fibred} and a
variant of several other results  concerning fibrations onto curves (see
\cite{isogenous} and
\cite{fibred}), valid more generally for quasi-projective varieties
(in this case the starting point is the closedness of logarithmic 
forms, proven by
Deligne in \cite{del}, which is used in order to obtain  extensions of the theorem 
of Castelnuovo
and De  Franchis to the non complete case,  see \cite{bauer} and 
\cite{arapura}).

\begin{teo}\label{orbfibre}
Let $X$ be a compact K\"ahler manifold and let $(b, m_1, \dots m_r)$
be a hyperbolic type. Then
there is a bijection between pencils $f : X \rightarrow B$ of
type $(b, m_1, \dots m_r)$ and epimorphisms
$\pi_1 (X) \rightarrow  \pi_1 (b,
m_1, \dots m_r)$ with finitely generated kernel.
\end{teo}
\Proof
One direction follows right away from proposition \ref{orbex}, so 
assume that we are
given such an epimorphism. Since $ \pi_1 (b, m_1, \dots m_r)$ is of 
hyperbolic type,
it contains a normal subgroup $H$ of finite index which is isomorphic
to a fundamental group $\Pi_g$ of  a compact curve of genus $ g \geq 2$.

Let $H'$ be the pull back of $H$ in $\pi_1(X)$ under the given surjection,
and let $X' \to X$ the corresponding Galois cover, with Galois group
$G \cong \pi_1 (b, m_1, \dots m_r) / H.$

By the isotropic subspace theorem, there is an irrational pencil
$ f' : X' \ra C$, where the genus of $C$ is at least $g$, corresponding to the
surjection $ \psi: \pi_1 (X') = H' \ra H  \cong \Pi_g$. The group $G$ 
acts on $X'$
leaving the associated cohomology subspace (${f'}^{*} (H^1 (C, \C) $) 
invariant,
whence $G$ acts on $C$ preserving the fibration, and we get a fibration
$ f : X \ra B := C/G$.

By theorem 4.3 of \cite{fibred}, since the kernel of
$\psi$ is finitely generated, it follows that $\psi = f'_*  : \pi_1 (X') \ra
  \Pi_g = \pi_1 (C) $. $G$ operates freely on
$X'$ and effectively on $C$: indeed $G$ acts nontrivially on
$\Pi_g $  by conjugation, since a hyperbolic group has trivial centre.
Thus we get an action of $ \pi_1 (b, m_1, \dots m_r)$ on the upper half
plane $\HH$ whose quotient equals $C/ G : = B$, which has genus $b$.

We use now  again a result from   theorem 4.3 of \cite{fibred}, namely, that
$f'$ has no multiple fibres. Since the projection $ C \to B$ is branched
in $r$ points with ramification indices equal to ($m_1, \dots m_r$),
it follows immediately that the orbifold fundamental group
of $f$ is isomorphic to  $ \pi_1 (b, m_1, \dots m_r)$.
\qed

\begin{oss}
The crucial property of Fuchsian groups which is used in \cite{fibred}
is the so called NINF property, i.e., that every normal nontrivial subgroup
of infinite index is not finitely generated. From this property
follows that, given a fibration $f : X' \ra C$, the kernel
of $f_* : \pi_1 (X') \ra \pi_1 (C)$ is finitely generated (in the 
hyperbolic case) if
and only if there are no multiple fibres.
\end{oss}

\noindent
\subsection{ Varieties isogenous to a product}

\begin{df}
A complex algebraic variety $X$ of dimension $n$ is said to be {\bf
    isogenous to a higher product} if and only if there is a finite \'etale
    cover $C_1 \times \ldots C_n \rightarrow X$, where $C_1, \ldots ,
    C_n$ are compact Riemann surfaces of respective genera $g_i :=
    g(C_i) \geq 2$.
\end{df}
In fact, $X$ is isogenous to a higher product if and only if there is
a finite \'etale Galois cover of $X$ isomorphic to a product of curves
of genera at least two, ie., $X \cong (C_1 \times \ldots C_n)/G$, where
$G$ is a finite group acting freely on $C_1 \times \ldots C_n$.

Moreover, one can prove that there exists a unique  minimal such
Galois realization $X \cong (C_1 \times \ldots C_n)/G$ (see \cite{isogenous}).

In proving this plays a key role a slightly more general fact:

\begin{oss}
The universal covering of a product of curves $C_1 \times \ldots C_n$ 
of hyperbolic type
as above is the polydisk $\HH^n$.

The group of automorphisms of $\HH^n$
is a semidirect product of the normal subgroup $ Aut (\HH)^n$ by the
symmetric group $\SSS_n$ (cf. \cite{ves} VIII, 1 pages 236-238). This result
is a consequence of three basic facts:

i) using the subgroup  $ Aut (\HH)^n$ we may reduce to
consider only automorphisms which leave the origin invariant,

ii) we use the Hurwitz trick to show that the tangent representation
is faithful: if $g(z) =  z + F_m (z) +  \dots $ is the Taylor
development at the origin and with m-th order term $F_m(z) \neq 0$,
then for the r-th iterate of $g$ we get $ z \to  z +r F_m (z) +  \dots $,
contradicting the Cauchy inequality for the r-th iterate when $ r >> 0$,

iii) using the circular invariance of the domain ( $ z \to \la z , |\la|= 1$),
one sees that the automorphisms which leave the origin invariant are linear,
since, if $ g(0) = 0$, then $ g (z)$ and $ \la^{-1} g (\la z )$ have the same
derivative  at the origin, whence by ii) they are equal.

A fortiori, the group of automorphisms of such a product,
$Aut (C_1 \times \ldots C_n)$ has as normal subgroup
$Aut (C_1 )\times \ldots Aut (C_n)$, and with quotient group a subgroup
of $\SSS_n$.

\end{oss}

The above remark leads to the following

\begin{df}
A variety isogenous to a product is said to be {\bf unmixed} if in its minimal
realization $ G  \subset Aut (C_1) \times \ldots Aut (C_n)$.
If $n=2$, the  condition of minimality is equivalent to requiring that
$ G \to Aut (C_i)$ is injective for $ i=1,2$.
\end{df}

The characterization of varieties $X$ isogenous to a (higher) product
becomes simpler in the surface case.
Hence, assume in the following $X = S$ to be a
surface: 
 then
\begin{teo}\label{fabiso}(see \cite{isogenous}). a) A
    projective smooth surface is isogenous to a higher product if and only if
the following two conditions are satisfied:

1) there is an exact sequence
$$
1 \rightarrow \Pi_{g_1} \times \Pi_{g_2} \rightarrow \pi = \pi_1(S)
\rightarrow G \rightarrow 1,
$$
where $G$ is a finite group and where $\Pi_{g_i}$ denotes the fundamental
group of a compact curve of genus $g_i \geq 2$;

2) $e(S) (= c_2(S)) = \frac{4}{|G|} (g_1-1)(g_2-1)$.

\noindent
b) Any surface $X$ with the
same topological Euler number and the same fundamental group as $S$
is diffeomorphic to $S$. The corresponding subset of the moduli space,
$\mathfrak{M}^{top}_S = \mathfrak{M}^{diff}_S$, corresponding to
surfaces orientedly homeomorphic,
resp. orientedly diffeomorphic to $S$, is either
irreducible and connected or it contains
two connected components which are exchanged by complex
conjugation.

In particular, if $X$ is orientedly diffeomorphic to $S$, then $X$ is
deformation equivalent to $S$ or to $\bar{S}$.
\end{teo}

{\em Sketch of the  Proof.}

The necessity of conditions 1) and 2) of a) is clear, since there is
an \'etale  Galois cover of $S$ which is a product, and then $ e(S) \cdot |G| =
e (C_1 \times C_2) = e (C_1 ) \cdot e ( C_2) = 4 (g_1-1)(g_2-1).$

Conversely, take the  \'etale  Galois cover $S'$ of $S$ with group 
$G$ corresponding
to the exact sequence 1). We need to show that $S'$ is isomorphic to a product.

By theorem \ref{orbfibre} the two projections of the direct product
$ \Pi_{g_1} \times \Pi_{g_2}$ yield two holomorphic maps to curves of 
respective
genera $g_1, g_2$, hence
we get a holomorphic map $ F : S' \to C_1 \times C_2$, such that $f_j 
: = p_j \circ F :
S' \to C_j$ is a fibration. Let $h_2$ be the genus of the fibres of $f_1$: then
since  $ \Pi_{g_2}$ is a quotient of the fundamental group of the 
fibre, it follows
right away that $ h_2 \geq g_2$.

We  use then the  classical (cf. \cite{bpv}, proposition 11.4 , page 97)

{\bf Theorem of Zeuthen-Segre} {\em Let $f: S \to B$ be a fibration
of an algebraic surface onto a curve of genus $b$, with fibres of 
genus $g$: then

$$ e (S) \geq 4 (g-1) (b-1), $$ equality holding iff all the fibres are smooth,
or , if $g=1$, all the fibres are multiple of smooth curves.}

Hence $e (S) \geq  4 (g_1-1)(h_2-1) \geq (g_1-1)(g_2-1) = e(S),$ 
equality holds,
$ h_2 = g_2$, all the fibres are smooth and $F $ is then an isomorphism.

Part b): we consider first the unmixed case. This means that the group $G$
does not mix the two  factors, whence the individual subgroups $ \Pi_{g_i}$
are normal in $\pi_1 (S)$, and moding out by the second of them one gets
the exact sequence
$$
1 \rightarrow \Pi_{g_1}   \rightarrow  \pi_1(S) / \Pi_{g_2}
\rightarrow G \rightarrow 1,
$$
which is easily seen to be the orbifold exact sequence for the quotient
map $ C_1 \to C_1 / G$. This immediately shows that the 
differentiable structure
of the action of $G$ on the product $C_1 \times C_2$ is determined, hence also
the differentiable structure of the quotient $S$ is determined by the 
exact sequence 1) in \ref{fabiso}.

We have now to choose  complex structures on the respective manifolds 
$C_i$, which
make the action of $G$ holomorphic. Note that the choice of a complex structure
implies the choice of an orientation, and that once we have fixed the 
isomorphism of
the fundamental group of $C_i$ with $ \Pi_{g_i}$ and we have chosen 
an orientation
(one of the two generators of $H^2 ( \Pi_{g_i} , \Z )$)  we have a 
marked Riemann
surface. Then the theory of Teichm\"uller spaces shows that the space 
of complex
structures on a marked Riemann surface of genus $g \geq 2$ is a 
complex manifold
$\sT_g$ of dimension $ 3 (g-1)$ diffeomorphic to a ball. The finite group
$G$, whose differentiable action is specified, acts on $\sT_g$, and the
fixed point set equals the set of complex structures for which the action is
holomorphic. The result follows then from Proposition 4.13 of \cite{isogenous},
which is a slight generalization of one of the solutions
(\cite{Tro}) of the Nielsen realization problem.

\label{Nielsen} 
\begin{prop} {\bf (Connectivity of Nielsen realization) } 
Given a differentiable action of
a finite group
$G$ on a fixed oriented and marked Riemann surface of genus $g$,  the 
fixed locus
Fix$($G) of
$G$ on
$\sT_g$ is non empty,  connected and indeed diffeomorphic to an euclidean
space.
\end{prop}

Let us first explain why the above proposition
implies part b) of the theorem (in the unmixed
case). Because the moduli space of such
surfaces is then the image of a surjective
holomorphic map from the union of 2 connected
complex manifolds. We get 2 such manifolds because of the choice of
orientations on both factors which together
must give the fixed orientation on our
algebraic surface. Now, if we change the choice
of orientations, the only admissible choice is
the one of reversing orientations on both
factors, which is exactly the result of complex
conjugation.

{\em Idea of proof}
Let us now comment on the underlying idea for the above proposition: 
as already said,
  Teichm\"uller space
$\sT_g$ is diffeomorphic to an Euclidean  space of dimension $6g -6$, 
and  admits a
Riemannian metric,  the Weil-Petersson metric, concerning which
Wolpert and Tromba proved the existence of
  a $C^2$-function $f$ on $\sT_g$ which is
proper, $G$-invariant, non negative ($f\ge  0$), and finally 
such that $f$ is 
strictly convex
for the given metric (i.e., strictly convex along the W-P geodesics).\\

Recall that, $G$ being a finite group, its action can be
linearized at the fixed points, in particular
Fix$(G)$ is a smooth  submanifold.

The idea is to use Morse theory for the function $f$ which is 
strictly convex, and
proper, thus  it  always has a minimum when restricted to a 
submanifold of  $\sT_g$

\begin{itemize}
\item
1)  There is a unique critical point $x_o$  for $f$ on $\sT_g$, which is
an absolute minimum on $\sT_g$ (thus   $\sT_g$ is diffeomorphic to an euclidean
space).
\item
2) If we are given a connected component
$M$ of  Fix$(G)$, then a critical point  $y_o$ for the restriction of 
$f$ to $M$ is also
a critical point for $f$ on  $\sT_g$: in fact $f$ is $G$ invariant, 
thus $df$ vanishes
on the normal space to $M$ at $y_o$.
\item
3) Thus every connected component
$M$ of  Fix$(G)$ contains $x_o$, and, Fix$(G)$ being smooth, it is connected.
  Fix$(G)$ is nonempty since $x_0$, being the
unique minimum, belongs to Fix$(G)$.
\item
4) Since  $f$ is strictly convex, and
proper on   Fix$(G)$, then by Morse theory  Fix$(G)$ is diffeomorphic to an
euclidean space.

\end{itemize}

\qed

In the mixed case  there is a subgroup $ G^o$ of index 2 
consisting of 
transformations which do
not mix the two factors, and a corresponding subgroup $\pi^o$ of $\pi 
= \pi_1 (S)$
of index 2, corresponding to an \'etale double cover $S'$ yielding a
surface of unmixed type. By the first part of the proof, it will 
suffice to show that,
once we have found a lifting isomorphism   of $\pi^o$ with a subgroup
$\Ga^o$ of $Aut (\HH) \times Aut (\HH)$, then the lifting isomorphism 
of $\pi$ with
a subgroup $\Ga$ of $Aut (\HH \times \HH)$ is uniquely determined.

The transformations of $\Ga^o$ are of the form
$ (x,y) \to (\ga_1 (x), \ga_2 (y))$. Pick any transformation in $ \Ga 
\setminus \Ga^o$:
it will be   a transformation of the form
  $  (a (y), b(x))$. Since it normalizes $ \Ga^o$, for each $\de \in 
\Ga^o$ there is
$\ga \in  \Ga^o$  such that
$$ a \ga_2 = \de_1 a ,  \  \ b \ga_1 = \de_2 b .$$
We claim that $a,b$ are uniquely determined. For instance, if $a'$ would also
satisfy $ a' \ga_2 = \de_1 a' $, we would obtain
$$  a ' a^{-1} = \de_1 ( a' a^{-1} ) \de_1 ^{-1}.$$
This would hold in particular for every $\de_1 \in \Pi_{g_1}  $, but
since only the identity centralizes such a Fuchsian group, we 
conclude that $ a' = a$.

\QED

\begin{oss}
A completely similar result holds in higher
     dimension, but the Zeuthen-Segre theorem allows
    an easier  formulation  in dimension two.

One can moreover weaken the hypothesis on the fundamental group, see Theorem B
of \cite{isogenous}.
\end{oss}

\subsection{Complex conjugation and  Real structures}

The interest of Theorem \ref{fabiso} lies in its constructive aspect.

Theorem \ref{fabiso} shows that in order to construct a whole 
connected component of the
moduli space of surfaces of general type, given by surfaces isogenous
to a product, it suffices, in the
   unmixed type, to provide the following data:

\begin{itemize}
\item
i) a finite group $G$
\item
  ii) two orbifold fundamental groups $ A_1: = \pi_1 (b_1, m_1, \dots m_r) $,
$A_2:= \pi_1 (b_2, n_1, \dots n_h) $

\item
iii) respective surjections $\rho_1 : A_1 \to G $, $\rho_2 : A_2 \to G $
such that
\item
iv) if we denote by $\Sigma_i$ the image under $\rho_i$ of the conjugates
of the powers of the generators of $A_i$ of finite order,
then
$$ \Sigma_1 \cap \Sigma_2 = \{ 1_G \} $$
\item
v) each surjection $\rho_i$ is order preserving, in the sense for instance that
  a generator of  $ A_1: = \pi_1 (b_1, m_1, \dots m_r) $ of finite order
$m_i$ has as image an element of the same order $m_i$.

\end{itemize}

In fact, if we take a curve $C'_1$ of genus $b_1$, and $r$ points on it,
to $\rho_1$ corresponds a Galois covering $C_i \to C'_i$ with
group $G$, and the elements of $G$ which have a fixed point on $C_i$
are exactly the elements of $\Sigma_i$. Therefore we have a diagonal action
of $G$ on $C_1 \times C_2$ (i.e., such that $ g (x,y) = 
(\rho_1(g)(x),\rho_2(g)(y))  $,
and condition iv) is exactly the condition that $G$ acts freely on 
$C_1 \times C_2$.

There is some arbitrariness in the above choice, namely, in the choice of
the isomorphism of the respective orbifold fundamental groups with 
$A_1, A_2$, and
moreover one can compose each $\rho_i$ simultaneously with the same
automorphism of $G$ (i.e., changing $G$ up to isomorphism).
Condition v) is technical, but important in order to calculate the genus
of the respective curves $C_i$.

In order to pass to the complex conjugate surface (this is an 
important issue in
Theorem \ref{fabiso}), it is clear that we take the conjugate curve
of each $C'_i$, and the conjugate points of the branch points, but we 
have to be
more careful in looking at what happens with the homomorphisms $\rho_i$.

For this reason, it is worthwhile to recall some basic facts
about complex conjugate structures and real structures.

\begin{df}
Let $X$ be an almost complex manifold, i.e., the pair of a 
differentiable manifold $M$
and an  almost complex structure $J $: then the complex conjugate 
almost complex
manifold $ \bar{X}$ is given by the pair $(M, -J)$. Assume now that 
$X$ is a {\bf
complex manifold}, i.e., that the almost complex structure is 
integrable. Then the same
occurs for $-J$, because, if $ \chi : U \to \C^n$ is a local chart for
$X$, then $\overline {\chi} : U \to \C^n$ is a local chart for $ \bar{X}$.

In the case where $X$ is a projective variety $ X \subset \PP^N$, 
then we easily see
that $ \bar{X}$ equals $\sigma (X)$, where $\sigma : \PP^N \to \PP^N$ is given
by complex conjugation, and the homogeneous ideal of $\bar{X} = 
\sigma (X)$ is the
complex conjugate of the homogeneous ideal $I_X$ of $X$, namely:
$$ I_{\bar{X} }  = \{ P \in \C [z_0 , \dots z_N ] | \overline {P 
(\bar z) } \in I_X  \}.
$$
\end{df}

\begin{df}
Given complex manifolds $X,Y$ let $\phi \colon X \to \bar{Y}$ be a 
holomorphic map.
Then the same map of differentiable manifolds  defines an {\bf antiholomorphic}
map $\bar {\phi} \colon X \to Y $  (also, equivalently,  an {\bf 
antiholomorphic}
map $\bar {\phi}^{**} \colon \bar {X} \to \bar{Y}$).

A map $\phi \colon X \to {Y}$ is said to be {\bf dianalytic} if it is either
holomorphic or antiholomorphic. $\phi$ determines also a dianalytic map
$\phi^{**} \colon \bar {X} \to \bar{Y}$ which is holomorphic iff
$\phi$ is holomorphic.
\end{df}

The reason to distinguish between the maps  $\phi $, $\bar {\phi} 
$  and $\bar {\phi}^{**}$ in the above
definition lies in the fact that maps between manifolds are expressed 
locally as maps in
local coordinates, so in these terms $\bar {\phi} (x) $  is indeed the
antiholomorphic function $\overline {\phi (x)} $, while  $\bar 
{\phi}^{**} (x) =
\phi (\bar{x})$.

With this setup notation, we can further proceed to define the concept of
a real structure on  a complex manifold.

\begin{df}
Let $X$ be a complex manifold.

1) The Klein Group of $X$, denoted by
$ \K l (X) $ or by $Dian (X)$, is the group of dianalytic
automorphisms of $X$.

2) A real structure on $X$ is an antiholomorphic automorphism
$\sigma : X \to X$ such that $\sigma ^2 = Id_ X$.
\end{df}

\begin{oss}
We have a sequence
$$ 0 \ra Bihol (X) : = Aut (X)  \ra Dian (X) : = \K l (X) \ra \Z / 2  \ra 0 $$
which is exact if and only if $X$ is biholomorphic to $\bar{X}$, and splits
if and only if $X$ admits a real structure.
\end{oss}

\begin{ex}
Consider the  anharmonic elliptic curve corresponding to the Gaussian integers:
$ X : = \C / ( \Z \oplus i \Z) $.

Obviously $X$ is real, since $ z \to \bar{z}$
is an antiholomorphic involution.

But there are infinitely many other real structures, since if we take an
antiholomorphism
$\sigma$ we can write   $\sigma (z) = i^r \bar{z} + \mu $,
$\mu = a + i b$, with  $ a,b \in \R / \Z$ and the
condition
$\sigma (\sigma (z)) \equiv
z (mod \ \Z \oplus i \Z)$ is equivalent to
$$ i^r  \bar \mu + \mu = n + im , \ n,m \in \Z \Leftrightarrow
  a + i b +  i^r a - i^{r+1} b = n + im$$

and has the following solutions :
\begin{itemize}
\item

$r= 0, a\in \{0, 1/2\},  b \ {\rm arbitrary} $
\item

$r= 1, a=  - b \ {\rm arbitrary} $
\item

$r= 2,  a \ {\rm arbitrary} , b\in \{0, 1/2\} $
\item

$r= 3, a=   b \ {\rm arbitrary} $.

\end{itemize}
\end{ex}

In the above example the group of biholomorphisms is infinite, and we have
an infinite number of real structures, but many of these are isomorphic, as
the number of isomorphism classes of real structures is equal to the number of
conjugacy classes (for $Aut(X)$) of such splitting involutions.

For instance, in the genus 0 case, there are only two conjugacy classes of
real structures on $\PP^1_{\C}$:
$$ \sigma ( z) = \bar {z}, \  \sigma ( z) = - \frac{1}{ \overline{z}}. $$

They are obviously distinguished by the fact that in the first case the set of
real points $ X(\R) = Fix (\sigma)$ equals $\PP^1_{\R}$, while in the 
second case we
have an empty set. The sign is important, because the real structure
$\sigma ( z) =  \frac{1}{ \overline{z}} $, which has  $\{ z| |z|=1\}$ as set of
real points, is conjugated to the first. Geometrically, in the first 
case we have
the circle of radius 1, $\{ (x,y,z) \in \PP^2_{\C} |  x^2 + y^2 + z^2 = 1\} $,
in the second the imaginary circle of radius $-1$,
$\{ (x,y,z) \in \PP^2_{\C} |  x^2 + y^2 + z^2 = - 1\} $.

It is clear from the above discussion that there can be curves $C$ which are
isomorphic to their conjugate, yet do not need to be real: this fact 
was discovered
by C. Earle, and shows that the set of real curves is only a semialgebraic set
of the complex moduli space,
because it does not coincide with the set  $\frak M_g (\R)$ of real 
points of $\frak
M_g$.

We want now to give some further easy example of this situation.

We observe preliminarily that $C$ is isomorphic to $\bar{C}$ if and 
only in there is
a finite group $G$ of automorphisms such that $C/G$ has a real structure
which lifts to an antiholomorphism of $C$ (in fact, if  $C \cong 
\bar{C}$ it suffices
to take $ Aut (C) = G$ if $ g(C) \geq 2$).

We shall denote this situation by saying that {\bf the covering $ C 
\to C/G$ is real}.

\begin{df}
We shall say that the covering $ C \to C/G$ is an {\bf n-angle 
covering } if $C/G \cong
\PP^1$ and the branch points set consists of n points.

We shall say that $C$ is an {\bf n-angle curve} if $ C \to C/ Aut (C)$ is an
n-angle covering.
\end{df}

\begin{oss}
a) Triangle coverings  furnish an example of a moduli space $(C,G)$, 
of the type
discussed above, which consists of a single point.

b) If $ C \to C/G$ is an {\bf n-angle covering } with n odd, then the induced
real structure on $C/G \cong \PP^1$ has a non empty set of real 
points (the branch locus
$B$ is indeed invariant), thus we may assume it to be  the standard 
complex conjugation
$ z \mapsto  \bar {z}$.
\end{oss}

\begin{ex}\label{nonreal}
We  construct  here examples of  families of real quadrangle covers 
$C \to C/G$ such
that $(C/G )(\R) = \emptyset$, and such that, for a general curve in the 
family, $ G = Aut
(C)$, and the curve $C$ is not real. The induced real structure on $ 
C/G \cong \PP^1$
is then $\sigma (z) = - \frac{1}{ \overline{z}}$, and the quotient 
$(C/G) / \sigma
\cong \PP^2_{\R}$.

We choose then as branch set $ B \subset \PP^1_{\C}$ the set $ \{\infty , 0, w,
- \frac{1}{ \overline{w}}   \}$, and denote by $0,u$ the 
corresponding image points
in $ \PP^2_{\R}$.

Observe now that
$$ \pi_1 ( \PP^2_{\R} \setminus \{0,u\}) = \langle a,b,x | ab = x^2 
\rangle \cong
\langle a,x  \rangle$$
and the \'etale double covering $ \PP^1_{\C} \to  \PP^2_{\R}$ corresponds to
the  quotient  obtained by setting $ a=b=1$, thus $ \pi_1 ( 
\PP^1_{\C} \setminus B)$
is the free group of rank 3
$$ \pi_1 ( \PP^1_{\C} \setminus B ) = \langle a, x^2, x^{-1} a x 
\rangle \cong \langle
a,b, a' := x^{-1} a x , b'  := x^{-1} bx | ab = a' b' \rangle.$$

We let $G'$ be the group $ (\Z/ 2n) \oplus (\Z/ m)$, and 
let $C$ be the Galois cover of
$\PP^2_{\R} $ branched in $ \{0,u\}$ corresponding to the epimorphisms
such that $ x \mapsto (1,0) , a   \mapsto (0,1) $. It follows that 
$C$ is a 4-angle
covering with group $ G \cong (2  \Z/ 2n \Z )\oplus (\Z/ m \Z)$. It is 
straightforward to
verify the following

{\bf Claim:  $G'$ contains no antiholomorphism of order 2, if $n$ is even.}

Thus it follows that $C$ is not real, provided that 
 $ G = Aut(C)$. To simplify things, let
$ n = 4, m = 2$. By Hurwitz' formula $C$ has genus $3$, and $ 8 = |G| 
= 4 (g-1)$.
Assume that $ Aut (C) \neq G$. If $G$ has index 2, then we get an involution
on $\PP^1$ preserving the branch set $B$. But the cross-ratio of the 
4 points equals
exactly $  - \frac{1}{|w|^2 } $, and this is not anharmonic for $w$ general
(i.e., $\neq 2,-1, 1/2$).
If  instead $ C \to C / Aut (C)$ is a triangle curve, then we get only a 
finite number
of curves, and again a finite set of values of $w$, which we can exclude.

Now, since $ | Aut (C) | > 8 (g-1)$, if $ C \to C / Aut (C)$ is not a 
triangle curve,
then the only possibility, by the Hurwitz'  formula, is that we have a quadrangle
cover with branching indices $ (2,2,2,3)$. But this is absurd, since 
a ramification
point of order 4 for $ C \to C / G$ must have a higher order of 
ramification for
the map $ C \to C / Aut (C)$.

\end{ex}

There is  however one important special case  when a curve isomorphic 
to its conjugate
must be real, we have namely the following

\begin{prop}
Let $ C \to C / G$ be a triangle cover which is real and has 
distinct branching
indices ($ m_1 < m_2 < m_3$) : then $C$ is real
(i.e., $C$ has a real structure).
\end{prop}

\Proof
Let $\sigma$ be the real structure on $C/G \cong \PP^1$. The 3 branch points of
the covering must be left fixed by $\sigma$, since the branching indices are
distinct (observe that $\mu (\sigma_* \ga_i)$ is conjugate to $\mu (\ga_i)$,
whence it has the same order).
Thus, without loss of generality we may assume that the three branch 
points are real,
and indeed equal to $\{ 0, 1, \infty \}$, while $ \sigma (z) = \bar{z}$.

Choose 2 as base point, and a basis of the fundamental group as in figure 4:

\begin{figure}[htbp]

% per fare la freccetta ai laccetti
\font\mate=cmmi8
\def\min{\hbox{\mate \char60}\hskip1pt}
\begin{center}
\begin{picture}(350,100)(0,-40)
%
%riga di base
%
\put(0,0){\line(1,0){300}}
\put(238,0){\line(1,0){80}}
%
% punti sulla retta
%
\put(20,0){\circle*{3}}
\put(17,4){0}
\put(90,0){\circle*{3}}
\put(87,4){1}
\put(160,0){\circle*{3}}
\put(156,3){2}
\put(228,0){\circle*{3}}
\put(221,3){$\infty$}
\put(300,0){\circle*{3}}
\put(296,4){-1}
{\thicklines
%
% Beta
%
\put(20,0){\circle{30}}
\put(36,0){\line(1,0){39}}
\put(90,0){\oval(31,31)[t]}
\put(105,0){\line(1,0){55}}
\put(17,14){$\min$}
\put(15,-26){$\beta$}
%
% Alpha
%
\put(300,0){\circle{30}}
\put(160,0){\line(1,0){124}}
\put(297,14){$\min$}
\put(295,-23){$\alpha$}
}
\end{picture}
\end{center}

\caption{The loops $\alpha$ and $\beta$.}
%\label{figura1}
\end{figure}
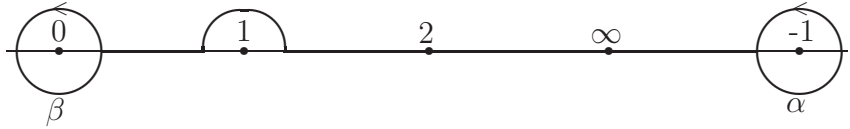

$$ \pi_1 ( \PP^1 \setminus \{ 0, 1, \infty \}, 2) = \langle \alpha,
\beta, \ga | \alpha \beta \ga = 1 \rangle, \ \sigma_* \alpha = \alpha^{-1},
  \sigma_* \ga = \ga^{-1}. $$

Now , $\sigma$ lifts if and only if the monodromy $\mu$ of the $G$-covering is
equivalent to the one of $\mu \circ \sigma_*$ by an inner automorphism
$ Int (\phi )$ of the symmetric group which yields a group automorphism
  $ \psi \colon G \to G$. Set $ a : = \mu (\alpha), b : = \mu (\beta)$. Then
these two elements generate $G$, and since $ \psi (a) = a^{-1}, \psi 
(b) = b^{-1}$
it follows that $\psi$ has order 2, as well as the corresponding covering
transformation. We have shown the existence of the desired real structure.

\QED

We shall now give a simple example of a nonreal triangle cover, based
on the following

\begin{lem}\label{nonconj}
Let $G$ be the symmetric group $\mathfrak{S}_n$ in $n \geq 7$ letters,
let  $ a : = (5,4,1) (2,6)$, $ c: = (1,2,3)(4,5,6, \dots, n)$.

Assume that $n$ is not divisible by $3$: then

1) there is no
automorphism $\psi$ of $G$ carrying $ a \ra a ^{-1}$, $c \ra c^{-1}$.

2) $\mathfrak{S}_n = <a,c>$.

3) The corresponding triangle cover is not real.

\end{lem}

\Proof 1) Since $n \neq 6$, every automorphism of $G$ is an inner one.
If there is a permutation $g$ conjugating $ a $ to $ a ^{-1}$, $c $ to
$ c^{-1}$, $g$ would leave each of the sets $\{1,2,3\}$,
       $\{4,5,\dots ,n\}$,
$\{1,4,5 \}$, $\{2,6 \}$ invariant. By looking
at their
intersections we conclude that $g$ leaves the elements
$1,2,3,6$ fixed. But then $gcg^{-1} \neq
c^{-1}$.

2) Observe that $ a^3$ is a transposition: hence, it suffices to show that the
group generated by $a$ and $c$ is 2-transitive. Transitivity being obvious,
let us consider the stabilizer of $3$.
Since $n$ is not divisible by $3$, the stabilizer of $3$
contains the cycle  $(4,5,6, \dots, n)$; since it contains the transposition
$(2,6)$ as well as $(5,4,1)$, this stabilizer is transitive on 
$\{1,2,4,\dots, n\}$.

3) We have $ ord (a) = 6$,$ ord (c) = 3 (n-3)$, $ord (b) = ord (ca) =
ord ((1,6,3) (4,2,7, \dots n)) = LCM (3, (n-4))$. Thus the orders are distinct
and the nonexistence of such a $\psi$ implies that the triangle cover 
is not real.

\QED

We can now go back to theorem \ref{fabiso}, where the surfaces homeomorphic
to a given surface isogenous to a product were forming one or two connected
components in the moduli space. The case of products of curves is an 
easy example where
we get one irreducible component, which is self conjugate. We show now
  the existence of countably many cases where there are two distinct connected
components.

\begin{teo}
Let $S = (C_1 \times C_2)/ G$ be a surface isogenous to a product of 
unmixed type, with
$ g_1 \neq g_2$. Then $S$ is deformation equivalent to $\bar{S}$ if and only
if $(C_j, G)$ is deformation equivalent to $(\overline{C_j}, G)$ for $ j=1,2$.
In particular, if $(C_1, G)$ is rigid,i.e., $C_1 \to C_1 / G$ is a 
triangle cover,
$S$ is deformation equivalent to
$\bar{S}$  only if $(C_1, G)$ is isomorphic to  $(\overline{C_1}, G)$.
There are infinitely many connected components  
of the moduli space of surfaces of general type for which  $S$ is
not  deformation
equivalent to $\bar{S}$.
\end{teo}

\Proof
  $\bar{S} =\overline{ (C_1 \times C_2)}/ G = (\overline{ (C_1} \times 
\overline{ C_2)}/
G $ and since $g_1 \neq g_2$ the normal subgroups $\Pi_{g_j}$ of the 
fundamental group
$ \pi_1 (\bar{S})$ are uniquely determined. Hence $\bar{S}$ belongs to the same
irreducible connected component containing $S$ (according to the key 
Proposition) if and
only if
$(\overline{C_j}, G)$ belongs to the same
irreducible connected component containing $(C_j, G)$.

We consider now  cases where $C_1 \to C_1 / G$ is a triangle cover, but not
isomorphic to  $(\overline{C_1}, G)$: then clearly $S$ is not deformation
equivalent to $\bar{S}$.

We let , for $ n \geq 7,  n \neq 0 (mod 3)$, $C_1 \to C_1 / G$  be the nonreal
triangle cover provided by Lemma \ref{nonconj}. Let $ g_1 $ be the 
genus of $C_1$,
  observe that $ 2 g_1 - 2 \leq (5/6) \  n!$
and consider an arbitrary integer $g \geq 2$ and a surjection $ \Pi_g 
\to \SSS_n$
(this always exists since  $ \Pi_g$ surjects onto a free group with 
$g$ generators).

The corresponding \'etale covering of a curve $C$ of genus $g$ is a curve $C_2$
with genus $ g_2 > g_1$ since $ 2 g_2 - 2 \geq (2 g -2) \  n!  \geq  2 \  n!$.
The surfaces $ S = C_1 \times C_2$ are our desired examples, the action of
$G = \SSS_n$ on the product is free since the action on the second 
factor is free.

\QED

Kharlamov and Kulikov gave ( \cite{k-k}, \cite{k-k2} ) rigid examples 
of surfaces $S$
which are not isomorphic to their complex conjugate, for instance they
considered a $( \Z/ 5)^2$ covering of the plane branched on
the 9 lines in the plane $\PP^2$ dual to the 9 flexes of a cubic, the 
Fermat cubic for
example. These examples have \'etale coverings which were constructed 
by Hirzebruch
(\cite{hirz}, see also \cite{bhh} ) in order to produce simple examples
of surfaces on the Bogomolov Miyaoka Yau line $ K^2 = 3 c_2$, which, 
by results of
Yau and Miyaoka (\cite{yau}, \cite{miy}) have the unit ball in $\C^2$ 
as universal
covering, whence they are strongly rigid according to a theorem of Mostow
(\cite{mostow}): this means that any surface homotopically equivalent 
to them is either
biholomorphic or antibiholomorphic to them.

Kharlamov and Kulikov prove that the Klein group of such a surface $S$
consists only of the above group $( \Z/ 5)^2$ of biholomorphic 
transformations,
for an appropriate choice of the  $( \Z/ 5)^2$ covering, such that to 
pairs of conjugate
lines correspond pairs of elements of the group which cannot be obtained from
each other by the action of a single automorphism of the group $( \Z/ 5)^2$.

In the next section we shall show how to obtain rigid examples with
surfaces isogenous to a product.

\bigskip

\subsection{Beauville surfaces}

\begin{df}
A surface $S$ isogenous to a higher product is called a {\bf Beauville
    surface} if and only if $S$ is rigid.
\end{df}

This definition is motivated by the fact that Beauville constructed
such a surface in \cite{Bea} , as a quotient $F \times F$ of two
Fermat curves of degree 5 (and genus 6). Rigidity was observed in 
\cite{isogenous}.

\begin{ex}{\bf ('The' Beauville surfaces)}
Let $F$ be the plane Fermat 5-ic $ \{ x^5 + y^5 + z^5 = 0\}$.
The group $( \Z/ 5)^2$ has a projective action obtained
by multiplying the coordinates by 5-th roots of unity. The
set of stabilizers is given by the multiples of $ a: = e_1, b:= 
e_2,c:=  e_1 + e_2$,
where $ e_1 (x,y,z) =   (\epsilon x,y,z)$, $ e_2 (x,y,z) =   (x, 
\epsilon y,z), \epsilon : = exp ( 2 \pi i /5)$.
In other words, $F$ is a triangle cover of $\PP^1$ with group $( \Z/ 5)^2$
and generators $ e_1, e_2, - (e_1 + e_2)$.
The set $\sigma$ of stabilizers is the union of 3 lines in the vector space
  $( \Z/ 5)^2$, corresponding to 3 points in $\PP^1_{\Z/ 5}$.
Hence, there is an automorphism $ \psi$ of $( \Z/ 5)^2$  such that
$\psi (\Sigma) \cap \Sigma = \{ 0 \}$.
Beauville lets then $( \Z/ 5)^2$ act on $F \times F$ by the action
$ g (P,Q) : = ( g P, \psi (g) Q)$, which is free and yields a surface
$S$ with $ K^2 _S = 8$, $p_g = q = 0$. It is easy to see that such a
surfaces is not only real, but defined over $\Q$. It was pointed out
in \cite{fano} that there are exactly two isomorphism classes of
such  Beauville surfaces.
\end{ex}

Let us now construct some Beauville surfaces which are not isomorphic to
their complex conjugate.

  To do so, we observe that the datum of an unmixed  Beauville surface 
amounts to a
purely group theoretical datum, of two systems
of generators $\{a,c\}$ and $\{a', c' \}$ for a finite group $G$ such that,
defining $b$ through the equation $ abc = 1$, and the stabilizer set 
$\Sigma(a,c)$ as
$$\cup _{i \in \N , g \in G} \{ ga^i g^{-1}, gb^i g^{-1}, gc^i g^{-1}\}  $$
the following condition must be satisfied, assuring that the diagonal 
action on the
product of the two  corresponding triangle curves is free
$$\Sigma(a,c) \cap \Sigma(a',c') = \{ 1_G \} . $$

\begin{ex}
Consider the symmetric group $\mathfrak{S}_n$ for $ n \equiv 2 (mod \ 3)$,
define elements $ a, c \in \mathfrak{S}_n $ as in Lemma \ref{nonconj},
and define further
$ a' : =  \sigma^{-1}$, $c' : = \tau
\sigma^2$, where
$\tau := (1,2)$
and $\sigma: = (1,2, \dots ,
n)$. It is obvious that $\mathfrak{S}_n = <a',c'>$.
    Assuming $n \geq 8$ and $n \equiv 2 (3)$,
it is easy to verify that $\Sigma(a,c) \cap \Sigma(a',c') = \{1\}$,
    since one observes that elements which are conjugate in 
$\mathfrak{S}_n$ have
the same type of cycle decomposition. The types in $\Sigma(a,c)$ are
derived from $(6)$, $ (3n-9)$, $(3n-12)$, (as for instance
$(3)$, $(2)$, $(n-4)$ and $(n-3)$) since we assume that
$3$ does neither divide $n$ nor $n-1$, whereas the types in
$\Sigma(a',c')$ are derived from $(n)$, $(n-1)$, or $(\frac{n-1}{2},
\frac{n+1}{2})$.

One sees therefore (since $g_1 \neq g_2$) that the pairs $(a,c), (a', c')$ determine 
Beauville surfaces which are not isomorphic to their complex conjugates.
\end{ex}

Our knowledge of Beauville surfaces is still rather unsatisfactory, 
for instance the
following question is not yet completely answered.
\begin{question}
Which groups $G$ can occur?
\end{question}

It is easy
  to see (cf. \cite{bcg}) that if the group $G$ is abelian,
then it can only be $ (\Z / n)^2$, where $G.C.D. (n,6) = 1$.

Together with I. Bauer and F. Grunewald, we proved in \cite{bcg}
(see also \cite{almeria}) the following results:

\begin{teo}
1) The following groups admit unmixed Beauville structures:

a) $\mathfrak{A}_n$ for large $n$,

b) $\mathfrak{S}_n$ for $n\in \NN$ with $n\ge 7$,

c) ${\bf SL}(2,\F_p)$, ${\bf PSL}(2,\F_p)$ for $p \neq 2,3,5$.
\end{teo}

After checking that all finite simple nonabelian groups of order 
$\leq 50000$,  with the
exception of $\mathfrak{A}_5$,
admit unmixed Beauville structures, we were led to the following
\begin{conj} All finite simple nonabelian groups except $\mathfrak{A}_5$
admit an unmixed Beauville structure.
\end{conj}

Beauville surfaces were extensively studied in
\cite{bcg} (cf. also \cite{almeria}) with special regard to the effect of complex 
conjugation  on them.

\begin{teo}
There are Beauville surfaces $S$ not biholomorphic to $\bar{S}$  with group
\begin{itemize}
\item[1)] the symmetric group $\mathfrak{S}_n$ for any $n\ge 7$,
\item[2)] the alternating group $\mathfrak{A}_n$ for $n\ge 16$ and
$n\equiv 0$ mod $4$,
$n\equiv 1$ mod $3$, $n\not\equiv 3,4$ mod $7$.
\end{itemize}
\end{teo}

We got also examples  of isolated real points in the moduli space
which do not correspond
to real surfaces :

\begin{teo}\label{nonreal}
Let $p>5$ be a prime with $p\equiv 1$ mod $4$, $p\not\equiv 2,4$ mod 5,
$p\not\equiv 5$ mod $13$ and $p\not\equiv 4$ mod $11$. Set $n:=3p+1$.
Then
there is  a Beauville surface $S$ with group $\mathfrak{A}_n$
which is biholomorphic to  its conjugate $\bar{S}$, but is not
real.
\end{teo}

Beauville surfaces of the mixed type also exist, but their 
construction turns out to be
quite more complicated (see \cite{bcg}). Indeed (cf. \cite{fedya})
the group of smallest order has order 512.

\newpage

\section{Lecture 5: Lefschetz pencils, braid and mapping class groups, and
diffeomorphism of ABC-surfaces.}
\label{fifth}

\subsection{Surgeries.}

The most common surgery is the connected sum, which we now describe.

Let $M$ be a manifold of real dimension $m$, thus
for each point $ p \in M$
there is an open set $ U_p$ containing $p$ and a  homeomorphism 
(local coordinate chart)
$ \psi_p : U_p \ra V_p \subset \R^m$ onto an open set $V_p$ of $\R^m$
such that (on its domain of definition)

$\psi_{p'} \circ \psi_p^{-1}$ is a:

\begin{itemize}
\item
Homeomorphism (onto its image) if   $M$ is a {\bf  topological manifold }
\item
Diffeomorphism (onto its image) if     $M$ is a {\bf differentiable manifold}
\item
Biholomorpism (onto its image) if   $M$ is a {\bf  complex manifold }
(in this last case $m = 2n$, $\R^m = \C^n$).

\end{itemize}

\begin{df}
The  operation of connected sum $ M_1 \sharp M_2$ can be done
for two differentiable or topological manifolds of the same dimension.

Choose respective points $p_i \in M_i$ and local charts
$$ \psi_{p_i}: U_{p_i}  \ra \cong B(0, \epsilon_i):
= \{ x \in \R^m | |x| < \epsilon_i \}.$$
Fix positive real numbers $ r_i < R_i < \epsilon_i $
such that $$(**) \  R_2 /  r_2 = R_1 /  r_1 $$ and set $M^*_i : =
M_i \setminus \psi_{p_i}^{-1} (\overline { B(0,
r_i)} )$: then $M^*_1$ and $M^*_2$ are glued together
  through the   diffeomorphism
$\psi : N_1 : = B(0, R_1) \setminus \overline { B(0,
r_1)}  \ra   N_2 : = B(0, R_2) \setminus \overline { B(0,
r_2)}$ such that
$\psi (x_1) =  \frac{R_2 r_1}{|x_1|}  \tau (x_1) $ where either
$\tau (x) = x$, or $\tau (x)$ is an orientation reversing linear isometry
(in the case  where the manifolds $M_i$ are oriented, we might prefer,
in order to furnish the connected sum  $ M_1 \sharp M_2$  of a compatible
orientation,
to have that $\psi$ be orientation preserving).

In other words the connected sum  $ M_1 \sharp M_2$ is the quotient 
space of the
disjoint union $ (M^*_1  ) \cup^o  (M^*_2 )$ through the equivalence relation
which identifies
$ y \in \psi_{p_1}^{-1} (N_1) ) $ to $  w \in \psi_{p_2}^{-1} (N_2) ) $ iff
$$ w = \psi_{p_2}^{-1} \circ \psi \circ \psi_{p_1} (y).$$

\end{df}

We have the following

{\bf Theorem}
{\em The result of the operation of connected sum is independent of 
the choices made.}

An elementary and detailed proof in the differentiable case (the one 
in which we are
more interested) can be found in \cite{b-j}, pages 101-110.

\begin{ex}
The most intuitive example (see Figure 5) is the one of two
  compact orientable Riemann surfaces $ M_1, M_2$ of respective genera $g_1,
g_2$: $ M_1 \sharp M_2$ has then genus $g_1 + g_2$.
In this case, however,
if $ M_1, M_2$ are endowed of a complex structure, we
can even define a connected sum as complex manifolds,
setting $\psi (z_1) =   e^{2 \pi i \theta}  \ \frac{R_2 r_1}{z_1} $.

Here, however, the complex structure is heavily dependent on the parameters
$p_1, p_2$, $ e^{2 \pi i \theta}$, 
and $  R_2  
r_1 = R_1   r_2 $.

In fact, if we set $ t : = R_2 r_1 e^{2 \pi i \theta} \in \C$, we see that
$ z_1 z_2 = t $, and if $t \to 0$ then it is not difficult to see that
the limit of  $ M_1 \sharp M_2$ is the singular curve obtained from $M_1, M_2$
by glueing the points $p_1, p_2$ to obtain the node $ z_1 z_2 = 0 $.

This interpretation shows that we get in this way all the curves near the
boundary of the moduli space $\frak M_g$. It is not clear to us in this moment
how big a subset of the moduli space one gets through  iterated connected sum
operations. One should however point out that many of the conjectures made
about the stable cohomology ring $H^* (\frak M_g, \Z )$ were suggested by
the possibility of interpreting the connected sum as a sort of 
$H$-space structure
on the union of all the moduli spaces  $\frak M_g$ (cf. \cite{mumshaf}).
\end{ex}

\begin{figure}[htbp]
\begin{center}
\scalebox{1}{\includegraphics{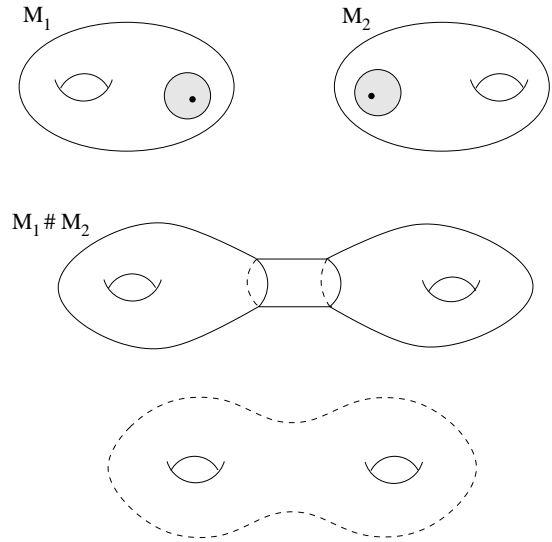}}
\end{center}
\caption{The Connected Sum}
\label{figura1}
\end{figure}

\begin{oss}

1) One cannot perform a connected sum operation for complex manifolds
of dimension $ >1$. The major point is that there is no biholomorphism
bringing the inside boundary of the ring domain $N_1$ to the outside
boundary of $N_2$. The reason for this goes under the name of holomorphic
convexity: if $ n \geq 2$ every holomorphic function on $N_1$  has, by Hartogs'
theorem, a holomorphic continuation to the ball $ B ( 0, R_1)$. While,
for each point $p$ in the outer boundary, there is a holomorphic
function $f$ on $N_1$ such that $ lim_{z \to p} | f(z)| = \infty $.

2) The  operation of connected sum makes the diffeomorphism classes of
manifolds of the same dimension $m$ a semigroup:
  associativity holds, and as neutral element we have the sphere
  $ S^m : = \{ x \in \R^{m+1}| |x| = 1\}$.

3) A manifold $M$ is said to be {\bf irreducible} if $ M \cong M_1 
\sharp  M_2$ implies
that either
$M_1$ or $M_2$ is homotopically equivalent to a sphere  $S^m$.
\end{oss}

A further example is the more general concept of

\begin{df}
{\bf (SURGERY )}

   For $i=1,2$,  let  $ N_i \subset  M_i$ be a  differentiable submanifold.

Then there exists (if $M_i= \R^N$ this is an easy consequence of the 
implicit function
theorem) an open set  $U_i \supset N_i$ which is diffeomorphic to the 
normal bundle
$\nu_{N_i}$ of the embedding $N_i \ra M_i$, and through a 
diffeomorphism which carries
$N_i$ onto the zero section of $\nu_{N_i}$.

Suppose now that we have diffeomorphisms $\phi : N_1 \ra N_2$, and
$\psi : (\nu_{N_1} - N_1)\ra (\nu_{N_2} - N_2)$, the latter compatible
with the  projections
$p_i : \nu_{N_i} \ra N_i$ (i.e.,  $p_2 \circ \psi = \phi \circ p_1$), and with 
the property
of being orientation reversing on the  fibres. We can then define as before a
manifold $ M_1 \sharp_{\psi} M_2$, the quotient of the disjoint  union
$ (M_1 - N_1 ) \cup^o  (M_2 - N_2)$ by the
equivalence relation identifying $(U_1 - N_1 )$ with $(U_2 - N_2 )$ through
the diffeomorphism induced by $\psi$.
\end{df}

\begin{oss}
This time the result of the operation depends upon the choice of
$\phi$ and $\psi$.
\end{oss}

The two surgeries described above combine together in the special situation
of the fibre sum.

\begin{df}
{\bf (FIBRE SUM )}

   For $i=1,2$,  let  $f_i \colon M_i \to  B_i$ be a proper surjective 
differentiable
map  between
differentiable manifolds, let $ p_i \in B_i$ be a noncritical value, 
and let $ N_i
\subset  M_i$ be the corresponding  smooth fibre $ N_i : = f_i ^{-1} (p_i)$.

Then there exists a natural trivialization  (up to a constant matrix) 
of the normal
bundle $\nu_{N_i}$ of the embedding $N_i \ra M_i$, and if we assume as before
  that we have a diffeomorphism $\phi : N_1 \ra N_2$
we can perform a surgery $ M : =  M_1 \sharp_{\phi} M_2$, and the new manifold
$M$ admits
  a proper surjective differentiable map onto the connected sum $ B : =  B_1
\sharp B_2$.
\end{df}

The possibility of variations on the same theme is large: for instance, given
  $f_i \colon M_i \to  B_i \ (i=1,2)$  proper surjective differentiable
maps  between
differentiable manifolds with boundary, assume that $ \partial M_i \to
\partial B_i$ is a fibre bundle, and there are  compatible  diffeomorphisms
$\phi : \partial B_1 \ra \partial B_2$ and $\psi : \partial M_1 \ra 
\partial M_2$:
then we can again define the fibre sum $ M : =  M_1 \sharp_{\psi} M_2$
which admits a  proper surjective differentiable
map onto $ B: = B_1 \sharp_{\phi} B_2$.

In the case where $( B_2, \partial B_2)$ is an euclidean  ball with a 
standard sphere as
boundary, and $ M_2 = F \times B_2$, the question about unicity (up 
to diffeomorphism)
of the surgery procedure is provided by a homotopy class.
Assume in fact that we have two attaching diffeomorphisms
$\psi, \psi ' : \partial M_1 \ra  F \times  \partial B_2$. Then from them we
construct  $ \Psi :=  \psi ' \circ \psi ^{-1}: F \times \partial B_2 
\ra  F \times
\partial B_2$, and we notice that $\Psi (x,t) =  (\Psi_1 (x,t) , \Psi_2 (t))$,
where $ \Psi_2 (t) = \phi ' \circ \phi ^{-1}$. We can then construct 
a classifying
map $ \chi \colon \partial B_2 \cong S^{n-1} \ra D iff (F)$ such that
$$\Psi_1 (x,t) = \chi ( \Psi_2 (t))(x).$$
We get in this way a free homotopy class $[\chi]$, on which the diffeomorphism
class of the surgery depends. If this homotopy class is a priori trivial,
then the result is independent of the choices made: this is the case 
for instance if
$F$ is a compact complex curve of genus $ g \geq 1$.

In order to understand better the unicity of these surgery operations,
and of their compositions, we therefore see the
  necessity of a good  understanding of isotopies of diffeomorphisms. 
To this topic
is devoted the next subsection.

\subsection{Braid and mapping class groups}

E. Artin introduced
the definition of the {\em braid group}
(cf. \cite{art1}, \cite{art}), thus allowing a remarkable extension of
Riemann's concept of monodromy of algebraic functions.  Braids are a
powerful tool, even if not so easy to handle, and especially 
appropriate for the study
of the differential topology of algebraic varieties, in particular of algebraic
surfaces.

\begin{oss}
We observe that the subsets $\{w_1, \ldots, w_n\} \subset \C$
of $n$ distinct points in $\C$ are in
one to one correspondence with monic polynomials $P(z) \in \C [z]$ of
degree $n$ with non vanishing discriminant $\delta (P)$.
\end{oss}

\begin{df}
Let $\C [z]^1_n$ be the affine space of monic polynomials of degree $n$.
Then the group
$$
\mathcal{B}_n := \pi_1 (\C [z]^1_n \setminus \{P |   \delta (P) =
0\}),
$$

i.e., the fundamental group of the space of monic polynomials of degree $n$
having $n$ distinct roots, is called {\em Artin's braid group}.
\end{df}

Usually, one takes as base point the polynomial $P(z) = (\prod_{i=1}^n
(z - i)) \in \C [z]^1_n$ (or the set $\{1, \ldots , n \}$).

To a closed (continuous) path $\alpha : [0,1] \rightarrow (\C [z]^1_n 
\setminus \{P
| \delta (P) = 0\})$ one can associate the subset  
$ B_{\alpha} : = \{ (z,t) \in \C \times \R ~~|
~~  \alpha_t(z) := \alpha (t) (z) = 0\}$ of $\R^3$, which gives a 
visually suggestive
representation of the associated braid.

It is however customary to view a braid as moving from up to down, that is,
to associate to $\alpha $ the set $ B'_{\alpha} : = \{ (z,t)|  (z, -t) \in B_{\alpha}\}.$

Figure \ref{figura3} below shows two realizations of the same braid.

\begin{figure}[htbp]
\begin{center}
\scalebox{1}{\includegraphics{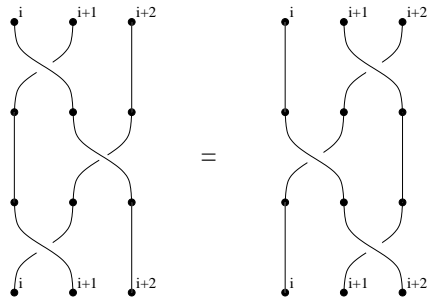}}
\end{center}
\caption{Relation $aba=bab$ on braids}
\label{figura3}
\end{figure}

\begin{oss}
  There is a lifting of $\alpha$ to $\C^n$, the space of ordered
$n$-tuples of roots of monic polynomials of degree $n$, hence there are
(continuous) functions $w_i(t)$
such that $w_i(0) = i$ and $\alpha_t(z) = \prod_{i=1}^n (z -
w_i(t))$.

It follows that to each braid  is   naturally associated a permutation
$\tau \in \mathfrak{S}_n$ given by $\tau(i) := w_i (1)$.
\end{oss}

Even if it is not a priori evident, a very powerful generalization of 
Artin's braid
group was  given by M. Dehn (cf. \cite{dehn}, we refer also to the book \cite
{birman}).

\begin{df}
Let $M$ be a differentiable manifold, then the {\bf mapping class
   group} (or {\em Dehn group}) of $M$  is the group
$$
Map({\rm M}) := \pi_0(Diff({\rm M})) = (Diff({\rm M}) / Diff^0({\rm M})),
$$

where $Diff^0({\rm M})$, the connected component of the identity, is 
the subgroup of
diffeomorphisms of
$M$ {\bf isotopic} to the identity (i.e., they are connected to the identity by
a path in $Diff({\rm M})$).
\end{df}

\begin{oss}
If $M$ is oriented then we often tacitly take $Diff^+({\rm M})$, the group
of orientation preserving diffeomorphisms of $M$ instead of $Diff({\rm
   M})$, in the definition of the mapping class group. But it is more accurate
to distinguish in this case  between $Map^+({\rm M})$ and $Map({\rm M})$.

If $M$ is a compact complex curve of genus $g$, then its mapping class group is
denoted by $Map_g$. The representation of $M = C_g$ as the $K 
(\pi,1)$ space $\HH / \Pi_g$,
i.e., as  a quotient of the (contractible) upper halfplane $\HH$ by the
free action of a Fuchsian group isomorphic to $ \Pi_g \cong \pi_1 
(C_g)$, immediately
yields the isomorphism $Map_g \cong Out ( \Pi_g) = Aut ( \Pi_g) / Int 
( \Pi_g).$

In this way the orbifold exact sequences considered in the previous lecture
$$
1 \rightarrow \Pi_{g_1}   \rightarrow  \pi_1^{orb}
\rightarrow G \rightarrow 1
$$
determine the topological action of $G$ since the homomorphism $ G \ra Map_g$
is obtain by considering, for $ g \in G$, the automorphisms obtained 
via  conjugation
by a lift $\tilde{g} \in \pi_1^{orb}$ of $g$.
\end{oss}

The relation between Artin's and Dehn's definition is the following:
\begin{teo}\label{artin}
The braid group $\mathcal{B}_n$ is isomorphic to the group
$$
\pi_0 (Map^{\infty} (\C \backslash \{1, \dots n\})),
$$

where $Map^{\infty} (\C \backslash \{1, \dots n\})$ is the group of
  diffeomorphisms which are the identity outside the disk with
center $0$ and radius $2n$.
\end{teo}

In this way Artin's standard generators $\sigma_i$ of
  $\mathcal{B}_n$ ($ i=1, \dots n-1$) can be represented by the 
so-called half-twists.

\begin{df}
The {\em half-twist} $\sigma_j$ is the diffeomorphism of $\C  \backslash \{1,
\dots n\}$ isotopic to the homeomorphism given by:

- rotation of $180$ degrees on the disk with center $ j +\frac{1}{2}$
and radius $\frac{1}{2}$,

- on a circle with the same center and radius $\frac{2+t}{4}$ the
  map $\sigma_j$ is the identity if $t\geq 1$ and rotation of
$180 (1-t)$ degrees, if $t\leq 1$.

\end{df}

\begin{figure}[htbp]
\begin{center}
\scalebox{1}{\includegraphics{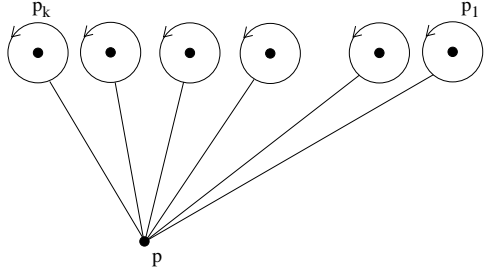}}
\end{center}
\caption{A geometric basis of $\pi_1( \C - \{1, \dots n\})$}
\label{figura2}
\end{figure}

Now, it is obvious from theorem \ref{artin} that $\mathcal{B}_n$ acts 
on the free group
$\pi_1( \C \backslash \{1, \dots n\})$, which has a
geometric basis (we take as base point the complex number $p:=
-2ni$)
$\gamma_1, \dots \gamma_n$ as illustrated in figure \ref{figura2}.

This action is called the {\em Hurwitz action of the braid group}
  and has the following algebraic description

\begin{itemize}
\item
$\sigma_i (\gamma_i) = \gamma_{i+1}$
\item
$\sigma_i (\gamma_i \gamma_{i+1}) = \gamma_i \gamma_{i+1}$, whence
$\sigma_i (\gamma_{i+1}) = \gamma_{i+1} ^{-1}\gamma_i \gamma_{i+1}$
\item
$\sigma_i (\gamma_j ) = \gamma_j$ for $j \neq i, i+1$.
\end{itemize}

Observe that the product
$\gamma_1  \gamma_2 \dots    \gamma_n$ is left invariant
  under this action.

\begin{df}
Let us consider a group $G$ and its cartesian product $G^n$.
  The map associating to each $(g_1, g_2, \dots ,g_n)$ the
product $g := g_1 g_2 \dots ,g_n \in G$ gives a partition of
$G^n$, whose subsets are called {\em factorizations} of an element
$g \in G$.

$\mathcal{B}_n$ acts on $G^n$ leaving invariant the partition,
and its orbits are called the {\bf Hurwitz equivalence classes of 
factorizations}.
\end{df}

We shall use the following notation for a factorization:
$ g_1 \circ g_2  \circ \dots  \circ g_n $, which should be carefully
distinguished from the product $g_1 g_2 \dots g_n$,
which yields an element of  $G$.

\begin{oss}
A broader equivalence relation for the set of factorizations
is obtained considering the equivalence relation
generated by Hurwitz equivalence and by
  {\bf simultaneous conjugation}. The latter,  using the following  notation
  $ a_b : = b^{-1} a b$, corresponds to the  action of $G$
on $G^n$ which carries  $ g_1 \circ g_2  \circ \dots  \circ g_n $ to
$ (g_1)_b \circ (g_2)_b  \circ \dots  \circ (g_n)_b$.

Observe that the latter action carries a factorization of  $g$ to a
factorization of the conjugate $g_b$ of $g$, hence we get equivalence
classes of factorizations for conjugacy classes of elements of $G$.
\end{oss}

The above equivalence relation plays an important role in several questions
concerning plane curves and algebraic surfaces, as we shall soon see.

Let us proceed for the meantime considering  another
interesting relation between the braid groups and the
Mapping class groups.

\begin{figure}[htbp]
\begin{center}
\scalebox{1}{\includegraphics{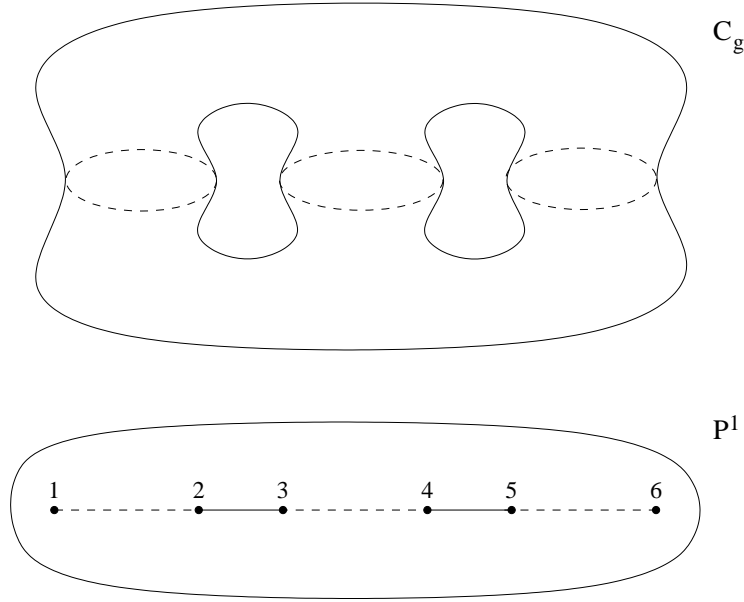}}
\end{center}
\caption{ Hyperelliptic curve of genus $2$}
\label{figura4}
\end{figure}

This relation is based on the topological  model
  provided by the hyperelliptic curve $C_g$ of equation
$$ w^2 = \prod_{i=1}^{2g + 2}  ( z - i)  $$
(see Figure \ref{figura4}  describing a   hyperelliptic curve of genus $g=2$).

Observe that, if $Y$ is the double unramified
covering of $( \PP^1 - \{1, \dots 2g+2\})$,   inverse image of $( 
\PP^1 - \{1, \dots
2g+2\})$ in $C_g$, $C_g$ is the natural compactification of $Y$ 
obtained by adding to
$Y$ the {\em ends} of $Y$
(i.e., in such a compactification one adds to $Y$  the following $ 
lim_{ K \subset
\subset Y}
\pi_0 (Y - K)$).

This description makes it clear that every homeomorphism of $( \PP^1 
- \{1, \dots
2g+2\})$ which leaves invariant the subgroup associated to
the covering $Y$  admits a lifting to a homeomorphism of $Y$,
whence also to a homeomorphism of its natural compactification $C_g$.

Such a lifting is not unique, since we can always compose with the
nontrivial automorphism of the covering.

We obtain in this way a central extension

$$ 1 \ra \Z /2 =  <H> \ra \MM ap_g^h  \ra \MM ap_{0,2g+2} \ra 1 $$

where
\begin{itemize}
\item
$H$ is the hyperelliptic involution  $ w \ra - w$
  (the
nontrivial automorphism of the covering)
\item
  $\MM ap_{0,2g+2}$ is the Dehn group of $( \PP^1 - \{1, \dots 2g+2\})$
\item
$\MM ap_g^h$ is called the hyperelliptic subgroup of the
mapping class group $\MM ap_g$, which consists of all the possible
liftings.

If  $ g \geq 3$, it is a proper subgroup of $\MM ap_g$.

\end{itemize}

While Artin's braid group
   $\BB_{2g+2}$ has the following  presentation:
$$ \langle \sigma_1,  \dots \sigma_{2g+1} |
  \sigma_i  \sigma_j =  \sigma_j
  \sigma_i  \ for  |i-j| \geq 2, \ \sigma_i  \sigma_{i+1}  \sigma_i =
\sigma_{i+1}  \sigma_i \sigma_{i+1} \rangle, $$

 Dehn's group of $( \PP^1 - \{1, \dots 2g+2\})$  $\MM
ap_{0,2g+2}$ has the presentation:
$$ \langle \sigma_1,  \dots \sigma_{2g+1} |   \sigma_1  \dots \sigma_{2g+1}
  \sigma_{2g+1} \dots  \sigma_1 = 1 , (\sigma_1  \dots 
\sigma_{2g+1})^{2g+2} =1,$$ $$
  \sigma_i  \sigma_j =  \sigma_j
  \sigma_i  \ for  |i-j| \geq 2, \ \sigma_i  \sigma_{i+1}  \sigma_i =
\sigma_{i+1}  \sigma_i \sigma_{i+1} \rangle, $$
finally the hyperelliptic mapping class group   $ \MM ap_g^h $ has 
the presentation:

$$ \langle \xi_1,  \dots \xi_{2g+1}, H |   \xi_1  \dots \xi_{2g+1}
   \xi_{2g+1} \dots  \xi_1 = H , H^2 = 1,(\xi_1  \dots 
\xi_{2g+1})^{2g+2} =1, $$ $$
   H \xi_i = \xi_i H \ \forall i,  \xi_i  \xi_j =  \xi_j
  \xi_i  \ for  |i-j| \geq 2, \ \xi_i  \xi_{i+1}  \xi_i =
\xi_{i+1}  \xi_i \xi_{i+1} \rangle. $$

We want to illustrate the geometry underlying these important formulae.
Observe that $\sigma_j$ yields a homeomorphism of the disk $U$ with 
centre $ j + 1/2$
and radius 3/4, which permutes the two points $j, j+1$.

\begin{figure}[htbp]
\begin{center}
\scalebox{1}{\includegraphics{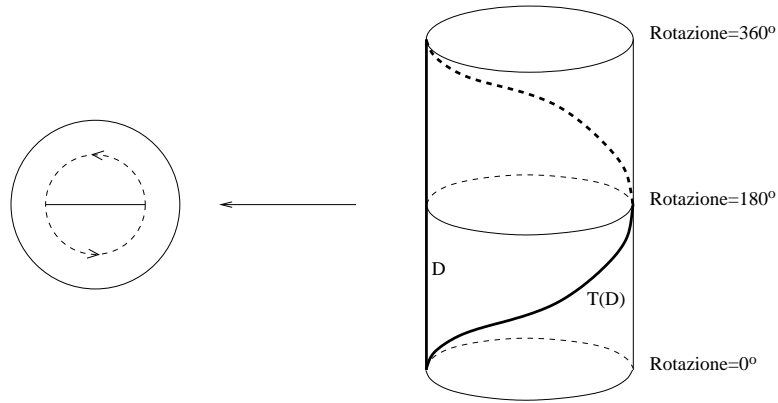}}
\end{center}
\caption{At the left, a half twist; at the right:its lift,
the Dehn-Twist-$T$, and its action on the segment $D$}
\label{figura5}
\end{figure}

Therefore there are two liftings of $\sigma_j$ to homeomorphisms of the
inverse image $V$ of $U$ in $C_g$: one defines then $\xi_j$ as the one
of the two liftings which acts as the identity on the boundary $\partial V$,
which is a union of two loops
(see figure \ref{figura5}).

$\xi_j$ is called the {\bf  Dehn twist} and
corresponds geometrically to the   diffeomorphism of a truncated
cylinder which is the identity on the boundary, a rotation by 180 
degrees on the
equator, and on each parallel at height   $t$ is a  rotation by $ t \  360 $
degrees ( where
$t \in [0,1]$).

One can define in the same way a Dehn twist for each  loop in
$C_g$ (i.e., a subvariety diffeomorphic to $S^1$):

\begin{df}
Let $C$ be an oriented Riemann surface. Then a {\em positive Dehn twist
   $T_{\alpha}$} with respect to a simple closed curve $\alpha$ on $C$
is an isotopy class of a diffeomorphism $h$ of $C$ which is equal to
the identity outside a neighbourhood of $\alpha$ orientedly 
homeomorphic to an annulus in the plane, while inside
the annulus $h$ rotates the inner boundary of the annulus by $360$ degrees to
the right and damps the rotation down to the identity at the outer
boundary.
\end{df}

Dehn's fundamental result (\cite{dehn}) was the following

\begin{teo}
The mapping class group  $\MM ap_g$ is generated by Dehn twists.
\end{teo}

Explicit presentations of $\MM ap_g$ have later been given by Hatcher and
Thurston (\cite{h-t}),and an improvement of the method lead
to the simplest available presentation, due to Wajnryb
(\cite{wajnryb}, see also \cite{waj2}).

We shall see in the next subsection how the Dehn twists are related 
to the theory of
Lefschetz fibrations.

\subsection{Lefschetz pencils and Lefschetz fibrations}

The method introduced by Lefschetz for the study of the topology of
algebraic varieties is the topological analogue of the method of
hyperplane sections and projections of the classical italian algebraic
geometers.

An excellent exposition of the theory of Lefschetz pencils is the
article by Andreotti and Frankel \cite{a-f}, that we try to briefly 
summarize here.

Let $X \subset \PP^N$ be  projective variety, which for simplicity we 
assume to be
smooth, and let
$L \cong \PP^{N-2} \subset \PP^N$ be a general linear subspace of 
codimension 2. $L$ is
the base locus of a pencil of hyperplanes $H_t, t \in \PP^1$, and the 
indeterminacy
locus of a rational map $\phi : \PP^N \setminus L \ra \PP^1$.

The intersection $Z : X \cap L$ is smooth, and the blow up of $X$ 
with centre $Z$
yields a smooth variety $X'$ with a morphism $f \colon X' \ra \PP^1$ whose
fibres are isomorphic to the hyperplane sections $Y_t := X \cap H_t$, while
the exceptional divisor is isomorphic to the product $Z \times \PP^1$ and on it
the morphism $f$ corresponds to the second projection.

\begin{df}
The dual variety $W^{\vee} \subset {\PP^N}^{\vee}$ of a projective 
variety $W$ is defined
as the closure of  the set of hyperplanes which contain the tangent 
space $TW_p$ at a
smooth point $p \in W$.
A pencil of hyperplanes $H_t, t \in \PP^1$, is said to be a {\bf 
Lefschetz pencil}
if the line $L'$ dual to the subspace $L$

1) does not intersect $W^{\vee}$ if $W^{\vee}$ is not a hypersurface

2) intersects $W^{\vee}$ transversally in $ \mu : = deg (W^{\vee})$ 
points otherwise.
\end{df}

An important theorem is the

{\bf Biduality theorem:} {\em $ (W^{\vee})^{\vee} = W$ .}

It follows from the above theorem and the previous definition that if 
$W^{\vee}$ is not
a hypersurface,
$f$ is a differentiable fibre bundle, while in case 2) all the fibres are
smooth, except $\mu$ fibres which correspond to tangent hyperplanes $H_{t_j}$.
And for these $Y_{t_j}$ has only one singular point $p_j$, which has
an ordinary quadratic singularity as a hypersurface in $X$ (i.e., there are
local holomorphic coordinates $(z_1, \dots z_n) $ for $X$ such that
locally at $p_h$ $$Y_{t_h} = \{ z| \sum_j z_j^2 = 0 \}).$$
Writing $z_j = u_j + i v_j$, the equation  $\sum_j z_j^2 = \rho$ for $\rho \in
\R$ reads out as $\sum_j u_j v_j = 0,\sum_j (u_j^2 - v_j^2)= \rho$. In vector
notation, and assuming $ \rho \in \R_{\geq0}$, we may rewrite as
$$ \langle u,v \rangle = 0 , | u|^2 = \rho + | v|^2 .$$

\begin{df}
The {\bf vanishing cycle} is the sphere $\Sigma_{t_h + \rho} $ of 
$Y_{t_h + \rho} $
given, for  $ \rho
\in \R_{ > 0}$, by $\{ u + i v | \ | u|^2 = \rho , v = 0 \}.$
\end{df}

The normal bundle of the vanishing cycle $\Sigma_t$ in $Y_t$ is easily seen,
in view of the above equations,
to be isomorphic to the tangent bundle to the sphere $S^{n-1}$, 
whence we can identify
a tubular neighbourhood of $\Sigma_t$ in $Y_t$ to the unit ball in
the tangent bundle of the sphere $S^{n-1}$.
We follow now the definition given in \cite{kas} of the corresponding Dehn twist.

\begin{df}
Identify the  sphere $\Sigma = S^{n-1} = \{ u |  |u|= 1\}$ to the zero section
of its  unit tangent bundle $Y  = \{ (u , v) | \langle u,v \rangle = 
0 , |u|= 1,
  | v| \leq 1 \}.$

Then the Dehn twist $T := T_{\Sigma}$ is the diffeomorphism of $Y$ such that,
if we let $\ga_{u,v} (t)$ be the geodesic on $S^{n-1}$ with initial point $u$,
initial velocity $v$, then
$$ T (u,v) : = - (  \ga_{u,v} (\pi |v|), \frac{d}{dt}  \ga_{u,v} (\pi |v|)). $$

We have then:
1) $T$ is the antipodal map on $\Sigma$

2) $T$ is the identity on the boundary  $\partial Y  = \{ (u , v) | \langle u,v
\rangle = 0 , |u|= 1 =
  | v| \}.$
\end{df}

One has  the

{\bf Picard- Lefschetz Theorem} {\em  The Dehn twist $T$ is the local monodromy
of the family $Y_t$ ( given by the level sets of the function $\sum_j 
z_j^2 $).}

Moreover, by the classical Ehresmann theorem, one sees that a singular fibre
$Y_{t_j}$ is obtained from a smooth fibre by substituting a neighbourhood of
the vanishing cycle $\Sigma$ with the contractible intersection of the complex
quadratic cone$\sum_j z_j^2 = 0$ with a ball around $p_j$. Hence

\begin{teo}
{\bf (Generalized Zeuthen Segre formula)} The number $\mu$ of 
singular fibres in a
Lefschetz pencil, i.e., the degree of the dual variety $X^{\vee}$, is 
expressed as a
sum of topological Euler numbers
$$ e (X) + e (Z) = 2 e(Y) + (-1)^n
\mu  ,$$
where $Y$ is a smooth hyperplane section, and $ Z = L \cap X$ 
is the base locus of the pencil.
\end{teo}
\Proof {\em (idea)}
Replacing $Z$ by $Z \times \PP^1$ we see that we replace $ e(Z)$ by $ 2 e (Z)$,
hence the left hand side expresses the Euler number of the blow up $X'$.

This number can however be computed from the mapping $f$: since the 
Euler number is
multiplicative for fibre bundles, we would have that this number were $ 2 e(Y)$
if there were no singular fibre. Since however for each singular 
fibre we replace
something homotopically equivalent to the sphere $ S^{n-1}$ by a 
contractible set,
we have to subtract $ (-1)^{n-1}$ for each singular fibre.

\QED

Lefschetz pencils were classically used to describe the homotopy
and homology groups of algebraic varieties.

The main point is that the finite part of $X'$, i.e., $X' - Y_{\infty}$,
has the socalled 'Lefschetz spine' as homotopy retract.

In order to explain what this means, assume, without loss of 
generality from the
differentiable viewpoint, that the fibres
$Y_0$ and $Y_{\infty}$ are smooth fibres, and that the singular 
fibres occur for
some roots of unity $t_j$, which we can order in counterclockwise order.

\begin{df}
Notation being as before,  define the relative vanishing cycle $\Delta_j$
as the union, over $t$ in the segment $[0, t_j]$, of the vanishing cycles
$\Sigma_{t,j}$: these are defined, for $t$ far away from $t_j$, using 
a trivialization
of the fibre bundle obtained restricting $ f \colon X' \ra \PP^1$ to the
half open segment $[0, t_j)$.

The {\bf  Lefschetz spine} of the Lefschetz pencil is the union of the fibre $Y_0$
with the $\mu$ relative vanishing cycles $\De_j$.
\end{df}

\begin{teo} {\bf (Lefschetz' theorems I and II)}
1) The Lefschetz spine is a deformation retract of $X' -  Y_{\infty}$.

2) The affine part $X -  Y_{\infty}$ has the homotopy type of a cell
complex of dimension $n$.

3) The inclusion $\iota \colon  Y_0 \ra X$ induces homology homomorphisms
$ H_i (\iota ) \colon H_i (Y_0, \Z) \to H_i (X, \Z) $ which are

3i) bijective for $ i < n-1$

3ii) surjective if $ i = n-1$; moreover

4) the kernel of $ H_{n-1} (\iota ) $ is generated by the vanishing cycles,
i.e., by the images of
$  H_{n-1} (\Sigma_{0,j}, \Z)$.
\end{teo}

{\em Comment on the proof:}

1) follows by using the Ehresmann's theorem outside of the singularities,
and by retracting locally a neighbourhood of the singularities partly 
on a smooth fibre
$ Y_t$, with $ t \in (0, t_j)$, and partly on the union of the 
vanishing cycles.
Then one goes back all the way to $Y_0$.

For 2) we simply observe that  $X -  Y_{\infty}$ has $(Y_0 
\setminus Z )\cup (\cup_j
\De_j)$ as deformation retract. Hence, it is homotopically equivalent 
to a cell complex
obtained by attaching $\mu$ $n$-cells to $Y_0 \setminus Z $, and 2) 
follows then by
induction on $n$.

3) and 4) are  more delicate and require some diagram chasing,  which
can be found in \cite{a-f}, and which we do not reproduce here.

\qed

In the 70's Moishezon and Kas realized (see e.g. \cite{moi77} and \cite{kas}),
after the work of Smale on the smoothing of handle attachments, that Lefschetz
fibrations could be used to investigate the differential topology of algebraic
varieties, and especially of algebraic surfaces.

For instance, they give a theoretical method, which we shall now 
explain,  for the
extremely difficult problem to decide whether two  algebraic surfaces 
which are not
deformation equivalent are in fact diffeomorphic (\cite{kas}).

\begin{df}
Let $M$ be a compact differentiable (or even symplectic) manifold of
real even dimension $2n$

A {\em Lefschetz fibration} is a differentiable map $ f : M \ra
\PP^1_{\C}$ which

a) is of maximal rank except for a finite number of critical
points $p_1, \dots p_m$ which have distinct critical values $b_1,
\dots b_m \in \PP^1_{\C}$,

b)  has the property that around $p_i$
there are complex coordinates $(z_1, \dots z_n) \in \C^n$  such that locally
$f =\sum_j z_j^2  + const.$
(in the symplectic case, we require the given coordinates to be
Darboux coordinates, i.e., such that the symplectic form
$\omega$ of $M$  corresponds to the natural constant coefficients symplectic
structure on $\C^n$).
\end{df}

\begin{oss}
1) A similar definition can be given if $M$ is a manifold with boundary,
replacing $\PP^1_{\C}$ by a disc $D \subset \C$.

2) An important theorem of Donaldson (\cite{donsympl}) asserts that
    for symplectic manifolds there exists (as for the case of projective
    manifolds) a Lefschetz pencil, i.e., a Lefschetz fibration $ f : M'
    \ra \PP^1_{\C}$ on a symplectic blow up $M'$ of $M$ (see \cite{mcduffsal}
for the definition of symplectic blow-up).

3) A Lefschetz fibration with smooth fibre $F_0 = f^{-1} (b_0)$  and 
with critical
values
$b_1, \dots b_m \in
\PP^1_{\C}$, once  a geometric basis
$\gamma_1 , \gamma_2 ,  \dots , \gamma_m$ of
$ \pi_1 ( \PP^1_{\C} \backslash \{ b_1, \dots , b_m \}, b_0)$ is
chosen, determines
a factorization of the identity in the mapping class group $\MM ap (F_0)$

$$
\tau_1 \circ \tau_2 \circ \dots \circ \tau_m = Id
$$
as a product of Dehn twists.

4) Assume further that
$b_0, b_1, \dots b_m \in \C = \PP^1 \setminus \{ \infty \}$: then the Lefschetz
fibration determines also a homotopy  class of an arc $\la  $  between $
\tau_1  \tau_2  \dots  \tau_m $ and the identity in $ Diff^0 (F_0) $. This
class is trivial when $F_0 = C_g$, a compact Riemann surface of genus
$ g \geq 1$.

5) More precisely, the Lefschetz fibration $f$ determines isotopy
classes of embeddings $\phi_j : S^{n-1} \ra F_0$ and of bundle isomorphisms
$\psi_j $ between the tangent bundle of $ S^{n-1}$ and the normal bundle
of the embedding $\phi_j$; $\tau_j$ corresponds then to the Dehn twist
for the embedding $\phi_j$.
\end{oss}

We are now ready to state the theorem of Kas (cf. \cite{kas}).

\begin{teo}
Two Lefschetz fibrations $(M, f)$,
$(M', f')$ are equivalent (i.e., there are two diffeomorphisms $u : M
\ra M', v : \PP^1 \ra \PP^1$ such that $ f' \ \circ u = v \circ f$) if
and only if

1) the corresponding invariants $$(\phi_1, \dots \phi_m),(\psi_1, 
\dots \psi_m);
(\phi '_1, \dots \phi '_m) ,(\psi '_1, \dots \psi '_m)$$ correspond to 
each other
via a diffeomorphism of $F_0$ and a diffeomorphism $v$ of $\PP^1$.
 This implies in particular

1') the two corresponding factorizations of the identity in
the mapping class group are equivalent (under the equivalence relation
generated by  Hurwitz equivalence and by simultaneous conjugation).

2) the respective homotopy classes  $\la, \la '$ correspond to each
other under the above equivalence.

Conversely, given  $(\phi_1, \dots \phi_m) (\psi_1, \dots \psi_m)$ such that
the corresponding Dehn twists $\tau_1 , \tau_2 ,\dots  \tau_m $
yield a factorization of the identity, and given a homotopy class $\la$
of a path connecting $\tau_1 \tau_2  \dots  \tau_m$ to the identity 
in $Diff (F_0)$,
there exists an associated  Lefschetz fibration.

If the fibre $F_0$  is a Riemann surface of genus $ g \geq 2$ then the
Lefschetz fibration is uniquely determined by the equivalence
class of a factorization of
the identity $$
\tau_1 \circ \tau_2 \circ \dots \circ \tau_m = Id
$$
as a product of Dehn twists.
\end{teo}

\begin{oss}
1) A similar result holds for Lefschetz fibrations over the disc and
    we get a factorization
$$
\tau_1 \circ \tau_2 \circ \dots \circ \tau_m = \phi
$$
of the monodromy $\phi$ of the fibration over the boundary of the disc $D$.

2) A Lefschetz fibration with fibre $C_g$ admits a symplectic 
structure if
 each Dehn twist in the factorization is positively oriented
(see  section 2 of \cite{abkp}).
\end{oss}

Assume that we are given two Lefschetz fibrations over $\PP^1_{\C}$:
then we can consider the fibre sum of these two fibrations, which depends
as we saw on a diffeomorphism chosen between two respective smooth fibers
(cf. \cite{g-s} for more details).

This operation translates (in view of the above quoted theorem of Kas)
  into the following definition of ``conjugated composition'' of factorization:

\begin{df}
Let $\tau_1 \circ \tau_2 \circ \dots \circ \tau_m = \phi $ and
$\tau'_1 \circ \tau'_2 \circ \dots \circ \tau'_r = \phi' $
be two factorizations: then their {\em composition conjugated by $\psi$  }
is the factorization
$$
\tau_1 \circ \tau_2 \circ \dots \tau_m \circ
   (\tau'_1 )_{\psi}\circ (\tau'_2)_{\psi} \circ \dots \circ
(\tau'_r)_{\psi}= \phi \circ (\phi')_{\psi}.
$$
\end{df}

\begin{oss}
1) If  ${\psi}$ and $\phi'$ commute, we obtain a factorization of $\phi \phi'$.

2) A particular case is the one where $\phi = \phi' = id$ and it corresponds to
    Lefschetz fibrations over $\PP^1$.
\end{oss}

No matter how beautiful the above results are, for a general $X$ projective
or $M$ symplectic, one has Lefschetz pencils, and not Lefschetz fibrations,
and a natural question is to which extent the surgery corresponding 
to the blowup does
indeed simplify the differentiable structure of the manifold. In the 
next subsection we
shall consider  results by Moishezon somehow related to this question.

\subsection{Simply connected algebraic surfaces: topology versus 
differential topology}

In the case of compact  topological manifolds of real dimension 4 the methods
of Morse theory and of simplification of cobordisms turned out to encounter
overwhelming difficulties, and only in 1982 M. Freedman (\cite{free}), using
new ideas in order to show the (topological) triviality of certain handles
introduced by Casson, was able to obtain a complete classification
of the simply connected compact  topological 4-manifolds.

Let $M$ be such a manifold, fix an orientation of $M$, and let
$$ q_M \colon H_2(M, \Z) \times H_2(M, \Z) \to \Z$$
be the intersection form, which is unimodular by Poincar\'e duality.

\begin{teo} {\bf Freedman's theorem}
Let $M$ be an oriented compact simply connected topological manifold: 
then $M$ is
determined by its intersection form and by the Kirby-Siebenmann invariant
$\alpha (M) \in \Z / 2$, which vanishes if and only if $ M \times [0,1]$
admits a differentiable structure.
\end{teo}

The basic invariants of $q_m$ are its signature $ \sigma (M) : = b^+ 
(M) - b^- (M)$,
and its parity ($q_m$ is said to be {\bf even} iff $q_m (x,x) \equiv 
0 \ ( mod 2)$ $
\forall x \in H_2(M, \Z)$).

A basic result by Serre (\cite{serre}) says that if $q_M$ is indefinite
then it is determined by its rank, signature and parity.

The corollary of Freedman's theorem for complex surfaces is the following

\begin{teo}\label{complex}
Let $S$ be a compact simply connected complex surface,
and let  $r$ be the divisibility index of the  canonical  class
$c_1( K_X) \in H^2 (X, \Z)$.

$S$ is said to be EVEN if $q_S$ is EVEN, and this holds  iff $ r 
\equiv 0(mod 2)$, else
$S$ is said to be ODD. Then
\begin{itemize}
\item
(EVEN ) If  $S$ is EVEN, then  $S$ is topologically a connected sum 
of copies of
  $ \PP ^1_{\C} \times \PP ^1_{\C}$ and of a  $K3$ surface
if the signature of the intersection form is negative, and of copies of
$ \PP ^1_{\C} \times \PP ^1_{\C}$ and of a  $K3$ surface with
opposed orientation in the case where the signature is positive.
\item
(ODD) $S$ is ODD: then  $S$ is topologically a connected sum of copies of
$\PP ^2_{\C} $ and of ${\PP ^2_{\C}}^{opp} .$
\end{itemize}
\end{teo}

\Proof $S$ has a differentiable structure, whence $\alpha (S) = 0$,
and the corollary follows from Serre's result if the intersection 
form is indefinite.

We shall now show that if the intersection form is definite, 
then $ S \cong \PP^2_{\C}$. Observe that $q=0$,
since
$S$ is simply connected, and therefore $ b^+ (S) = 2 p_g + 1$,
in particular  the
intersection form is positive, $b_2 = 2 p_g + 1$,
hence $ e (S) = 2 \chi (S) + 1$, and $K^2_S = 10 \chi (S) -1$
by Noether's formula.

By the Yau Miyaoka inequality $q=0$ implies $K^2_S \leq 9 \chi (S) $,
whence $\chi (S) \leq 1 $ and $ p_g = 0$. 

Therefore $ \chi(S) = 1$, and $ K_S^2 = 9$. Applying again Yau's theorem (\cite{yau})
we see that  $ S = \PP^2 _{\C}$. In fact, if   $S$ were of general type 
its universal cover would be the unit ball in $\C^2$, 
contradicting simple connectivity.

\QED

\begin{oss}
${\PP ^2_{\C}}^{opp} $ is the manifold
${\PP ^2_{\C}}$ with opposed  orientation.

A K3 surface is  (as we already mentioned) a surface  $S$ orientedly
diffeomorphic to a nonsingular surface
$X$ of degree 4 in $\PP^3_{\C}$,
for instance
$$ X = \{ (x_0,x_1,x_2,x_3)
\in \PP^3_{\C} | x_0^4 + x_1^4 + x_2^4 + x_3^4 = 0\}.$$
(by a theorem of Kodaira, cf. \cite{KodIII},  $S$ is also deformation
equivalent to such a surface
$X$).
\end{oss}

Not only $\PP^2_{\C}$ is the only algebraic surface with a definite 
intersection form,
but Donaldson showed that a result of a similar flavour holds for differentiable
manifolds, i.e., if we have a positive definite intersection form, then
we have topologically a connected sum of copies of
$\PP ^2_{\C} $.

There are several restrictions for the intersection forms of differentiable
manifolds, the oldest one being Rokhlin's theorem stating that  the 
intersection form
in the even case is divisible by 16. Donaldson gave other restrictions
for the intersection forms of differentiable 4-manifolds ( see \cite{d-k}), but
the socalled 11/8 conjecture is still unproven: it states that if the 
intersection form
is even, then we have topologically a connected sum as in the case (EVEN)
of theorem \ref{complex}.

More important is the fact that  Donaldson's theory has made clear
in the 80's  (\cite{don1},
\cite{don2},
\cite{don3},\cite{don4}) how drastically homeomorphism and
diffeomorphism differ in dimension 4,
and especially for algebraic surfaces.

Later on, the  Seiberg-Witten theory showed with simpler methods the 
following result
(cf. \cite{witten} o \cite{mor}):

\begin{teo}
Any diffeomorphism between minimal surfaces (a fortiori, an even 
surface is necessarily
  minimal) $S, S'$
carries
  $c_1( K_S)$ either to  $c_1( K_{S'})$ or to  $- c_1( K_{S'})$
\end{teo}

\begin{cor}
The divisibility index  $r$ of the  canonical  class
$c_1( K_S) \in H^2 (S, \Z)$ is a differentiable invariant of $S$.
\end{cor}

Since only the parity $ r (mod 2)$ of the canonical class is a 
topological invariant
it is then not difficult to construct examples of simply connected algebraic
surfaces  which are homeomorphic but not diffeomorphic
  (see \cite{cat3}).

Let us illustrate these examples, obtained as simple bidouble covers of
$\PP^1 \times \PP^1$.

These surfaces are contained in the geometric vector bundle whose sheaf
of holomorphic sections is $ \hol_{\PP^1 \times \PP^1} (a,b) \bigoplus
\hol_{\PP^1 \times \PP^1} (c,d) $ and are
described there by the following two equations:
$$
z^2 = f(x,y),
$$
$$
w^2 = g(x,y),
$$
where $f$ and $g$ are bihomogeneous polynomials of respective 
bidegrees $(2a,2b)$,
$(2c,2d)$ ($f$ is a section of  $ \hol_{\PP^1 \times \PP^1} (2a,2b)$, $g$
is a section of  $ \hol_{\PP^1 \times \PP^1} (2c,2d)$).

These Galois
  covers of $\PP^1 \times \PP^1$, with Galois group $ (\Z / 2 \Z )^2$,
are smooth if and only if the two curves
$C : = \{ f=0 \}$ and $D : = \{ g=0 \}$ in
$\PP^1 \times \PP^1$ are smooth and intersect transversally.

The holomorphic invariants can be easily calculated, since, if $ p 
\colon  X \ra
\PP:= \PP^1 \times \PP^1$ is the finite Galois cover, then
$$ p_* \hol_X \cong  \hol_{\PP^1 \times \PP^1} \oplus  z  \hol_{\PP^1 
\times \PP^1}
(-a, -b) \oplus  w  \hol_{\PP^1 \times \PP^1}
(-c, -d) \oplus  zw  \hol_{\PP^1 \times \PP^1}
(-a - c, -b -d).$$

Hence  $ h^1 (\hol_X) = 0$, whereas $ h^2 (\hol_X) = (a-1)(b-1) +(c-1)(d-1) +
(a+c-1)(b+d-1) $.
Assume that $X$ is smooth: then the ramification formula yields
$$ \hol_X (K_X) =  \hol_X (p^*  K_{\PP^1 \times \PP^1} + R) = p^*
( \hol_{\PP^1 \times \PP^1} (a+c-2,b+d-2)) $$
since $ R = div (z) + div (w)$. In particular, $K_X^2 = 8 (a+c-2)(b+d-2)  $
and the holomorphic invariants of  such coverings depend only upon
the numbers $(a+b-2)(c+d-2)$ and $ab + cd$.

\begin{teo}
Let $S, S' $ be  smooth bidouble covers of $\PP^1 \times \PP^1$ of 
respective types
$ (a,b)(c,d), (a',b')(c',d')$.

Then $S$ is of general type for $ a+c \geq 3 , b+d \geq 3$, and is 
simply connected.
Moreover, the divisibility $r (S)$ of the canonical class $ K_S$ is equal to
$ G.C.D. ( (a+c-2),(b+d-2))$.

$S$ and $S'$ are (orientedly) homeomorphic if and only if $r(S) 
\equiv r (S') ( mod 2)
$  and
$$  (a+b-2)(c+d-2) = (a'+b'-2)(c'+d'-2) \ {\rm
and}\ ab + cd = a'b' + c'd'.$$
$S$ and $S'$ are not diffeomorphic if  $r(S) \neq r(S')$, and for each integer
$h$, we can find such surfaces $S_1, \dots S_h$ which are pairwise homeomorphic
but not diffeomorphic.
\end{teo}

{\em Idea of the proof.}
Set for simplicity $ u :=  (a+c-2), v : = (b+d-2)$  so that $\hol_S (K_S )= p^*
( \hol_{\PP^1 \times \PP^1} (u ,v))$ is ample whenever $ u,v \geq 1$.

The property that $S$ is simply connected (cf. \cite{cat1} for 
details) follows once one
establishes that the fundamental group $ \pi_1 ( (\PP^1 \times \PP^1) 
\setminus (C \cup
D))$ is abelian. To establish that the group is abelian, since it is generated
by a product of simple geometric loops winding once around a smooth point
of $ C \cup D$, it suffices to show that these loops are central.
But this follows from considering a Lefschetz pencil having $C$ 
(respectively, $D$)
as a fibre (in fact, an $S^1$ bundle over a punctured Riemann surface 
is trivial).

Since  this group is abelian, it is generated by two elements $\ga_C, \ga_D$
which are simple geometric loops winding once around $C$, resp. $D$.
The fundamental group $ \pi_1 (S \setminus R)$ is then generated by
$ 2 \ga_C$ and $ 2 \ga_D$, but these two elements lie in the kernel
of the surjection $ \pi_1 (S \setminus R) \to  \pi_1 (S )$ and we conclude the
triviality of this latter group.

The argument for the divisibility of $K_S$ is more delicate, and we refer
to \cite{cat3} for the proof of the key lemma asserting that
$ p^* (H^2 ( \PP^1 \times \PP^1 , \Z )) = H^2 ( S, \Z)^{G} $
where $G$ is the Galois group $ G = (\Z / 2)^2 $ (the proof uses 
arguments of group
cohomology and is valid in greater generality). Thus, the divisibility of $K_S$
equals the one of $ c_1 ( \hol_{\PP^1 \times \PP^1} (u ,v))$, i.e., 
$G.C.D. (u,v)$.

Now, resorting to Freedman's theorem, it suffices to observe that rank 
and signature
of the intersection form are given by $ e(S) - 2, \sigma (S)$, and these,
as we saw in the first lecture, equal $ 12 \chi (S) - K^2_S , K^2_S - 
8 \chi (S) $.
In this case $K^2_S = 8 uv$, $ \chi (S) = u v + (ab + cd)$.

There remain to find $h$ such surfaces, and for this purpose, we use 
Bombieri's argument
(appendix to \cite{cat1}): namely, let $u'_i v'_i = 6^n$ be $h$ distinct
factorizations and, for a positive number $T$, set $ u_i : = T u'_i, 
v_i : = T v'_i$.
It is clear that $G.C.D. (u_i, v_i) = T (G.C.D. (u'_i, v'_i)) $ and 
these G.C.D.'s
are distinct since the given factorizations are distinct (as 
unordered factorizations),
and they are even integers if each $u'_i, v'_i$ is even.

It suffices to show that there are integers $w_i, z_i$ such that, setting
$ a_i : = (u_i + w_i)/2 + 1, c_i : = (u_i - w_i)/2 + 1, b_i : = (v_i 
- z_i)/2 + 1,
d_i : = (v_i + z_i)/2 + 1 $, then $ a_i b_i + c_i d_i = constant$ and
the required inequalities $ a_i, b_i, c_i, d_i \geq 3$ are verified.

This can be done by the box principle.

\qed

It is important to contrast the existence of homeomorphic but not diffeomorphic
algebraic surfaces to an important theorem  established at the 
beginning of the study of
$4$-manifolds by C.T.C. Wall  ( \cite{wall}):

\begin{teo} {\bf (C.T.C. Wall)}
Given two simply connected differentiable 4-manifolds $M, M'$ with
isomorphic intersection forms, then there exists an integer $k$ such that
the iterated connected sums
$ M \sharp k (\PP^1 \times \PP^1)$ and $ M' \sharp k (\PP^1 \times \PP^1)$
are diffeomorphic.
\end{teo}

\begin{oss}
1) If we take ${\PP ^2_{\C}}^{opp} $, i.e., $\PP^2$ with opposite orientation,
then the selfintersection of a line equals $-1$, just as for the 
exceptional curve
of a blow up. It is easy to see that blowing up a point of a smooth 
complex surface
$S$ is the same differentiable operation as taking the connected sum
$ S \ \sharp  \ {\PP ^2_{\C}}^{opp} $.

2) Recall that the blowup of the plane $\PP^2$ in two points is isomorphic to
the quadric $\PP^1 \times \PP^1$ blown up in a point. Whence, for the 
connected sum
calculus, $ M \sharp (\PP^1 \times \PP^1) \ \sharp {\PP ^2_{\C}}^{opp} \cong
M \sharp (\PP^2 ) \ \sharp  2 {\PP ^2_{\C}}^{opp}$. From Wall's theorem follows
then (consider Wall's theorem for $M^{opp}$) that for any  simply 
connected 4-manifold
$M$  there are integers
$k,p,q$ such that
  $M \sharp (k+1) (\PP^2 ) \ \sharp  (k) {\PP ^2_{\C}}^{opp} \cong
  p (\PP^2 ) \ \sharp  (q) {\PP ^2_{\C}}^{opp} $.
\end{oss}

The moral of Wall's theorem was  that homeomorphism of simply 
connected 4-manifolds
implies {\bf stable} diffeomorphism (i.e., after iterated connected 
sum with some basic
manifolds as $(\PP^1 \times \PP^1)$ or, with both $\PP^2 ,   {\PP
^2_{\C}}^{opp} $).

The natural question was then how many such connected sums were 
indeed needed, and
if there were needed at all. As we saw, the Donaldson and Seiberg 
Witten invariants
show that some connected sum is indeed needed.

Boris Moishezon, in collaboration with Mandelbaum, studied the 
question in detail
(\cite{moi77}, \cite{mm1}, \cite{mm2})  for
many concrete examples of algebraic surfaces, and gave the following

\begin{df}
A differentiable simply connected 4-manifold $M$ is {\bf completely 
decomposable}
if there are integers $p, q$ with  $M  \cong
  p (\PP^2 ) \ \sharp  (q) {\PP ^2_{\C}}^{opp} $, and { \bf almost completely
decomposable} if  $M  \sharp
   (\PP^2 )$ is completely decomposable (note that the operation 
yields a manifold
with odd intersection form, and if $M$ is an algebraic surface $\neq 
\PP^2$, then
we get an indefinite intersection form.
\end{df}

Moishezon and Mandelbaum (\cite{mm1}) proved almost complete 
decomposability for smooth
hypersurfaces in $\PP^3$, and Moishezon proved (\cite{moi77}) almost complete
decomposability for simply connected elliptic surfaces. Observe that 
rational surfaces
are obviously completely decomposable, and therefore one is only left 
with simply
connected surfaces of general type, for which as far as I know the 
question of almost
complete decomposability is still unresolved.

Donaldson's work clarified the importance of the connected sum with $\PP^2$,
showing the following results (cf. \cite{d-k} pages 26-27)

\begin{teo} {\bf( Donaldson)}
If $M_1, M_2$ are simply connected differentiable 4-manifolds with $ b^+(M_i) >
0$, then the Donaldson polynomial invariants $ q_k \in S^d (H^2 (M, 
\Z)$ are all zero
for $ M = M_1 \sharp M_2$. If instead $M$ is an algebraic surfaces, then
the Donaldson polynomials $q_k$ are $ \neq 0 $ for large $k$. In particular,
an algebraic surface cannot be diffeomorphic to a connected sum $M_1 
\sharp M_2$
with $M_1, M_2$ as above ( i.e., with $ b^+(M_i) > 0 $).

\end{teo}

\subsection{ABC surfaces}

This subsection is devoted to the diffeomorphism type of certain series of
families of bidouble covers,
   depending on 3 integer parameters (a,b,c) (cf.
\cite{cat02}, \cite{c-w}).

Let us make some elementary remark, which will be useful in order to 
understand concretely
the last part of the forthcoming definition.

Consider the projective line $\PP^1$ with homogeneous coordinates $(x_0, x_1)$
and with nonhomogeneous coordinate $ x : = x_1 / x_0$..
Then the homogeneous polynomials of degree $m$  $ F (x_0,x_1)$ are exactly
the space of holomorphic sections of $\hol_{\PP^1} (m)$: in fact to such an $F$
corresponds the pair of holomorphic functions $f_0 (x) : = \frac{F 
(x_0,x_1)}{x_0^m}$
on $U_0 : = \PP^1 \setminus \{\infty\}$, and $f_1 (1/x) : = \frac{F 
(x_0,x_1)}{x_1^m}$
on $U_1 : = \PP^1 \setminus \{ 0 \}$. They satisfy the cocycle condition
$ f_0 (x)  x^m = f_1 (1/x) $.

We assumed here $m$ to be a positive integer, because  $\hol_{\PP^1} (- m)$ has
no holomorphic sections, if $ m > 0$. On the other hand, sheaf theory 
(the exponential sequence
and the partition of unity argument) teaches us that the cocycle $x^{-m}$ for
  $\hol_{\PP^1} (- m)$ is cohomologous, if we use differentiable functions,
to $\bar{x}^m$ (indeed $x^{-m} = \frac{ \bar{x}^m}{|x|^{2m} }$, a formula which
hints at  the homotopy $ \frac{ \bar{x}^m}{|x|^{2m t} }$ of
the two cocycles).

This shows in particular that the  polynomials  $ F (\bar{x}_0,\bar{x}_1)$
which are homogeneous  of degree $m$ are differentiable sections of 
$\hol_{\PP^1} (- m)$.

Since sometimes we shall need to multiply together sections of 
$\hol_{\PP^1} (- m)$
with sections of $\hol_{\PP^1} ( m)$,and get a global function, we 
need the cocycles
to be the inverses of each other. This is not a big problem, since on 
a circle of radius
$R$ we have $ \bar{x} x = R^2$. Hence to a polynomial $ F 
(\bar{x}_0,\bar{x}_1)$
we associate the two functions
$$f_0 (\bar{x}) : =  \frac{F (\bar{x}_0,\bar{x}_1)}{\bar{x}_0^m} \ 
{\rm on} \  \{ x | |x| \leq R \}  $$
$$f_1 (1/ \bar{x}) : = R^{2m} \  \frac{F 
(\bar{x}_0,\bar{x}_1)}{\bar{x}_1^m}\ {\rm on} \  \{ x | |x|
\geq R\}
$$
and this trick allows to carry out local computations comfortably.

Let us go now to the main definition:

\begin{df}
An $(a,b,c) $ surface is the minimal resolution of singularities of a 
simple bidouble
cover
$S$ of
$(\PP^1 \times \PP^1)$ of type $((2a, 2b),(2c,2b)$ having at most 
Rational Double Points
as singularities.

An  $(a,b,c)^{nd} $ surface is defined more generally as
(the minimal resolution of singularities of) a {\bf natural deformation} of
an  $(a,b,c)$ surface with R.D.P.'s : i.e., the canonical model of an 
$(a,b,c)^{nd} $
surface is embedded in the total space
of the direct sum of 2 line bundles  $L_1, L_2$ (whose corresponding sheaves of
sections are
${\hol}_{\PP^1\times \PP^1}(a,b), {\hol}_{\PP^1\times \PP^1}(c,b)$), 
and defined
there by a pair of equations

$$  (***) \ \ \  \ z_{a,b}^2 =  f_{2a,2b} (x,y)  + w_{c,b} \phi_{2a-c,b}(x,y)$$
  $$ w_{c,b}^2 = g_{2c,2b}(x,y) +  z_{a,b}  \psi_{2c-a,b}(x,y)$$

where f,g ,$\phi, \psi$, are bihomogeneous polynomials ,  belonging to
   respective vector spaces of sections of line bundles:
   $ f \in H^0({\PP^1\times \PP^1}, {\hol}_{\PP^1\times \PP^1}(2a,2b)) ,
\phi \in H^0({\PP^1\times \PP^1}, {\hol}_{\PP^1\times \PP^1}(2a-c,b))$ and
$ g \in H^0({\PP^1\times \PP^1}, {\hol}_{\PP^1\times \PP^1}(2c,2d)),
\psi \in H^0({\PP^1\times \PP^1}, {\hol}_{\PP^1\times \PP^1}(2c-a,b))$.

A {\bf perturbation} of an  $(a,b,c) $ surface is an oriented  smooth 
4-manifold defined
by equations as $(***)$, but where the sections $\phi, \psi$ are
differentiable, and we have a {\bf dianalytic perturbation} if
  $\phi, \psi$ are polynomials in the variables $x_i, y_j, \overline{x_i},
  \overline{y_j}$, according to the respective positivity or negativity of
the entries of the bidegree.

%We shall moreover assume, for technical reasons , that $ a \geq 2c + 1$,
%   $ b \geq c+2$ , and that a,b,c are even and $\geq 3$.

\end{df}

\begin{oss}
By the previous formulae ,

1)$(a,b,c)$ surfaces have the same invariants
$\chi (S)= 2 (a+c-2) (b-1) +  b (a+c) ,  K^2_{S} = 16 (a+c-2) (b-1) $ .

2) the divisibility of their canonical class is $G.C.D. ((a+c-2), 2 (b-1))  $.

3) Moreover, we saw that $(a,b,c)$ surfaces are simply connected, thus

4) Once we fix $b$ and the sum $(a+c) = s $, the corresponding 
$(a,b,c)$ surfaces
are all homeomorphic.
\end{oss}

As a matter of fact, once we fix $b$ and the sum $(a+c)$,  the surfaces in the
respective families
   are homeomorphic by a homeomorphism carrying the canonical class
to the canonical class. This fact is a consequence
   of the following  proposition, which we learnt from \cite{man2}

\begin{prop}
Let $S$, $S'$ be simply connected minimal surfaces of general type
such that $\chi(S) = \chi(S') \geq 2$ , $K^2_S = K^2_{S'}$, and moreover such
   that the divisibility indices of $K_S $ and $K_{S'}$ are the same.

Then there exists a homeomorphism $F$ between $S $ and $S'$,
unique up to isotopy,
carrying $K_{S'}$ to $K_S $.
\end{prop}

\begin{proof}
By Freedman's theorem (\cite{free}, cf. especially
\cite{f-q}, page 162) for each isometry
$h: H_2(S, \Z) \rightarrow  H_2(S', \Z)$ there exists a homeomorphism
$F$ between $S $ and $S'$,
unique up to isotopy, such that $F_{*} = h$. In fact, $S$ and $S'$ are
smooth $4$-manifolds, whence
   the Kirby-Siebenmann invariant vanishes.

Our hypotheses that $\chi(S) = \chi(S')$ , $K^2_S = K^2_{S'}$
and that $K_S, K_{S'}$ have the same divisibility
imply that the two lattices $ H_2(S, \Z)$, $ H_2(S', \Z)$ have the
same rank, signature and parity, whence they are isometric
since $S, S'$ are algebraic surfaces.
  Finally, by Wall's theorem (\cite{wall}) (cf. also
\cite{man2}, page 93) such isometry $h$ exists since the vectors
corresponding to the respective canonical classes have the same
divisibility and by Wu's theorem they are characteristic:
   in fact Wall's condition
$ b_2 - |\sigma| \geq 4$
($\sigma$ being the signature of the intersection form)
is equivalent to $\chi \geq 2$.
\end{proof}

\medskip

We come now to the main result of this section (see \cite{c-w} for details)

\begin{teo}\label{diffabc}
Let $S$ be an $(a,b,c)$ - surface and $S'$ be an
$(a+1,b,c-1)$-surface. Moreover, assume that $a,b,c-1 \geq 2$. Then
$S$ and $S'$ are diffeomorphic.
\end{teo}

{\em Idea of the Proof.}

Before we dwell into the proof, let us explain the geometric argument
which led me to  conjecture the above theorem in  1997.

Assume that the polynomials $f, g$ define curves $C, D$ which are 
union of vertical and
horizontal lines.
Fix for simplicity affine coordinates in $\PP^1$. Then we may assume, 
without loss of
generality, that the curve $C$ is constituted by the horizontal
lines  $ y = 1,  \dots y =  2b $, and by the vertical lines
$ x= 2, \dots x = 2a + 1$, while the curve $D$ is formed by the horizontal
lines  $ y = - 1,  \dots y =  - 2b $, and by the vertical lines $ 
x=0$, $x= 1/4$,
$ x= 2a + 2, \dots x = 2a + 2c -1$. The corresponding surface $X$ has 
double points as
singularities, and its minimal resolution is a deformation of a smooth
$(a,b,c)$ - surface (by the cited results of Brieskorn and Tjurina).

Likewise, we let $X'$ be the singular surface corresponding to
the curve $C'$  constituted by the horizontal
lines  $ y = 1,  \dots y =  2b $, and by the vertical lines $ x=0$, $x= 1/4$,
$ x= 2, \dots x = 2a + 1$, and to the curve $D'$  formed by the horizontal
lines  $ y = - 1,  \dots y =  - 2b $, and by the ($2c-2$) vertical lines
$ x= 2a + 2, \dots x = 2a + 2c -1$.

We can split $X$ as the union $ X_0 \cup X_{\infty}$, where
$ X_0 := \{(x,y,z,w) | \ | x| \leq 1 \}$,
$ X_{\infty} := \{(x,y,z,w) | \ | x| \geq 1 \}$, and
similarly
  $ X' = X'_0 \cup X'_{\infty}$.

By our construction, we see immediately that $X'_{\infty} = X_{\infty}$, while
there is a natural diffeomorphism $\Phi$ of $ X_0 \cong X'_0$.

It suffices in fact to set $\Phi (x,y,z,w) = (x, -y, w, z)$.

The conclusion is that both $S$ and $S'$ are obtained glueing the same
two 4-manifolds with boundary $S_0 ,S_{\infty}$ glueing the boundary
$\partial X_0 =\partial  X_{\infty}$ once through the identity,
and another time through the diffeomorphism $\Phi$. It will follow
that the two 4-manifolds are diffeomorphic if the  diffeomorphism 
$\Phi |_{\partial
S_0  }$ admits an extension to a diffeomorphism of $S_0$.

1) The relation with Lefschetz fibrations comes from the form of
$\Phi$, since $\Phi$ does not affect the variable $x$, but it is
essentially given by a diffeomorphism $ \Psi$ of the fibre over $ x=1$,
$$\Psi (y,z,w) = ( -y, w, z).$$

Now, the projection  of an $(a,b,c)$ surface onto $\PP^1$ via the
coordinate $x$ is not a Lefschetz fibration, even if $f,g$ are general, since
each time one of the two curves $C,D$ has a vertical tangent , we 
shall have two nodes
on the corresponding fibre.
But a smooth general natural deformation

\begin{eqnarray}  z^2 &=&
   f(x,y) + w \phi(x,y)\\
   w^2 &=&
   g(x,y) + z \psi(x,y) \nonumber ,\end{eqnarray}

would do the game if  $\phi \neq 0$ (i.e., $2a-c > 0$) and
  $\psi \neq 0$ (i.e., $2c-a > 0$).

Otherwise, it is enough to take a perturbation as in the previous
definition (a dianalytic one suffices),
and we can realize both surfaces $S$ and $S'$  as symplectic Lefschetz
    fibrations (cf. also  \cite{donsympl}, \cite{g-s}).

2) The above argument about $S, S'$ being the glueing of the same
two manifolds with boundary $S_0 , S_{\infty}$ translates directly
into the property that the corresponding Lefschetz fibrations over $\PP^1$
are fibre
    sums of the same pair of Lefschetz fibrations over the
respective complex discs $\{ x | \ | x| \leq 1 \}$,$\{ x| \ | x| \geq 1 \}$.

3) Once the first fibre sum is presented as composition of two
    factorizations and the second as twisted by the  'rotation' $\Psi$,
(i.e., as we saw, the same composition of
    factorizations, where the second is conjugated by  $\Psi$),
in order to prove that
    the two fibre sums are equivalent, it suffices to apply a
very simple lemma, which can be found in
  \cite{auroux02}, and that we reproduce here because of its beauty

\begin{lem}\label{auroux} {\bf ( Auroux)} Let $\tau$ be a Dehn twist and let
$F $  be a factorization of a central element
$ \phi \in \MM ap_g$,
$ \tau_1 \circ \tau_2 \circ \dots \circ \tau_m = \phi$.

If there is a factorization $F'$ such that $F$ is Hurwitz equivalent
to $\tau \circ F'$, then $ (F)_{\tau}$ is Hurwitz equivalent
to $F$.

In particular, if $F$ is a factorization of the identity,
$\Psi = \Pi _h \tau'_h$, and $\forall h \ \exists F'_h$ such that
$F \cong \tau'_h \circ F'_h$, then the fibre sum with the Lefschetz
pencil associated with $F$ yields the same Lefschetz pencil
as the fibre sum twisted by $\Psi$.

\end{lem}
\Proof

If $\cong$ denotes Hurwitz equivalence, then
$$ (F)_{\tau} \cong \tau \circ (F')_{\tau} \cong F' \circ \tau
\cong (\tau)_{(F')^{-1}} \circ F' =  \tau \circ F' \cong F.$$

\QED

\begin{cor} \label{cor}
Notation as above, assume that $F$:
$ \tau_1 \circ \tau_2 \circ \dots \circ \tau_m = \phi$
is a factorization of the Identity
and that $\Psi$ is a product of some Dehn twists $\tau_i$ appearing
in $F$. Then the fibre sum with the Lefschetz
pencil associated with $F$ yields the same result
as the same fibre sum twisted by $\Psi$.
\end{cor}

\Proof
We need only to verify that for each $h$ , there is $F'_h$ such that
$F \cong \tau_h \circ F'_h$.

But this is immediately obtained by applying $h-1$ Hurwitz moves,
the first one between $\tau_{h-1}$ and $\tau_h$, and proceeding further
to the left till we obtain $\tau_h$ as first factor.

\qed

4)  It suffices now to show that the
diffeomorphism $\Psi$ is in the subgroup of the mapping
class group generated by the Dehn twists which appear in the first
factorization.

\bigskip
Figure \ref{figura6} below shows the fibre $C$ of the fibration in the case
   $2b = 6$: it is a bidouble cover of $\PP^1$, which we can assume to
   be given by the equations $z^2 = F(y)$, $w^2 = F(-y)$, where the
   roots of $F$ are the integers $1, \ldots, 2b$.

Moreover, one sees that the monodromy of the fibration at the boundary
of the disc is trivial, and we saw that the map $\Psi$ is the diffeomorphism of
order $2$ given by $y \mapsto -y$, $z \mapsto w$, $w \mapsto z$, which
in our figure is given as a rotation of $180$ degrees around an axis
inclined in direction north-west.

The figure shows a dihedral symmetry, where the automorphism of order
$4$ is given by $y \mapsto -y$, $z \mapsto -w$, $w \mapsto z$.

\begin{figure}[htbp]
\begin{center}
\scalebox{1}{\includegraphics{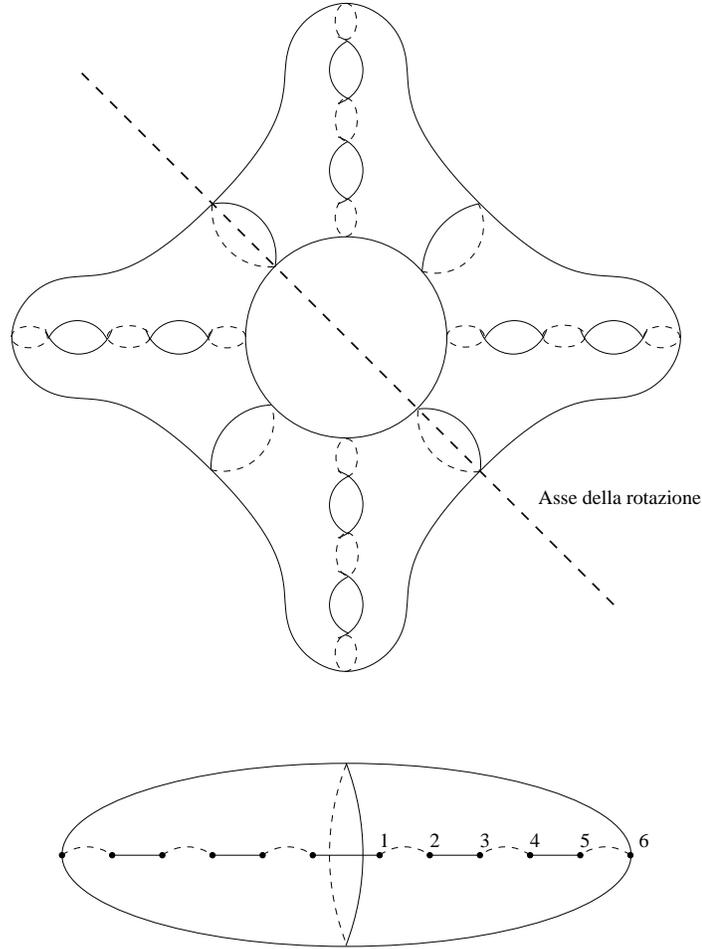}}
\end{center}
\caption{The curve $C$ with a dihedral symmetry}
\label{figura6}
\end{figure}

5) A first   part of the proof, which we skip here, consists in identifying
the Dehn twists which appear in the first factorization.

It turns out that, among the Dehn twists which appear in the first 
factorization,
there are those which correspond to the inverse images of the segments
between two consecutive integers (cf. figure \ref{figura6}). These 
circles can be
organized on the curve $C$ in six chains (not disjoint) and finally one
  reduces oneself to the computational heart of the proof: showing 
that the isotopy class
of $\Psi$ is the same as the product $\Psi '$ of the six Coxeter 
elements associated to
such chains.

We recall here that, given a chain of curves
$\alpha_1, \dots \alpha_n$ on a Riemann surface,  the {\em Coxeter 
element associated
to the chain} is the product
$$
\Delta:  = (T_{\alpha_1})(T_{\alpha_2} T_{\alpha_1}) \ldots
(T_{\alpha_n} T_{\alpha_{n-1}} \ldots T_{\alpha_1})
$$
of the Dehn twists associated to the curves of the chain.

\bigskip

In order to finally prove that $\Psi '$ (the product of such
  Coxeter elements) and $\Psi$ are isotopic, one observes that if one
removes the above cited chains of circles from the curve $C$, one
obtains $4$ connected components which are diffeomorphic to
circles. By a result of Epstein it is then sufficient to verify that
$\Psi$ and $\Psi '$ send each such curve to a pair of isotopic curves:
this last step needs a list of lengthy (though easy) verifications,
for which it is necessary to have explicit drawings.

For details we refer to the original paper \cite{c-w}.

\bigskip

\newpage

\section{Epilogue: Deformation,  diffeomorphism  and symplectomorphism type
of surfaces of general type.}

As we repeatedly said, one of the fundamental problems in the theory of complex
algebraic surfaces is to understand the moduli spaces of surfaces of general
type, and in particular their connected components, which, as we
saw in the third lecture, parametrize
the deformation equivalence classes of minimal surfaces of general
type, or equivalently of their canonical models.

We remarked that deformation equivalence of two minimal models
$S, S'$ implies their canonical
symplectomorphism and a fortiori an oriented diffeomorphism
 preserving the canonical class (a fortiori, a homeomorphism
with such a property).

In the late eighties Friedman and Morgan (cf. \cite{f-m}) made the 
bold conjecture
that two algebraic surfaces are diffeomorphic if and only if they are
deformation equivalent. We will abbreviate this conjecture
by the acronym {\bf def = diff}.
Indeed, I should
point out that I had made the opposite conjecture in the early eighties
(cf. \cite{katata}).

  Later in this section we shall briefly describe the first
counterexamples, due to  M. Manetti (cf. \cite{man4}):
these have the small disadvantage of providing nonsimplyconnected
surfaces, but the great advantage of yielding non deformation 
equivalent surfaces
which are canonically symplectomorphic (see  \cite{cat02} and 
\cite{cat06} for more
details).

We already described in Lecture 4 some easy counterexamples  
to this conjecture (cf. \cite{cat4},
  \cite{k-k},\cite{bcg}),  given by pairs of complex conjugate
surfaces, which are not deformation equivalent to their 
complex conjugate
surface.

We might say that, although describing some interesting phenomena, 
the counterexamples
contained in the cited papers by
Catanese, Kharlamov-Kulikov, Bauer-Catanese-Grunewald
  are 'cheap', since the diffeomorphism carries the canonical class to 
its opposite.
I was  recently informed (\cite{friedman}) by R. Friedman that also he 
and Morgan were
aware of such  'complex conjugate' counterexamples, but for
the case of some elliptic surfaces having an infinite 
fundamental group.

After the examples by Manetti it was however still possible to weaken 
the conjecture
{\bf def = diff} in the following way.

\begin{question} Is the speculation {\bf def = diff} true if one
requires the diffeomorphism $\phi : S \rightarrow S'$ to send the
first Chern class $c_1(K_S) \in H^2 (S, \mathbb{Z})$ in  $c_1 (K_{S'})$
and moreover one requires the surfaces to be simply connected?
\end{question}

  But even this weaker question turned
out to have a negative answer,
as it was shown in our joint work with Wajnryb (\cite{c-w}).

\begin{teo} (\cite{c-w})\label{cw}

For each natural number $h$ there are simply connected surfaces $S_1,
\ldots , S_h$ which are pairwise diffeomorphic, but not deformation equivalent.
\end{teo}

The following remark shows that the statement of the theorem implies
a negative answer to the above question.

\begin{oss}
If two surfaces are deformation equivalent, then there exists a diffeomorphism
sending the canonical class $c_1(K_S) \in H^2(S, \mathbb{Z})$ to the
canonical class $c_1(K_{S'})$. On the other hand, by the cited result of
Seiberg - Witten theory we know that a diffeomorphism sends the
canonical class of a minimal surface $S$ to $\pm
c_1(K_{S'})$. Therefore, if one gives at least three surfaces, which
are pairwise diffeomorphic, one finds at least two surfaces with the
property that there exists a diffeomorphism between them sending the
canonical class of one to the canonical class of the other.
\end{oss}

\subsection{Deformations in the large of ABC surfaces. }

The above surfaces $S_1, \ldots , S_h$ in theorem \ref{cw} belong to
the class of the so-called $(a,b,c)$-
surfaces, whose diffeomorphism type was shown in the previous Lecture
to depend only upon the integers $(a+c)$ and $b$.

The above theorem \ref{cw} is thus implied by the  following result:

\begin{teo} \label{nondef}Let  $S$, $S'$ be simple bidouble
covers of  ${\Bbb P}^1
\times  {\Bbb P}^1  $ of respective
   types ((2a, 2b),(2c,2b), and (2a + 2k, 2b),(2c - 2k,2b) , and assume

\begin{itemize}
\item
   (I) $a,b,c, k$
are strictly  positive even integers with $ a, b, c-k \geq 4$
\item
(II)  $ a \geq 2c + 1$,
\item
(III) $ b \geq c+2$  and either
\item
(IV1) $ b \geq 2a + 2k -1$ or\\
(IV2) $a \geq b + 2$
\end{itemize}
Then  $S$ and  $S'$ are not deformation equivalent.
\end{teo}

  The  theorem uses techniques which have been developed
in a series of papers by the  author and by Manetti
(\cite{cat1}, \cite{cat2}, \cite{cat3}, \cite{man1},
  \cite{man3}). They use essentially the local deformation
theory a' la Kuranishi  for the canonical models, normal 
degenerations of smooth
surfaces and a study of
quotient singularities of rational double points and of their smoothings
(this method was used in \cite{cat2} in order to study the closure
in the moduli space of a class of bidouble covers of $\PP^1 \times \PP^1$
satisfying other types of inequalities).

Although the proof can be found in \cite{cat02}, \cite{c-w}, and in the
Lecture Notes by Manetti in this volume, I  believe it worthwhile
to sketch the main ideas and arguments of the proof.

{\em Main arguments  of the Proof.}

These are the three main steps of the proof.

Step I : determination of a subset $\frak N_{a,b,c}$ of the moduli space.

Step II: proof that $\frak N_{a,b,c}$ is an open set.

Step III: proof that $\frak N_{a,b,c}$ is a closed set.

Let us first of all explain the relevance of hypothesis 2) for step III.
If we consider the natural deformations of $(a,b,c)$ surfaces, which are
parametrized by a quadruple of polynomials $(f, g, \phi, \psi)$ and given
by the two equations
$$
z^2 = f(x,y) + w \phi (x,y),
$$
$$
w^2 = g(x,y) + z \psi (x,y),
$$
we observe that $f$ and $g$ are polynomials of respective bidegrees $(2a,2b)$,
$(2c,2b)$, while $\phi$ and $\psi$ have respective bidegrees
$(2a-c,b)$, $(2c-a,b)$. Hence
  $a \geq 2c+1$, implies that $\psi
\equiv 0$, therefore every small deformation preserves the structure
of an iterated double cover. This means that the quotient $Y$ of
our canonical model $X$ by the involution $ z \mapsto - z$
admits an involution $ w \mapsto -w$, whose quotient is
indeed  $\PP^1 \times \PP^1$.

This fact will play a special role in the study of limits of such
$(a,b,c)^{nd}$ surfaces, showing that this iterated double cover
structure  passes in a
suitable way to the limit, hence $\frak N_{a,b,c}$
is a closed subset of the moduli space.

{\bf Step I.}

The family $(\frak N_{a,b,c})$ consists of all the (minimal resolutions of the)
  natural deformations of simple bidouble covers  of the
Segre-Hirzebruch
   surfaces $\F_{2h}  $ which have only
   Rational Double Points as singularities and are of type ((2a, 2b),(2c,2b).

In order to explain what this means, let us
recall, as in \cite{cat0} pages 105-111,
that a basis of the Picard group of $\F_{2h}  $ is provided, for $ h \geq 1$,
by the fibre $F$ of the projection to $\PP^1$, and by $F' : =
\sigma_{\infty} + h F$,
where $\sigma_{\infty}$ is the unique section with negative
self-intersection $ = - 2h$. Observe that $F^2 = {F'}^2 = 0, F F' =
1$, and that
$F$ is nef, while $F' \cdot \sigma_{\infty} = - h$.

We set $\sigma_0 :=  \sigma_{\infty} + 2 h F$, so that  $\sigma_{\infty}
   \sigma_{0} = 0$, and we observe (cf. Lemma 2.7 of \cite{cat0})
that $ |m \sigma_0 + n F |  $  has no base point if and only if $m,n \geq 0$.
Moreover, $ |m \sigma_0 + n F |  $ contains  $\sigma_{\infty}$ with
multiplicity
$ \geq 2$ if $ n < - 2h$.

At this moment, the above remarks and the inequalities (II), (III), (IV) can be
used to imply that all natural deformations have the structure of an 
iterated double
covering, since their canonical models are defined by the following 
two equations:

$$
z^2 = f(x,y) + w \phi (x,y),
$$
$$
w^2 = g(x,y) .
$$

{\bf Step II.}

A key point here is to look only
at the deformation theory of the canonical models.

To prove that the family of canonical models $ (\frak N_{a,b,c})$ yields
an open set in the moduli space it suffices to show that,
   for each surface $X$, the Kodaira Spencer map is surjective.

In fact, one can  see  as in in \cite{cat0} that the family  $ (\frak 
N_{a,b,c})$ is
parametrized
   by a smooth variety which surjects onto $H^1 ( \Theta_{\F} )$.

Observe  that the tangent space to the Deformations of $X$ is provided by
$ {\rm Ext}^1_{\hol_X} (\Omega^1_X, \hol_X) $ .

Denoting by $\pi : X \ra  \F :=
\F_{2h} $ the projection map and differentiating equations
(7) we get an exact sequence for $\Omega^1_X$
$$ o \ra \pi^* (\Omega^1_{\F}) \ra \Omega^1_X \ra \hol_{R_z} (- R_z) \oplus
\hol_{R_w} (- R_w)  \ra 0 $$ as in  (1.7) of \cite{man1}, where $
{R_z} = div (z),
{R_w} = div(w)$.

Applying the derived exact sequence for
$ {\rm Hom}_{\hol_X} ( \dots, \hol_X) $ we obtain the same exact sequence as
Theorem (2.7) of \cite{cat0}, and (1.9) of \cite{man1},
namely:

$$(**) \ 0 \ra H^0 (\Theta_{X}) \ra H^0 (\pi^* \Theta_{\F})  \ra
    H^0 (\hol_{R_z} (2 R_z)) \oplus
H^0 (\hol_{R_w} (2 R_w)) \ra $$
$$ \ra{\rm Ext}^1_{\hol_X} (\Omega^1_X, \hol_X) \ra
H^1 (\pi^* \Theta_{\F}) $$ .

There is now some technical argument, quite similar to
  the one given in
\cite{cat0}, and  where our inequalities are used
in order to show that $H^1 (\pi^* \Theta_{\F}) =
H^1 ( \Theta_{\F} \otimes \pi_* (\hol_X)) $ equals $
H^1 ( \Theta_{\F} )$:
  we refer to \cite{c-w}  for details.

Summarizing the proof of step II, we  observe
that the smooth parameter space of our family surjects onto $H^1 (
\Theta_{\F})$,
and its kernel, provided by the natural deformations
with fixed base $\F_{2h}$, surjects onto
$ H^0 (\hol_{R_z} (2 R_z)) \oplus
H^0 (\hol_{R_w} (2 R_w))  $. Thus the Kodaira Spencer is onto and we get an
open set in the moduli space.

{\bf Step III.}

We want now to show that our family  $\frak N_{a,b,c}$ yields a 
closed set in the moduli
space.

It is clear at this moment that we obtained an irreducible component
of the moduli space. Let us consider the surface over the generic point of
the base space of our family:
then it has  $\Z / 2$ in the automorphism group
   ( sending $z \ra -z$, as already mentioned).

As shown in \cite{cat0}, this automorphism acts then biregularly on the
canonical model $X_0$ of each surface corresponding  to a point in
the closure of our open set.
This holds in fact more generally for the action of
any finite group $G$: the representation of $G$ on $H^0 (S, \hol (5 K_S))$
depends on discrete data, whence it is fixed in a family, and then the
set of fixed points in the pseudomoduli space $ \{ X | g (X) = X \ 
\forall g \in G\}$
is a closed set.

We use now the methods of \cite{cat2} and \cite{man3},
and more specifically we can apply Theorem 4.1 of  \cite{man3} to conclude
with

{\bf Claim III .1}
If $X_0$ is a  canonical model which is a limit of canonical models 
$X_t$ of surfaces
$S_t$ in our family, then the quotient $Y_0$ of $X_0$ by the subgroup 
$\Z / 2 \subset
Aut(X_0)$ mentioned above is a surface with Rational Double Points.

{\bf Claim III .2}
The family of such quotients $ Y_t$ has a  $\Z / 2$-action over
the generic point, and dividing by it we get ( cf. \cite [Theorem 4.10] {man3})
as quotient $Z_0$ a Hirzebruch surface. Thus our surface $X_0$
is also an iterated double cover of some
$\F_{2h}  $, hence it belongs to the family we constructed.

{\em Argument for claim III.1}
Since smooth canonical models are dense, we may assume that $X_0$ is a limit
of a 1-parameter family $X_t$ of smooth canonical models; for the same reason
we may assume that the quotient $Y_0$ is the limit of smooth surfaces
  $ Y_t = X_t / (\Z / 2)$ (of general
type if $ c, b \geq 3$).

Whence,

1) $Y_0$ has singularities which are quotient of Rational Double Points by
$(\Z / 2)$

2) $Y_t$ is a smoothing of $Y_0$, and since we assume the integers 
$c,b$ to be even,
the canonical divisor of $Y_t$ is 2-divisible.

Now, using Theorem \ref{autrdp}, the involutions acting on RDP's can 
be classified
( cf. \cite{cat2} for this and the following), and it turns out that the
quotient singularities are again RDP's, with 2 possible exceptions:

Type (c): the singularity of $Y_0$ is a quotient singularity of type
$ \frac{1}{4k+2} (1, 2k) $, and $X_0$ is the $A_{2k}$ singularity,
quotient by the subgroup $ 2 \Z / (4k+2) \Z$.

Type (e): the singularity of $Y_0$ is a quotient singularity of type
$ \frac{1}{4k+4} (1, 2k + 1) $, and $X_0$ is the $A_{2k+1}$ singularity,
quotient by the subgroup $ 2 \Z / (4k+4) \Z$.

The versal families of deformations of the above singularities have been
described by Riemenschneider in \cite{riem}, who showed:

(C) In the case of type (c), the base space is smooth, and it yields 
a smoothing
admitting a simultaneous resolution.

(E) In the case of type (e), the base space consists of two smooth components
intersecting transversally, $ T_1 \cup T_2$. $T_1$ yields a smoothing
admitting a simultaneous resolution (we denote this case by ` case (E1)').

Hypothesis 2), of 2-divisibility of the canonical divisor of $Y_t$, is used
in two ways. The first consequence is that the intersection form on
$H^2 (Y_t, \Z)$ is even; since however the Milnor fibre of the smoothing
is contained in $Y_t$, it follows that no 2-cycle in the Milnor fibre
can have odd selfintersection number. This then excludes case (C), and also
case (E1) for $ k \geq 1$.

In case (E2) we have a socalled $\Z$-Gorenstein smoothing, namely,
the $T_2$ family is the quotient of the hypersurface
$$(***) \ \  uv - z^{2n} = \Sigma_{h=0}^{1} t_h z ^{hn}$$
by the involution sending $(u,v,z) \mapsto (-u,-v,-z)$.

The result is that the Milnor fibre has a double \'etale cover which 
is the Milnor
fibre of $ A_{n-1}$ ( $n = k+1$), in particular its fundamental group
equals $ \Z/ 2$. The universal cover corresponds to the cohomology class
of the canonical divisor. This however contradicts condition 2), and
case E2) is excluded too.

For case E1) $ k=0$ we argue similarly: the involution acts trivially
on the parameter $t$, and in the central fibre it has an isolated fixed point.
Because of simultaneous resolution, the total space $ \cup_t X_t$ may 
be taken to
be smooth, and then the set of fixed points for the involution is a curve
mapping isomorphically on the parameter space $ \{ t\}$. Then the Milnor fibre
should have a  double cover ramified exactly in one point, but this is absurd
since by van Kampen's theorem the point complement is simply connected.

\qed

{\em Argument for claim III.2}

Here, $Z_t : = (Y_t / \Z/2) \cong \PP^1 \times \PP^1 = \F_0$ and 
again the canonical
divisor is 2-divisible. Whence, the same argument as before applies, 
showing that $Z_0$
has necessarily Rational Double Points as singularities.
But again, since the Milnor fibre embeds in $\PP^1 \times \PP^1 = \F_0$,
the intersection form must have negativity at most 1, and be even.
This leaves only the possibility of an $A_1$ singularity. This case 
can be again
excluded by the same argument given for the case  E1) $ k=0$ above.

\qed

\QED

{\em Proof that Theorem \ref{nondef} implies Theorem \ref{cw}.}

It suffices to show what we took up to now for granted:
the irreducible component $\mathfrak N _{a,b,c}$
  uniquely determines the
numbers $a,b,c$ up to the obvious permutations: $a \leftrightarrow c$,
   and , if $ a=c$, the possibilities of exchanging $a$ with $b$.

It was shown more generally in \cite{cat1} theorem 3.8  that the natural
deformations of bidouble covers of type $(2a,2b)(2c,2d)$ yield an irreducible
component of the moduli space, and that these are distinct modulo the obvious
permutations (exchange  type $(2a,2b)(2c,2d)$ with type  $(2c,2d) (2a,2b)$
and with type $(2b,2a)(2d,2c)$). This follows from geometrical properties
of the canonical map at the generic point.

However, the easiest way to see that the irreducible component $\mathfrak N _{a,b,c}$
   determines the numbers $a,b,c$, under the given inequalities II0, 
III), IV) is to
observe that the dimension of $\mathfrak N _{a,b,c}$ equals
$  M:= (b+1) (4 a + c + 3) + 2 b (a+c+1) - 8$.
Recall in fact that $K^2/16 = (a+c-2)(b-1)$, and
$  (8 \chi - K^2 )/8 = b (a+c)$: setting  $\alpha = a+c, \beta = 2b$, 
we get that
$\alpha , \beta $ are then the roots of a quadratic equation,
so they are determined up to exchange, and uniquely if we  restrict our
numbers either to the inequality $ a \geq 2b$ or to the inequality
$b \geq a$..

   Finally  $ M =( \frac\beta2 + 1) ( \alpha + 3) +  \beta (\alpha + 1)
- 8  + 3 a ( \frac\beta2 + 1)$
  then determines  $a$, whence the ordered triple $(a,b,c)$.

\qed

\begin{oss}
If, as in \cite{cat02}, we assume

(IV2) $ a \geq b+2$,

then the connected component $\frak N_{a,b,c}$ of the moduli space
contains only iterated double covers of ${\PP}^1 \times  {\PP}^1  $.
\end{oss}

\bigskip

\subsection{Manetti  surfaces}

Manetti in \cite{man4} considers  surfaces which are desingularization of
certain $(\Z/2)^r$ covers $X$ of  rational surfaces $Y$ which are 
blowup of   the
quadric $ Q := \PP^1 \times \PP^1 $ at $n$ points $P_1, \dots P_n$.

His construction is made rather complicated, not only by the desire 
to construct an
arbitrarily high number of surfaces which are pairwise diffeomorphic but not
deformation equivalent, but also by the  crucial target  to obtain that 
every small
deformation is again such a Galois   $(\Z/2)^r$ cover. This 
requirement makes the
construction not very explicit (Lemma 3.6 ibidem).

Let us briefly recall the structure of normal finite $(\Z/2)^r$ 
covers with smooth base
$Y$ (compare \cite{Pardini},  \cite{man4}, and also \cite{bc} for a 
description in terms of
the monodromy homomorphism).

We denote by $G = (\Z/2)^r$ the Galois group, and by $\sigma$ an 
element of $G$.
We denote by $G^{\vee} := Hom (G, \C^* )$ the dual group of 
characters, $G^{\vee} \cong
(\Z/2)^r$, and by $\chi$ an element of $G^{\vee}$. As for any flat 
finite abelian
covering
$ f \colon X \to Y$ we have

$$f_* {\hol}_X =\bigoplus_{\chi \in G^{\vee} } {\hol}_Y
(- L_{\chi}) =  {\hol}_Y \oplus  ( \bigoplus_{\chi \in G^{\vee} 
\setminus \{0 \}}
{\hol}_Y (- L_{\chi})).$$

To each element of the Galois group $ \sigma \in G$ one associates a divisor
$D_{\sigma}$,such that $ 2 D_{\sigma}$ is the direct image 
divisor $f_* (R_{\sigma})$,
$R_{\sigma}$ being the divisorial part of the set of fixed points for $\sigma$.

Let $x_{\sigma}$ be a section such that $ div (x_{\sigma}) = D_{\sigma}$: then
the algebra structure on $f_* {\hol}_X$ is given by the following
symmetric bilinear multiplication maps:

$$ {\hol}_Y (- L_{\chi}) \otimes {\hol}_Y (- L_{\eta}) \to {\hol}_Y 
(- L_{\chi +
\eta})$$
associated to the section $$ x_{\chi, \eta} \in  H^0 (Y,{\hol}_Y
( L_{\eta}  +  L_{\chi} -  L_{\chi +
\eta}) ), \  x_{\chi, \eta} : = \prod _{\chi (\sigma) = \eta (\sigma) 
= 1} x_{\sigma}
.$$
Associativity follows since, given  characters $\chi, \eta, \theta$,
$\{ \sigma|( \chi + \eta ) (\sigma) = \theta (\sigma) = 1 \}$ is the 
disjoint union
of $\{ \sigma| \chi (\sigma) = \theta (\sigma) = 1,  \eta  (\sigma)= 0 \}$
and of $\{ \sigma| \eta (\sigma) = \theta (\sigma) = 1,  \chi  (\sigma)= 0 \}$,
so that $$ {\hol}_Y (- L_{\chi}) \otimes {\hol}_Y (- L_{\eta})
\otimes {\hol}_Y (- L_{\theta}) \to {\hol}_Y (- L_{\chi +
\eta + \theta})$$ is given by the section $\prod _{ \sigma \in \Sigma 
} x_{\sigma} ,$
where $$\Sigma : = \{ \sigma | \chi (\sigma) = \eta (\sigma) = 1, {\rm or} \
\chi (\sigma) = \theta (\sigma) = 1 , {\rm or} \
\eta (\sigma) = \theta (\sigma) = 1\}.$$

In particular, the covering $ f \colon X \to Y$ is embedded in the 
vector bundle
$\VV$ which is the direct sum of the line bundles whose sheaves of 
sections are the
$ {\hol}_Y (- L_{\chi})$, and is there defined by equations
$$ z_{\chi}  z_{\eta} = z_{\chi + \eta} \prod _{\chi (\sigma) = \eta 
(\sigma) = 1}
x_{\sigma}.$$

Noteworthy is the special case $ \chi = \eta$, where $ \chi + \eta$ 
is the trivial
character $1$,  and $ z_1 = 1$.

In particular, let $\chi_1, \dots \chi_r$ be a basis of $G^{\vee} \cong
(\Z/2)^r$, and set $ z_i := z_{\chi_i}$. We get then the $r$ equations
$$ (\sharp ) \ z_i ^2 =  \prod _{\chi_i (\sigma)  = 1}
x_{\sigma}.$$
These equations determine the field extension, hence one gets $X$ as 
the normalization
of the Galois cover given by $(\sharp )$.

We can summarize the above discussion in the following

\begin{prop}\label{data}
A normal finite $G \cong (\Z/2)^r$ covering of smooth variety $Y$ is completely
determined by the datum of

1) reduced effective divisors $  D_{\sigma}$, $ \forall {\sigma} \in G$,
which have no common components

2) divisor classes $ L_1, \dots L_r$, for $\chi_1, \dots \chi_r$  a 
basis of $G^{\vee}
$, such that we  have the following linear equivalence

3) $$  2 L_i \equiv  \sum _{\chi_i (\sigma)  = 1}
D_{\sigma}. $$

Conversely, given the datum of 1) and 2), if 3) holds, we obtain a normal
scheme $X$ with a finite $G \cong (\Z/2)^r$ covering $f \colon X \to Y$.

\end{prop}

{\em Idea of the proof}

It suffices to determine the divisors $L_{\chi}$ for the other 
elements of $G^{\vee} $.
But since any $\chi$ is a sum of basis elements, it suffices to 
exploit the fact that
the linear equivalences
$$  L_{\chi +
\eta} \equiv  L_{\eta}  +  L_{\chi}  \ -  \sum _{\chi (\sigma) = \eta 
(\sigma) = 1}
D_{\sigma} $$
must hold, and apply induction. Since the covering is well defined as
the normalization
of the Galois cover given by $(\sharp )$, each $L_{\chi}$ is well defined.
Then the above formulae determine explicitly the ring structure
of $ f_{*} \hol_X$, hence $X$.

\qed

A natural question is of course when the scheme $X$ is a variety, 
i.e., $X$ being
normal, when $X$ is connected, or equivalently irreducible. The 
obvious answer is that
$X$ is irreducible if and only if the monodromy homomorphism

$$\mu \colon H_1 (Y \setminus (\cup_{\sigma} D_{\sigma}) ,\Z) \to G   $$
is surjective.

\begin{oss}
As a matter of fact, we know, from the cited theorem of Grauert and Remmert,
that $\mu$ determines the covering. It is therefore worthwhile to see how
$\mu$ determines the datum of 1) and 2).

Write for this purpose the branch locus $D : = \sum_{\sigma} D_{\sigma}$
as a sum of irreducible components $ D_i$. To each $D_i$ corresponds a
simple geometric loop $\ga_i$ around $D_i$, and we set
$ \sigma_i : = \mu (\ga_i)$. Then we have that
$D_{\sigma} : =  \sum_{\sigma_i = \sigma } D_i$. For each character $\chi$,
yielding a double covering associated to the composition $ \chi \circ \mu$,
we must find a divisor class $L_{\chi}$ such that $ 2 L_{\chi} \equiv
\sum_{\chi (\sigma ) = 1} D_{\sigma}$.

Consider the exact sequence
$$ H^{2n-2} (Y, \Z) \to H^{2n-2} (D, \Z) = \oplus_i \Z [D_i] \to H_1 
(Y \setminus D, \Z)
\to  H_1 (Y, \Z)  \to 0$$
and the similar one with $\Z$ replaced by $\Z / 2$. Denote by 
$\Delta$ the subgroup
image of $ \oplus_i \Z  / 2 [D_i]$. The restriction of $\mu$ to $\Delta$ is
completely determined by the knowledge of the $\sigma_i$ 's, and we have
$$ 0 \to \Delta \to  H_1 (Y \setminus D, \Z  / 2)
\to  H_1 (Y, \Z  / 2)  \to 0 .$$

Dualizing, we get $$ 0 \to H^1 (Y, \Z  / 2)   \to  H^1 (Y \setminus D, \Z  / 2)
\to  Hom (\Delta, \Z / 2)  \to 0 .$$

The datum of $\mu$,  extending $\mu | _{\Delta}$ is then seen to correspond
to an affine  space over the vector space $H^1 (Y, \Z  / 2)$: and since
$H^1 (Y, \Z  / 2)$ classifies divisor classes of 2-torsion on $Y$, we infer
that the different choices of  $ L_{\chi} $ such that $ 2 L_{\chi} \equiv
\sum_{\chi (\sigma ) = 1} D_{\sigma}$ correspond bijectively to all the possible 
choices for $\mu$.

\end{oss}

\begin{cor}
Same notation as in proposition \ref{data}. Then the scheme $X$ is 
irreducible if
$\{ \sigma | D_{\sigma} > 0  \}$ generates $G$.

\end{cor}

\Proof
We have seen that if $D_{\sigma} \geq D_i \neq 0$, then $ \mu (\ga_i) 
= \sigma$,
whence we infer that $\mu$ is surjective.

\qed

An important role plays again here the concept of {\bf natural deformations.}
This concept was introduced for bidouble covers in \cite{cat1}, 
definition 2.8, and
extended to the case of abelian covers in \cite{Pardini}, definition 5.1.
However, the two definitions do not coincide, because Pardini takes a 
much larger
parameter space. We propose therefore to call Pardini's case the case of
{\bf extended natural deformations.}

\begin{df}
Let $ f \colon X \to Y$ be a  finite $G \cong (\Z/2)^r$ covering with 
$Y$ smooth and
$X$ normal, so that $X$ is embedded in the vector bundle
$\VV$ defined above and is  defined by equations
$$ z_{\chi}  z_{\eta} = z_{\chi + \eta} \prod _{\chi (\sigma) = \eta 
(\sigma) = 1}
x_{\sigma}.$$
Let $\psi_{\sigma, \chi} $ be a section $\psi_{\sigma, \chi} \in H^0
(Y,\hol_Y (D_{\sigma} - L_{\chi} ) $, given $ \forall \sigma \in G, 
\chi \in G^{\vee}.$
To such a collection we associate an {\bf extended natural deformation},
namely, the subscheme of $\VV$   defined by equations
$$ z_{\chi}  z_{\eta} = z_{\chi + \eta} \prod _{\chi (\sigma) = \eta 
(\sigma) = 1}
( \sum_{\theta} \psi_{\sigma, \theta } \cdot z_{\theta}).$$

We have instead a (restricted) {\bf natural deformation} if we 
restrict ourselves
to the $\theta$'s such that $\theta (\sigma) = 0$,and we consider 
only an equation of
the form
$$ z_{\chi}  z_{\eta} = z_{\chi + \eta} \prod _{\chi (\sigma) = \eta 
(\sigma) = 1}
( \sum_{\theta (\sigma) = 0} \psi_{\sigma, \theta } \cdot z_{\theta}).$$

\end{df}

The deformation results which we explained in the last lecture for simple
bidouble covers work out also for $G \cong (\Z/2)^r$ which are {\bf 
locally simple},
i.e., enjoy the property that for each point $y \in Y$ the $\sigma$'s
such that $y \in D_{\sigma}$ are a linear independent set. This is a 
good notion
since (compare \cite{cat1}, proposition 1.1) if also $X$ is smooth 
the covering is
indeed locally simple.

One has the following result (see \cite{man4}, section 3)

\begin{prop}\label{natdef}
Let $ f : X \to Y$ be a locally simple  $G \cong (\Z/2)^r$ covering 
with $Y$ smooth
and $X$ normal. Then we have the exact sequence

$$
    \oplus_{ \chi (\sigma ) = 0}
( H^0 (\hol_{D_{\sigma}} (D_{\sigma} - L_{\chi}))) \ra
  {\rm Ext}^1_{\hol_X} (\Omega^1_X, \hol_X) \ra
  {\rm Ext}^1_{\hol_X} ( f^* \Omega^1_Y, \hol_X)  .$$
In particular, every small deformation of $X$ is a natural deformation if

i) $  H^1 (\hol_Y ( - L_{\chi})) = 0$

ii) $ {\rm Ext}^1_{\hol_X} ( f^* \Omega^1_Y, \hol_X) = 0.$

If moreover

iii) $ H^0 (\hol_Y (D_{\sigma} - L_{\chi})) = 0$ $ \forall \sigma \in 
G,  \chi \in
G^{\vee},$

every small deformation of $X$ is again a  $G \cong (\Z/2)^r$ covering.
\end{prop}

{\em Comment on the proof.}

In the above proposition condition i) ensures that
  $ H^0 (\hol_Y (D_{\sigma} - L_{\chi})) \to  H^0 (\hol_{D_{\sigma}} 
(D_{\sigma} -
L_{\chi}))$ is surjective.

Condition ii) and the above diagram imply then that the natural deformations
are parametrized by a smooth manifold and have surjective Kodaira Spencer map,
whence they induce all the infinitesimal deformations.

\qed

In Manetti's application one needs an extension of the above result. In fact
  ii) does not hold, since the manifold $Y$ is not rigid (one can move 
the points
$P_1, \dots P_n$ which are blown up in the quadric $Q$). But the 
moral is the same, in
the sense that one can show that all the small deformations of $X$ 
are $G$-coverings
of a small deformation of $Y$.

Before we proceed to the description of the Manetti surfaces, we consider
some simpler surfaces, which however clearly illustrate one of the features
of Manetti's construction.

\begin{df}
A singular bidouble  Manetti surface of type $(a,b)$ and triple of order $n$ is
a  singular bidouble cover of $ Q : = \PP^1 \times \PP^1 $
branched on three smooth curves $C_1,C_2, C_3$ belonging to the linear system
of sections of the sheaf $ \hol_Q (a,b)$ and which intersect in
$ n$ points  $p_1, \dots p_n$, with distinct tangents.

A smooth bidouble  Manetti surface of type $(a,b)$ and triple of order $n$ is
the minimal resolution of singularities $S$ of such a surface $X$ as above.
\end{df}

  \begin{oss}
1) With such a branch locus, a  Galois group  of type  $ G =  ( \Z/2) 
^r$ can be only
$ G =  ( \Z/2) ^3$ or $ G =  (\Z/2) ^2$ (we can exclude the 
uninteresting case  $
G =  ( \Z/2) $). The case $ r = 3$ can only occur if the class of the 
three curves $
(a,b)$ is divisible by two since, as we said, the homology group of 
the complement
$ Q \setminus (\cup_i C_i ) $ is
the  cokernel of the map $ H^2 ( Q,
\Z) \to
\oplus_1^3 (\Z C_i$). The case $ r=3$ is however uninteresting, since 
in this case
the elements $\phi (\ga_i)$ are a  basis, thus over each point $p_i$ we have a
  nodal singularity of the covering surface, which obviously makes us 
remain in the same
  moduli space as the one where the three curves have no intersection 
points whatsoever.

2) Assume that $r=2$, and consider the case where the monodromy
  $\mu$ is such that the  $\mu (\ga_i)$'s  are the three nontrivial elements
of the group $ G =  ( \Z/2) ^2$.

Let $p= p_i$ be a point where the three smooth curves $C_1, C_2, C_3$ intersect
with
  distinct tangents: then over the point $p$ there is a singularity $ (X,x)$
of the type considered in example \ref{Riem}, namely, a quotient singularity
which is analytically the cone over a rational curve of degree 4.

If we blow up the point $p$, and get an exceptional divisor $E$, the 
loop $\ga$  around
the exceptional divisor $E$ is homologous to the sum  of the three 
loops $\ga_1,
\ga_2, \ga_3$ around the respective three curves $C_1, C_2, C_3$.
Hence it must hold  $ \mu (\ga) = \sum_i  \mu (\ga_i) = 0$, and the pull back of 
the
covering does not have $E$ in the branch locus. The inverse image $A$ of $E$
is a $( \Z/2) ^2$ covering of $E$ branched in three points, and we 
conclude that
$A$ is a smooth rational curve of self-intersection $-4$.
\end{oss}

One sees
(compare \cite{bidouble})  that

\begin{prop}

Let $X$ be a  singular bidouble Manetti surface of type $(a,b)$ and 
triple of order $n$:
then if $S$ is the minimal resolution of the singularities $x_1, 
\dots x_n$ of $X$,
then $S$ has the following invariants:

$ K^2_S = 18 ab - 24 (a+b) + 32 - n  $

$ \chi (S) = 4 + 3 (ab -a -b). $

Moreover $S$ is simply connected if $(a,b)$ is not divisible by 2.

\end{prop}

{\em Idea of the proof}
For $ n =0$ these are the standard formulae since $2 K_S = f^* (3a -4,3b-4)$,
and $ \chi (\hol_Q (-a,-b)) = 1 + 1/2 (a (b-2) + b (a-2))$.

For $n > 0$, each singular point $x_n$ lowers $ K^2_S$ by 1, but leaves
$\chi (S)$ invariant. In fact again we have  $2 K_X = f^* (3a -4,3b-4)$,
but $2 K_S =  2 K_X - \sum_i A_i$. For $\chi (S)$, one observes that
$x_i$ is a rational singularity, whence $ \chi (\hol_X) = \chi (\hol_S)$.

It was proven in \cite{cat1} that $S$ is simply connected for $n=0$
when $(a,b)$ is not divisible by 2 (in the contrary case the fundamental group
equals $ \Z / 2$.) Let us then assume that $ n \geq 1$.

Consider now a 1-parameter family $C_{3,t}$, $ t \in T$, such that 
for $ t \neq 0$
$C_{3,t}$ intersects $C_1, C_2$ transversally, while $C_{3,0} = C_3$.
We get a corresponding family $  X_t$ of bidouble covers such that $X_t$
is smooth for $ t \neq 0$ and, as we just saw, simply connected.
Then $S$ is obtained from $X_t, t \neq 0$ replacing the Milnor fibres by
  tubular neighbourhoods of the exceptional divisors $A_i$, $ i =1, \dots n$.
Since $A_i$ is smooth rational, these neighbourhoods are simply connected,
and the result follows then easily by the first van Kampen theorem,
which implies that $ \pi_1 (S)$ is a quotient of $\pi_1 (X_t) , t \neq 0$.

\qed

The important fact is that the above smooth bidouble Manetti surfaces 
of type $(a,b)$
and triple of order $n$ are parametrized, for $ b = la, l \geq 2, n = 
l a (2 a -c), 0 <
2c < a $, by a disconnected parameter space (\cite{man4}, corollary 2.12:
observe that we treat here only the case of $ k=3$ curves).

We cannot discuss here the method of proof, which relies on the 
socalled Brill Noether
theory of special divisors: we only mention that Manetti considers 
the two components
arising form the respective cases where $\hol_{C_1} (p_1 + \dots p_n) 
\cong \hol_{C_1}
(a-c,b)$, $\hol_{C_1} (p_1 + \dots p_n) \cong \hol_{C_1}
(a,b - lc)$, and shows that the closures of these loci yield two 
distinct connected
components.

Unfortunately, one sees easily that smooth bidouble Manetti surfaces admit
natural deformations which are not Galois coverings of the blowup $Y$ of $Q$
in the points $p_1, \dots p_n$, hence Manetti is forced to take  more 
complicated
$G \cong (\Z / 2)^r$ coverings (compare  section
6 of \cite{man4},  especially  page  68, but compare also the crucial 
lemma 3.6).

The Galois group is chosen as
$ G =  ( \Z/2) ^r$,where $ r : = 2  + n + 5$ (once more we make the simplifying
choice $ k=1$ in  6.1 and foll. of \cite{man4}).

\begin{df}
1) Let $G_1 : = ( \Z/2) ^2$,$G_2 : = ( \Z/2) ^n$,  $ G' : =  G_1 
\oplus G_2 \oplus   (
\Z/2) ^4$, $ G : = G'   \oplus   ( \Z/2)$.

2) Let $ D : G' \to Pic (Y)$ be the mapping sending

\begin{itemize}
\item
The three nonzero elements of $G_1$ to the classes of the proper transforms
of the curves $C_i$, i.e., of $ \pi^* (C_j) - \sum_i A_i$
\item
the canonical basis of $G_2$ to the classes of the exceptional divisors $A_i$
\item
the first two elements of the canonical basis of $(\Z/2) ^4$ to the 
pull back of
the class of $\hol_Q (1,0)$, the last two to the pull back of
the class of $\hol_Q (0,1)$
\item
the other elements of $G'$ to the zero class.

\end{itemize}
\end{df}

With the above setting one has (lemma 3.6 of \cite{man4})

\begin{prop}
There is an extension of the map $ D : G' \to Pic (Y)$ to $ D : G \to Pic (Y)$,
and a map $ L : G^{\vee} \to Pic (Y) $,  $\chi \mapsto  L_{\chi}$ such that

i) the cover conditions $ 2 L_{\chi} \equiv
\sum_{\chi (\sigma ) = 1} D_{\sigma}$  are satisfied

ii) $ -  D_{\sigma} +  L_{\chi}$ is an ample divisor

iii) $   D_{\sigma}$ is an ample divisor for $ \sigma \in G \setminus G'$

\end{prop}

\begin{df}

Let now $S$ be a $G$-covering of $Y$ associated to the choice of
some effective divisors $D_{\sigma}$ in the given classes.
$S$ is said to be a Manetti surface.

\end{df}

  For simplicity we assume now
that these divisors $D_{\sigma}$ are smooth and intersect 
transversally, so that $S$ is
smooth.

Condition iii) guarantees that $S$ is connected, while condition
ii) and an extension of the argument of proposition \ref{natdef} shows
that all the small deformations are $G$-coverings of such a rational 
surface $Y$,
blowup of $Q$ at $n$ points.

We are going now only to very briefly sketch the rest of the arguments:

{\bf Step A} It is possible to choose one of the  $   D_{\sigma}$'s to be so
positive that the group of automorphisms of a generic such surface $S$ is just
the group $G$.

{\bf Step B} Using the natural action of $G$ on any such surface, and 
using again
arguments similar to the ones described in Step III of the last lecture,
one sees that we get a closed set of the moduli space.

{\bf Step C} The families of surfaces thus described fibre over the 
corresponding
families of smooth bidouble Manetti surfaces: since for the latter one has more
than one connected component, the same holds for the Manetti surfaces.

In the next section we shall show that the Manetti surfaces 
corresponding to a fixed
choice of the extension $D$  are canonically symplectomorphic.

In particular, they are a strong counterexample to the Def=Diff question.

\bigskip

\subsection{Deformation and canonical symplectomorphism}

We start discussing a simpler case:

\begin{teo}
Let $S$ and $S'$ be the respective minimal resolutions of the 
singularities of two
singular bidouble Manetti surfaces $X, X'$ of type $(a,b)$, both 
triple of the same order $n$:
then $S$ and $S'$ are diffeomorphic, and indeed symplectomorphic for
their canonical symplectic structure.
\end{teo}

\Proof
In order to set up our notation, we denote by $C_1, C_2, C_3 $ the 
three smooth branch curves
for $ p : X \to Q$, and denote by $p_1,  .. , p_n$ the points where 
these three curves
intersect (with distinct tangents): similarly the covering  $ p' : X' \to Q$
determines $C'_1, C'_2, C'_3 $  and $p'_1,  .. , p'_n$. Let $Y$ be 
the blow up of the
quadric $ Q$ at the n points $p_1,  .. , p_n$, so that $S$ is a 
smooth bidouble cover
of $Y$, similarly $ S'$ of $Y'$.

Without loss of generality we may assume that $C_1, C_2$ intersect 
transversally in
$ 2 ab$ points, and similarly $C'_1, C'_2$.

We want to apply theorem \ref{families} to $S, S'$ (i.e., the $X, X'$ 
of theorem \ref{families}
are our $S, S'$). Let $ \hat{C}_3$ be a general curve in the pencil spanned by
$C_1, C_2$, and consider the pencil $ C (t) = t C_3 + ( 1-t)  \hat{C}_3$.
For each value of $t$, $C_1, C_2 , C (t)$ meet in $p_1,  .. , p_n$, 
while for $t=0$
they meet in $ 2 ab$ points, again with distinct tangents by our 
generality assumption.
We omit the other finitely many $t$'s for which the intersection 
points are more
than $n$, or the tangents are not distinct.
After blowing up $p_1,  .. , p_n$ and taking the corresponding bidouble covers,
we obtain a family $S_t$ with $S_1 = S$, and such that $S_0$ has 
exactly $2 ab - n : = h$
singular points, quadruple of the type considered in example \ref{Riem}.

Similarly, we have a family $ S'_t$, and
we must find an equisingular family $Z_u, u \in U $, containing $S_0 $ 
and $S'_0$.

Let $\PP$ be the linear system $\PP ( H^0 ( Q , \hol_Q (a,b))$, and 
consider a general curve
in the Grassmannian $ Gr ( 1, \PP)$, giving a one dimensional family 
$ C_1[w], C_2 [w], w \in W$,
of pairs of points of  $\PP$ such that $ C_1[w]$ and $ C_2 [w]$ 
intersect transversally
in $ 2 ab$ points of $Q$.

Now, the covering of $W$ given  by $$ \{(w, p_1(w),  \dots p_n (w)) 
|p_1 (w), \dots p_n (w)
\in C_1[w] \cap  C_2 [w], \  p_i(w)
\neq  p_j(w) {\rm for} \ i \neq j  \}$$
is irreducible. This  is a consequence of the General Position 
Theorem (see \cite{acgh}, page 112)
stating that if $C$ is a smooth projective curve, then for each integer $n$ the
subset $C^n_{dep} \subset C^n$, $$  C^n_{dep} : = \{ (p_1, \dots p_n )|  p_i
\neq  p_j{\rm for} \ i \neq j, \  p_1, \dots p_n \ {\rm are \ linearly \
dependent} \}$$
is smooth and irreducible.

We obtain then a one dimensional family with irreducible basis $U$ of 
rational surfaces
$ Y (u)$, obtained blowing up $Q$ in the $n$ points $p_1 (w(u)), \dots 
p_n (w(u))$,
and a corresponding family $Z_u$ of singular bidouble covers of $Y(u)$, each
with $ 2ab - n$ singularities of the same type described above.

We have then the situation of theorem \ref{families}, whence it follows that
$S, S'$, endowed with their canonical symplectic structures, are 
symplectomorphic.

\qed

The same argument , mutatis mutandis, shows (compare \cite{cat02}, 
\cite{cat06})

\begin{teo}
  Manetti surfaces of the same type
(same integers $ a,b,n,r= 2n + 7$, same divisor classes $[ D_{\sigma}]$)
are canonically symplectomorphic.
\end{teo}

Manetti indeed gave the following counterexample to the Def= Diff question:

\begin{teo} {\bf (Manetti)}
For each integer $h > 0$ there exists a surface of general type $ S$
with first Betti number $ b_1 (S) = 0$,
such that the subset of the moduli space corresponding to surfaces which
are orientedly diffeomorphic to $S$ contains at least $h$ connected components.
\end{teo}

\begin{oss}
Manetti proved the diffeomorphism of the surfaces which are here 
called Manetti surfaces
using some results of Bonahon (\cite{bonahon}) on the diffeotopies of 
lens spaces.

We have given a more direct proof also because of the application to
canonical symplectomorphism.

\end{oss}

\begin{cor}
For each integer $h > 0$ there exist  surfaces of general type $ S_1, 
\dots S_h$
with first Betti number $ b_1 (S_j) = 0$, socalled {\bf Manetti surfaces},
which are canonically symplectomorphic, but which belong to $h$ distinct
connected components of the moduli space of surfaces diffeomorphic to $S_1$.
\end{cor}

\medskip

In spite of the fact that we begin to have quite a variety of examples
and counterexamples, there are quite a few interesting open questions,
the first one concerns the existence of  simply connected surfaces which
are canonically symplectomorphic, but not deformation equivalent:

\begin{question} Are the diffeomorphic $(a,b,c)$-surfaces of theorem
\ref{diffabc}, endowed with their canonical symplectic structure,
indeed symplectomorphic?
\end{question}

\begin{oss}
  A possible way of showing that the answer to the question above is
    yes (and therefore exhibiting symplectomorphic simply connected
    surfaces which are not deformation equivalent) goes through
  the analysis of the braid monodromy of the branch curve of the
``perturbed'' quadruple
covering of $\PP^1 \times \PP^1$ (the  composition of the perturbed 
covering with
the first projection $\PP^1 \times \PP^1 \to  \PP^1$ yields the 
Lefschetz fibration).
One would like to see whether  the involution
$\iota$ on
$\PP^1$, $\iota(y) = -y$ can be written as the product of braids which
show up in the factorization.

This approach turns out to be more difficult than the
corresponding analysis which has been made in the mapping class
group, because the braid monodromy contains very many
'tangency' factors which do not come from local contributions
to the regeneration of
the branch curve from the union of the curves $f=0, g=0$ counted twice.

\end{oss}

\begin{question} Are there (minimal) surfaces of general type which 
are orientedly
  diffeomorphic through a diffeomorphism carrying the canonical class 
to the canonical
class, but, endowed with their canonical symplectic structure,
are not
canonically   symplectomorphic?

Are there such examples in the simply connected case ?
\end{question}

The difficult question is then: how to show that diffeomorphic surfaces
(diffeomorphic through a diffeomorphism carrying the canonical class to 
the canonical
class) are not symplectomorphic ?

We shall briefly comment on this in the next section, referring the reader
to the  other Lecture Notes in this volume (for instance, the one by 
Auroux and Smith)
  for more details.

\subsection{Braid monodromy and Chisini' problem.}

Let $B \subset \PP^2_{\C}$ be a plane algebraic curve of degree $d$,
and let $P$ be a general point not on $B$. Then the pencil of lines
$L_t$ passing through  $P$ determines a one parameter family of
$d$-uples of points of $\C \cong L_t \backslash \{P\}$, namely, $ L_t \cap B$.

Blowing up the point $P$ we get the projection
$\F_1 \to \PP^1$, whence the  braid at infinity is a full rotation,
corresponding to  the generator of the (infinite cyclic)
center of the braid group $\mathcal{B}_d$, $$(\Delta^2_d) : =
(\sigma_{d-1}\sigma_{d-2} \dots \sigma_{1})^d.$$

Therefore one gets a factorization of $\Delta^2_d $ in the braid
group $\mathcal{B}_d$, and the equivalence class of the factorization
does neither depend on the point $P$ (if $P$ is chosen to be general), nor does
it depend on $B$, if $B$ varies in an equisingular family of curves.

Chisini was mainly interested in the case of {\em cuspidal} curves
(compare \cite{chis1}, \cite{chis2}), mainly because these are the
branch curves of a generic projection $f : S \ra \PP^2_{\C}$, for any
smooth projective surface $S \subset \PP^r$.

  More precisely, a {\em generic
projection} $f : S \ra \PP^2_{\C}$ is a covering whose branch
curve $B$ has only nodes and cusps as singularities, and
moreover is such that the local monodromy around
a smooth point of the branch curve is a transposition.

Maps with those properties are called {\em generic coverings}:
for these the local monodromies are only  $\Z/2 = \mathfrak{S}_2$
(at the smooth points of the branch curve $B$),
$\mathfrak{S}_3$ at the cusps, and $\Z/2 \times \Z/2$
at the nodes.

In such a case we have a {\em cuspidal} factorization, i.e.,  all factors
are powers of a half twist, with respective exponents $1,2,3$.

Chisini posed the following daring

\begin{conj} (Chisini's conjecture.)

Given two generic coverings $f : S \ra \PP^2_{\CC}$, $f' : S' \ra
\PP^2_{\CC}$, both  of degree at least $5$, assume that
they have the same branch curve
$B$. Is it then true that $f$ and $f'$ are equivalent?
\end{conj}

Observe that the condition on the degree is necessary, since a
counterexample for $d \leq 4$ is furnished by
the dual curve $B$ of a smooth plane cubic (as
already known to Chisini, cf. \cite{chis1}). Chisini in fact observed that
there are two generic coverings, of respective degrees $3$ and $4$,
and with the given branch curve. Combinatorially, we have a triple
of transpositions corresponding in one case to the sides of a triangle ($d=3$,
and the monodromy permutes the vertices of the triangle), and in the other case
to the three medians of the triangle ( $d=4$,
and the monodromy permutes the vertices of the triangle plus the barycentre).

While establishing a very weak form of the conjecture (\cite{catchis}).
I remarked that the dual curve $B$ of a smooth plane cubic is also the
branch curve for three nonequivalent generic covers of the plane
from the Veronese surface (they are distinct since they determine 
three distinct
divisors of 2-torsion on the cubic).

The conjecture seems now to have been  finally proven in the strongest 
form, after that it
was first proven  by Kulikov (cf. \cite{kulikov1}) under a rather 
complicated assumption,
and that shortly later Nemirovskii  \cite{nemiroski} noticed (just by using the
Miyaoka-Yau inequality) that Kulikov's  assumption was implied by the
simpler assumption $d \geq 12$.

Kulikov proved now  (\cite{kulikov2}) the following

\begin{teo}
{\bf (Kulikov)}
Two generic projections with the same cuspidal branch
curve $B$  are isomorphic unless if the projection $ p : S \to \PP^2$
of one of them is just
a linear projection of the Veronese surface.

\end{teo}

The above statement concerns a fundamental property of the fundamental
group of the complement $ \PP^2 \setminus B$, namely to admit only one
conjugacy class of surjections onto a symmetric group $\SSS_n$, satisfying the 
properties of a generic covering.

In turn, the fundamental group $\pi_1 (  \PP^2 \setminus B )$ is completely
determined by the braid monodromy of $B$, i.e., the above equivalence class
(modulo Hurwitz equivalence and simultaneous conjugation)
of the above factorization of $\Delta^2_d $. So, a classical question was:
which are the braid monodromies of cuspidal curves?

Chisini found some necessary conditions, and proposed some argument in order
to show the sufficiency of these conditions, which can be reformulated as

{\bf Chisini' s problem:} (cf.\cite{chis2}).

Given a cuspidal factorization, which is regenerable to the
   factorization of a smooth plane curve, is there a cuspidal curve
   which induces the given factorization?

{\em Regenerable} means that there is a factorization (in the
equivalence class) such that, after replacing each factor $\sigma^i$
($i=2,3$) by the $i$ corresponding factors (e.g. , $\sigma^3$ is
replaced by $\sigma \circ \sigma \circ \sigma$) one obtains the
factorization belonging to a non singular plane curve.

\bigskip
A negative answer to the problem of Chisini was given by
  B. Moishezon in \cite{moi}.

\begin{oss}
1) Moishezon proves that there exist infinitely many non equivalent
cuspidal factorizations observing that
$\pi_1(\PP^2_{\C} \backslash B)$ is an invariant defined in terms
of the
factorization alone and constructing infinitely many 
non isomorphic such groups. On the other hand, the family of
cuspidal curves  of a fixed
degree is an algebraic set, hence it has a finite number of 
connected components.
These two statements together give a negative answer to the above cited
problem of Chisini.

The examples of Moishezon have been recently reinterpreted
in \cite{adk}, with a simpler
treatment, in terms of symplectic surgeries.

\end{oss}

Now, as conjectured by Moishezon, a cuspidal factorization together with
a generic monodromy with values in $\mathfrak{S}_n$ induces a covering
$ M \ra \PP^2_{\CC}$, where the fourmanifold $M$ has a unique symplectic
structure (up to symplectomorphism) with
class equal to the pull back of the Fubini Study form on $\PP^2$
(see  for instance \cite{a-k}).

What is more interesting (and much more difficult) is however the converse.

Extending Donaldson's techniques (for proving the existence of
symplectic Lefschetz fibrations) Auroux and Katzarkov
(\cite{a-k}) proved that each symplectic $4$-manifold is in a
natural way 'asymptotically' realized by such a generic covering.

They show that, given a symplectic fourmanifold $ (M, \omega)$ with
$ [\omega ] \in H^2 (M, \Z)$, there exists a multiple $m$ of a line bundle
$L$ with $ c_1(L) = [\omega]$ and three general sections
$s_0, s_1, s_2 $ of $ L^{\otimes m}$, which are $\epsilon$-holomorphic
with many of their derivatives (that a section $s$ is $\epsilon$-holomorphic
  means, very roughly speaking,
that once one has chosen a compatible almost complex structure,
$ | \bar \partial  s | < \epsilon  \ | \partial s| $)
yielding a finite covering of the plane  $\PP^2$ which
is generic and with branch curve a symplectic subvariety whose
singularities are only nodes and cusps.

The only price they have to pay is to allow also negative nodes, i.e.,
nodes which in local holomorphic coordinates are defined by the equation
$$ ( y -\bar x) ( y + \bar x) = 0. $$

The corresponding factorization in the braid group  contains then
only factors which are conjugates of $ \sigma_1^j$, with $ j = -2, 1, 2,3$.

Moreover, the factorization is not unique, because it may happen that
two consecutive nodes, one positive and one negative, may disappear,
and the corresponding two factors disappear from the factorization.
In particular,  $\pi_1(\PP^2_{\C} \backslash B)$ is no longer an invariant
and the authors
propose to use an appropriate quotient of $\pi_1(\PP^2_{\C} \backslash B)$ in
order to produce invariants of symplectic structures.

It seems however that, in the computations done up to now, even the  groups
$\pi_1(\PP^2_{\C} \backslash B)$ allow only to detect
homology invariants of the projected fourmanifold (\cite{a-d-k-y}).

Let us now return to the world of surfaces of general type.

Suppose we have a surface $S$ of general type and a pluricanonical
embedding $\psi_m \colon X \to \PP^N$ of the canonical model $X$ of $S$.
Then a generic linear projection of
the pluricanonical image to
$\PP^2_{\C}$ yields, if $ S \cong X$, a generic covering $S \to \PP^2_{\C}$
(else the singularities of $X$ create further singularities for the
branch curve $B$ and other local coverings).

By the positive solution of Chisini's conjecture, the branch curve 
$B$ determines
the surface $S$ uniquely (up to isomorphism). We get moreover the equivalence
class of the braid monodromy factorization, and this does not change
if $S$ varies in a connected family of surfaces with $K_S$ ample (i.e., the
surfaces equal their canonical models).

Motivated by this observation of Moishezon,   Moishezon and Teicher in a
series of technically difficult papers ( see e.g. \cite{m-t})
tried to calculate fundamental groups of complements
  $\pi_1(\PP^2_{\C} \backslash B)$, with the intention of distinguishing
connected components of the moduli spaces of surfaces of general type.

Indeed, it is clear that these groups are invariants of the
connected components of the open set of the moduli space
corresponding to surfaces with ample canonical divisor $K_S$.
Whether one could even distinguish connected components
of moduli spaces would in my opinion deserve a detailed argument,
in view of the fact that several irreducible components consist
of surfaces whose canonical divisor is not ample (see for instance 
\cite{cat5} for
several series of examples).

But it may be that the information carried by $\pi_1(\PP^2_{\C} \backslash B)$
be too scanty, so one could look at further combinatorial invariants,
rather than the class of the braid monodromy factorization for $B$.

In fact a  generic linear projection of
the pluricanonical image to  $\PP^3_{\C}$ gives a surface $\Sigma$
with a double curve $\Gamma'$. Now, projecting further to $\PP^2_{\\C}$
  we do not only get the branch curve $B$, but also a plane curve
$\Gamma$, image of $\Gamma'$.

Even if Chisini's conjecture tells us that from the holomorphic point
of view $B$ determines the surface $S$ and therefore the curve
$\Gamma$, it does not follow that the fundamental group
$\pi_1(\PP^2_{\CC} \backslash B)$ determines the group
  $\pi_1(\PP^2_{\CC} \backslash (B \cup \Gamma))$.

It would be interesting to calculate this second fundamental group,
even in special cases.

Moreover, generalizing a proposal done
by Moishezon in \cite{moi2}, one can observe that the monodromy
of the restriction of the covering $\Sigma \to \PP^2$
to  $\PP^2_{\CC} \backslash (B \cup \Gamma))$
is more refined, since it takes values in a braid group $\sB_n$ , rather than
in a symmetric group $\SSS_n$.

One could proceed similarly also for the generic projections of 
symplectic fourmanifolds.

But in the symplectic case one does not have the advantage of 
knowing a priori
an explicit  number  $m \leq 5$ such that $\psi_m$ is a pluricanonical
embedding for the general surface $S$ in the moduli space.

\bigskip

{\bf Acknowledgements.}
I would like to thank the  I.H.P. for its hospitality in november 
2005, and  the I.H.E.S. for
its hospitality in march 2006, which allowed me to  write down some
  parts  of these
  Lecture Notes.

I am grateful to Fabio Tonoli for his invaluable help with the figures,
and to Ingrid Bauer, Michael L\"onne, Roberto Pignatelli and S\"onke Rollenske
for comments on a preliminary version. 

\bigskip

\newpage

\begin{footnotesize}
\noindent

\end{footnotesize}

\vfill

\noindent
{\bf Author's address:}

\bigskip

\noindent
Prof. Fabrizio Catanese\\
Lehrstuhl Mathematik VIII\\
Universit\"at Bayreuth, NWII\\
  D-95440 Bayreuth, Germany

e-mail: Fabrizio.Catanese@uni-bayreuth.de

\end{document}